# Publication of my Master Thesis
## *"Using Selmer Groups to compute Mordell-Weil Groups of Elliptic Curves"*

Anika Behrens

December 19, 2018

This document is a publication of my master thesis entitled *"Using Selmer Groups to compute Mordell-Weil Groups of Elliptic Curves"* of 2016. Please note that this thesis aimed to be a *summary and reproduction* of the state of the art of that time. Content and results (in particular the mathematical results given in Part II) refer to the theory presented in the book "The Arithmetic of Elliptic Curves, 2nd Edition" by Joseph H. Silverman, especially Chapter X "Computing the Mordell-Weil Group".

This publication contains the **master thesis** (80 pages) and a **correction document** (2 pages) attached at the end.

It was evaluated with grade 1.0 (excellent).

# Using Selmer Groups to compute Mordell-Weil Groups of Elliptic Curves

**Master Thesis**

University of Bremen

Anika Behrens
Matrikel Nr. 2208335

September 2016

*First Referee:*
   Prof. Dr. Jens GAMST

*Second Referee:*
   Prof. Dr. Eberhard OELJEKLAUS

# Contents

# Acknowledgements

I would like to express deep gratitude to my advisor Prof. Dr. Jens Gamst. He guided me through my reading course and master module.

I would also like to thank Prof. Dr. Eberhard Oeljeklaus for agreeing to be the second reader of my thesis.

Special thanks to Prof. Joseph H. Silverman for replying so quickly to my questions.

Further, I want to thank my family and best friends for any support during my studies.

# Part I
# General Introduction and Preparation

In this thesis, I want to show how Selmer groups can be used to determine the Mordell-Weil group of elliptic curves over a number field $K$. In *Part I* of this thesis I give an introduction to the main definitions and properties. The Mordell-Weil group $E(K)$ of an elliptic curve $E/K$ is defined to be the set of $K$-rational points of E, i.e. the points with the coordinates contained in $K$. This set is a subgroup of the group $E(\overline{K})$ of points of $E$. Because it is always possible to compute the Mordell-Weil group $E(K)$ from the weak Mordell-Weil group $E(K)/mE(K)$ (see e.g. [SIL1] Exercise 8.18), I want to concentrate on the question of how to compute the weak Mordell-Weil group. Unfortunately, one has not found a procedure that is guaranteed to give generators for $E(K)/mE(K)$ by a finite amount of computation. Let K be a number field. The Mordell-Weil Theorem states that

$$E(K) = E(K)_{\text{tors}} \times \mathbb{Z}^r,$$

where $r$ is the *rank* of $E$ and $E(K)_{\text{tors}}$ is the *torsion subgroup*, i.e. the group of points of finite order in $E(K)$. The group $E(K)_{\text{tors}}$ is finite and well understood. In particular, for $K = \mathbb{Q}$ Mazur's theorem tells us that $E(\mathbb{Q})_{\text{tors}}$ is isomorphic to one of fifteen possible groups. So, one tries to find a way to determine the rank $r$ of $E$, which is the major problem. The procedure described in this thesis shows how to transfer the computation of the weak Mordell-Weil group $E(K)/mE(K)$ to the existence or non-existence of a rational point on certain curves, called *homogeneous spaces*. If one can find some completion $K_v$ of $K$ such that the homogeneous space has no points in $K_v$, then it follows that it has no points in $K$.

Suppose that we have another elliptic curve $E'K$ isogeneous to $E'/K$ via an isogeny $\phi : E \to E'$. The Selmer group $S^\phi(E/K)$ and Shafarevich-Tate group Ш$(E/K)$ are defined to be

$$\begin{aligned} S^\phi(E/K) &:= \ker\{H^1(G_{\overline{K}/K}, E[\phi]) \to \prod_{v \in M_K} WC(E/K_v)\} \\ \text{Ш}(E/K) &:= \ker\{WC(E/K) \to \prod_{v \in M_K} WC(E/K_v)\}, \end{aligned}$$

where $M_K$ is a complete set of inequivalent absolute values on $K$.

Then, there is an exact sequence

$$0 \to E'(K)/\phi(E(K)) \to S^{(\phi)}(E/K) \to \text{Ш}(E/K)[\phi] \to 0.$$

I will prove that $S^{(\phi)}(E/K)$ is finite and give a procedure to compute the Selmer group if $\phi$ is an isogeny of degree 2. Having computed the Selmer groups for $\phi : E \to E'$ and its dual isogeny $\hat{\phi} : E' \to E$, one tries to find generators for $E'(K)/\phi(E(K))$ and $E(K)/\hat{\phi}(E'(K))$, such that the exact sequence

$$0 \to \frac{E'(K)[\hat{\phi}]}{\phi(E(K)[m])} \to \frac{E'(K)}{\phi(E(\mathbb{Q}))} \xrightarrow{\hat{\phi}} \frac{E(K)}{mE(K)} \to \frac{E(K)}{\hat{\phi}(E'(K))} \to 0$$

gives generators for $E(K)/mE(K)$.

Finally, under the assumption that the Shafarevich-Tate group is finite, the rank of Elliptic curves over $\mathbb{Q}$ with j-invariant 1728 is fully determined in certain cases, which is done in *Part II* of this thesis. An outlook on further research is given in *Part III*.

# 1 Basic Definitions and Properties

In the following, let $K$ be a *perfect field*, let $\overline{K}$ be a fixed *algebraic closure* of $K$, and let char$K$ denote the *characteristic* of the field $K$.

I want to put *elliptic curves* into a more abstract context of *algebraic varieties* and give basic definitions and properties that are necessary for understanding the concrete topic of this thesis. This also includes an introduction to the main definitions used in Galois and cohomology theory. So, I will not prove any properties here, but give a brief summary of the most important objects and relations given in [SIL1] Appendix B, Chapter I to III, and VIII.

## 1.1 Galois and Cohomology Theory

### 1.1.1 Cohomology of finite Groups

Let $G$ be a finite group and let $M$ be a $G$-Module, i.e. an abelian group s.t. $G$ acts on $M$, denoted by $m \mapsto m^\sigma$, for $\sigma \in G$ and $m \in M$.
The *0th cohomology group* is the submodule of $M$ consisting of all $G$-invariant elements:

$$H^0(G, M) = M^G = \{m \in M : m^\sigma = m \text{ for all } \sigma \in G\}.$$

The *group of 1-cochains (from G to M)* is the set of maps

$$C^1(G, M) = \{\xi : G \to M\}.$$

The *group of 1-cocycles (from G to M)* is defined to be

$$Z^1(G, M) = \{\xi \in C^1(G, M) : \xi_{\sigma\tau} = \xi_\sigma^\tau + \xi_\tau \text{ for all } \sigma, \tau \in G\}.$$

The *group of 1-coboundaries (from G to M)* is given by

$$B^1(G, M) = \{\xi \in C^1(G, M) : \text{ there exists } m \in M \text{ s.t.} \xi_\sigma = m^\sigma - m \text{ for all } \sigma \in G\}.$$

The *1st cohomology group* of the $G$-module $M$ is the quotient group

$$H^1(G, M) = \frac{Z^1(G, M)}{B^1(G, M)}.$$

### 1.1.2 Inverse Limit and profinite Groups

Let $(I, \leq)$ be a partially ordered set, let $(G_i)_i \in I$ be a familily of groups, and let $(f_{ij})_{i,j \in I, i \leq j}$ be a family of homomorphisms $f_{ij} : G_j \to G_i$ s.t. $f_{ii} = id_{G_i}$ and $f_{ij} \circ f_{jk} = f_{ik}$ for all $i, j, k \in I$ with $i \leq j \leq k$. Then we call $((G_i)_{i \in I}, (f_{ij})_{i,j \in I, i \leq j})$ an *inverse system* over $I$ and the inverse limit is defined to be

$$\varprojlim_{i \in I} A_i = \left\{ \vec{a} \in \prod_{i \in I} A_i \ \Big| \ a_i = f_{ij}(a_j) \text{ for all } i, j \in I \text{ with } i \leq j \right\},$$

where $\prod$ denotes the direct product of groups.

A group is called *profinite*, if it is given by an inverse limit of finite groups.

### 1.1.3 Galois Group and Krull-Topology

Let $L/K$ be a finite *Galois extension*, i.e. a *normal* and *separable* field extension.

Then the Galois group $G_{L/K}$ is defined to be

$$G_{L/K} := \left\{ \sigma \in \text{Isom}(L) \Big| k^\sigma = k \text{ for all } k \in K \right\},$$

where $\text{Isom}(L)$ denotes the automorphism group of $L$, i.e. the group of bijective homomorphisms from $L$ to $L$.

Further, we define the profinite group

$$G_{\overline{K}/K} := \varprojlim_{\substack{L/K \\ \text{finite} \\ \text{Galois}}} G_{L/K}.$$

This group can be equipped with a topology, called *Krull-Topology*, where a *neighbourhood base* of *id* is given by

$$\{G_{L/M} \mid K \subseteq M \subseteq L \text{ and } [M : K] < \infty\}.$$

### 1.1.4 Galois Cohomology

In the case of a profinite Galois group, we need additional restrictions on the modules and cocycles: For a $G_{\overline{K}/K}$-Module, the action of $G$ on $M$ needs to be *continuous* for the Krull-Topology on $G_{\overline{K}/K}$ and for the discrete topology on $M$. We recall that the discrete topology consists of all subsets as open sets. *Continuous* means that for every open set the inverse image is an open subset. If $M$ is equipped with the discrete topology, coboundaries $B^1(G_{\overline{K}/K}, M)$ are automatically continuous. But for cocycles we need to define $Z^1_{\text{cont}}(G_{\overline{K}/K}, M)$ to be the set of cocycles *continuous* for the Krull-Topology on $G_{\overline{K}/K}$ and for the discrete topology on $M$. Then

$$\begin{aligned} H^0(G_{\overline{K}/K}, M) &:= \{m \in M : m^\sigma = m \text{ for all } \sigma \in G_{\overline{K}/K}\} \text{ and} \\ H^1(G_{\overline{K}/K}, M) &:= \frac{Z^1_{\text{cont}}(G_{\overline{K}/K}, M)}{B^1(G_{\overline{K}/K}, M)}. \end{aligned}$$

### 1.1.5 Exact Sequences

A sequence $G_1 \xrightarrow{f_1} G_2 \xrightarrow{f_2} \ldots \xrightarrow{f_{n-1}} G_n$ of groups and homomorphisms is called *exact*, if $im(f_{i-1}) = ker(f_i)$ for all $1 < i < n$.

Let $A, B,$ and $C$ denote groups and let $G$ act on $A$, $B$ and $C$. Let

$$0 \to A \xrightarrow{\phi} B \xrightarrow{\psi} C \to 0$$

be a *short exact sequence*, i.e. $\phi : A \to B$ and $\psi : B \to C$ are homomorphisms such that $im(\phi) = ker(\psi)$, $\phi$ is injective, and $\psi$ is surjective. Then there is a *long exact sequence* of cohomology groups

$$\begin{array}{ccccccc}
0 & \longrightarrow & H^0(G,A) & \longrightarrow & H^0(G,B) & \longrightarrow & H^0(G,C) \\
& & & & & & \downarrow \delta \\
& & H^1(G,A) & \longrightarrow & H^1(G,B) & \longrightarrow & H^1(G,C)
\end{array},$$

where the connecting homomorphism $\delta : H^0(G,C) \to H^1(G,A)$ can be constructed as follows: Let $c \in H^0(G,C)$, i.e. $c \in C$, such that $c^\sigma = c$ for all $\sigma \in G$. Since $\psi : B \to C$ is surjective, one can choose $b \in B$, such that $\psi(b) = c$. We can define a cochain $\xi \in C^1(G,B)$, i.e. $\xi : G \to B$, by $\xi_\sigma = b^\sigma - b$. Because of $c = c_\sigma$, we have $\psi(\xi_\sigma) = 0$, and $\xi \in Z^1(G,A)$. Finally, $\delta(c)$ is defined to be the cohomology class of $\xi$.

### 1.1.6 Nonabelian Cohomology of finite Groups

See also [SIL1] Appendix B Section 3.

Let $G$ be a finite group and let $M$ be a group s.t. $G$ acts on $M$. In order to indicate that $M$ does not need to be abelian here, we write the group law multiplicatively, now.
The *0th cohomology group* is the submodule of $M$ consisting of all $G$-invariant elements:

$$H^0(G,M) = M^G = \{m \in M : m^\sigma = m \text{ for all } \sigma \in G\}.$$

The *set of 1-cocycles (from G to M)* is defined to be

$$Z^1(G,M) = \{\xi \in C^1(G,M) : \xi_{\sigma\tau} = (\xi_\sigma)^\tau + \xi_\tau \text{ for all } \sigma, \tau \in G\},$$

which does not form a group, in general.
The *1st cohomology set* of $M$ is the quotient

$$H^1(G,M) = Z^1(G,M)/\sim,$$

where $\sim$ is the following equivalence relation on the set of 1-cocycles

$$\xi \sim \zeta \quad :\Leftrightarrow \quad \text{there exists } m \in M \text{ s.t. } m^\sigma \xi_\sigma = \zeta_\sigma m \text{ for all } \sigma \in G.$$

Such cocycles $\xi$ and $\zeta$ are called *cohomologous*.

### 1.1.7 Nonabelian profinite Galois Cohomology

Again, in the following we do not suppose that $M$ is abelian. But we assume that $G_{\overline{K}/K}$ acts *continuously* on $M$, i.e. for each $m \in M$ the *stabilizer*

$$stab(m) := \{\sigma \in G_{\overline{K}/K} : m^\sigma = m\}$$

is a subgroup of finite index in $G_{\overline{K}/K}$.
Then $H^0(G_{\overline{K}/K}, M)$ and $H^1(G_{\overline{K}/K}, M)$ are defined as in the abelian case, respecting the profinite topology on $G_{\overline{K}/K}$ and discrete topology on $M$.

## 1.2 Algebraic Varieties in Affine $n$-Space

(see also [SIL1], Chapter I, Section 1)

### 1.2.1 Affine $n$-Space

*Affine n-space (over K)* is given by

$$\mathbb{A}^n = \mathbb{A}^n(\overline{K}) := \{P = (x_1, ... x_n) : \text{ all } x_i \in \overline{K}\}$$

The *set of K-rational points* of $\mathbb{A}^n$ is defined to be

$$\mathbb{A}^n(K) := \{P = (x_1, ... x_n) : \text{ all } x_i \in K\}.$$

### 1.2.2 Affine Algebraic Sets

An *affine algebraic set* is defined to be a set of the form

$$\{P \in \mathbb{A}^n : f(P) = 0 \text{ for all } f \in I\},$$

where $I \subset \overline{K}[X]$ is an ideal in the polynomial ring in variables $X = (X_1, ..., X_n)$ over $\overline{K}$.

The *ideal* of an affine algebraic set $V$ is given by

$$I(V) = \{f \in \overline{K}[X] : f(P) = 0 \text{ for all } P \in V\}.$$

An algebraic set $V$ is *defined over K*, denoted by $V/K$, if $I(V)$ can be generated by elements in $K[X]$.

If $V/K$, then *set of K-rational points* of $V$ is defined to be $V(K) = V \cap \mathbb{A}^n(K)$.

### 1.2.3 Affine Varieties

An *affine variety* is an affine algebraic set $V$ s.t. $I(V)$ is a prime ideal in $\overline{K}[X]$.

If $V/K$ is an affine variety, then the *affine coordinate ring* of $V$ is defined to be $K[V] = \frac{K[X]}{I(V) \cap K[X]}$ and the *(rational) function field* of $V/K$, denoted by $K(V)$, is the quotient field of $K[V]$.

Analogously, for any field extension $L/K$ one can define $L[V] := \frac{L[X]}{I(V) \cap L[X]}$ and $L(V)$ is the quotient field of $L[V]$.

The *dimension* of an algebraic variety, denoted by $dim(V)$, is the transcendence degree of the field extension $\overline{K}(V)/\overline{K}$. For more details see [SIL1] I.1, p.4, Ex.1.12.

An algebraic variety is called *nonsingular* or *smooth* at $P \in V$ if for any set of generators $(f_1, ..., f_m)$ for $I(V)$ the rank of the $m \times n$ *Jacobi*-matrix is

$$rank\left(\left(\frac{\partial f_i}{\partial X_j}(P)\right)\right) = n - dimV.$$

An algebraic variety is called *nonsingular* or *smooth* if it is *nonsingular* at every point.

The *local ring of V at P* is defined to be

$$\overline{K}[V]_P := \{F \in \overline{K}(V) : \exists f, g \in \overline{K}[V], g(P) \neq 0, F(P) = f(P)/g(P)\}.$$

The functions in $\overline{K}[V]_P$ are called *regular* or *defined* at $P$.

If $V$ is a nonsingular curve, the local ring is a *discrete valuation ring*, i.e. it has a unique non-zero maximal ideal and every ideal is principal, i.e., can be generated by a single element.

A generator of the maximal ideal is called an *uniformizing element*. It is unique up to multiplication by a unit. Every nonzero element then is of the form $x = \alpha t^k$ with a unit $\alpha$ and a power $t^k$ of an uniformizing element. The number $k \in \mathbb{N}$ is called the *order* of $x$.

## 1.3 Algebraic Varieties in Projective $n$-Space

(see also [SIL1], Chapter I, Section 2)

### 1.3.1 Projective $n$-Space

*Projective n-space (over K)* is given by the quotient

$$\mathbb{P}^n = \mathbb{P}^n(K) = \left(\mathbb{A}^{n+1}(K) \setminus \{0\}\right) / \sim,$$

where $\sim$ is the equivalence relation on $\mathbb{A}^{n+1} \setminus \{0\}$ defined by

$(x_0, ..., x_n) \sim (y_0, ..., y_n) \quad :\Leftrightarrow \quad$ there exists a $\lambda \in \overline{K}^*$ such that $x_i = \lambda y_i$ for all $i$.

For any point $(x_0, ..., x_n) \in \mathbb{A}^{n+1} \setminus \{0\}$ we denote its equivalence class by $[x_0, ..., x_n] \in \mathbb{P}^n$.

The *set of K-rational points* of $\mathbb{A}^n$ is defined to be

$$\mathbb{P}^n(K) := \{[x_0, ...x_n] \in \mathbb{P}^n : \text{ all } x_i \in K\}.$$

Let $\overline{K}[X] = \overline{K}[X_0, ..., X_n]$ be the polynomial ring in $n+1$ variables. A polynomial $f \in \overline{K}[X]$ is called *homogeneous* of degree $d \in \mathbb{N}_0$, if

$$f(\lambda X_0, ..., X_n) = \lambda^d f(X_0, ..., X_n) \text{ for all } \lambda \in \overline{K}.$$

In particular, $f = 0$ is homogeneous of degree $d$ for any $d \in \mathbb{N}_0$.
Monomials $X^\alpha = X_0^{d_0} \ldots X_n^{d_n}$ with $|d| = d_0 + \ldots + d_n$ are homogeneous of degree $d$.

### 1.3.2 Projective Algebraic Sets

A *projective algebraic set* is defined to be a set of the form

$$\{P \in \mathbb{P}^n : f(P) = 0 \text{ for all homogeneous } f \in I\},$$

where $I \subset \overline{K}[X]$ is an ideal in the polynomial ring with coefficients in $\overline{K}$ a $n+1 \in \mathbb{N}$ variables $X = X_0, ..., X_n$.

The *homogeneous ideal* of a projective algebraic set $V$ is given by

$$I(V) = \{f \in \overline{K}[X] : f \text{ is homogeneous and } f(P) = 0 \text{ for all } P \in V\}.$$

A projective algebraic set $V$ is *defined over $K$*, denoted by $V/K$, if $I(V)$ can be generated by *homogeneous* polynomials in $K[X]$.

If $V/K$ is a projective algebraic set, then the *set of $K$-rational points* of $V$ is defined to be $V(K) = V \cap \mathbb{P}^n(K)$.

For each $i \in \{0, \ldots, n\}$ there is a natural embedding

$$\phi_i : \mathbb{A}^n \to \mathbb{P}^n$$
$$(y_1, \ldots, y_n) \mapsto [y_1, \ldots, y_{i-1}, 1, y_i, \ldots y_n].$$

We define the hyperplane $H_i := \{[x_0, \ldots, x_n] \in \mathbb{P}^n : x_i = 0\}$. There is a natural bijection

$$\phi_i^{-1} : \mathbb{P}^n \setminus H_i \to \mathbb{A}^n$$
$$[x_0, \ldots, x_n] \mapsto \left(\frac{x_0}{x_i}, \ldots, \frac{x_{i-1}}{x_i}, \frac{x_{i+1}}{x_i}, \ldots, \frac{x_n}{x_i}\right).$$

For $g(X_0, \ldots, X_n) \in \overline{K}[X]$ the polynomial

$$g_*(Y_1, \ldots, Y_n) = g(Y_1, \ldots, Y_{i-1}, 1, Y_{i+1}, \ldots, Y_n)$$

is called the *dehomogenization of $g$ with respect to $X_i$*.

For $f(Y_1, \ldots, Y_n) \in \overline{K}[Y]$ the polynomial

$$f^*(X_0, \ldots, X_n) = X_i^{\deg(f)} f\left(\frac{X_0}{X_i}, \ldots, \frac{X_{i-1}}{X_i}, \frac{X_{i+1}}{X_i}, \ldots, \frac{X_n}{X_i}\right)$$

is called the *homogenization of $f$ with respect to $X_i$*.

For a *projective* algebraic set $V$ we have

$$V \cap \mathbb{A}^n := \phi_i^{-1}(V \setminus H_i) = \{P \in \mathbb{A}^n : g_*(P) = 0 \text{ for all } g \in I(V)\},$$

which is an *affine* algebraic set.

For an *affine* algebraic set $V$ the *projective* algebraic set

$$\overline{V} := \{P \in \mathbb{P}^n : f^*(P) = 0 \text{ for all } f \in I(V)\}$$

is called the *projective closure* of $V$.

### 1.3.3 Projective Algebraic Varieties

A *projective algebraic variety* is a projective algebraic set $V$ s.t. $I(V)$ is a prime ideal in $\overline{K}[X]$.

*Dimension*, *coordinate ring*, *function field*, *(non-)singularity*, the *local ring*, and *regular functions* at a point $P$ of a projective variety can be defined analogously as in the case of affine varieties.

## 1.4 Galois Action on Algebraic Varieties

The Galois group $G_{\overline{K}/K}$ acts on

- $\overline{K}$ by definition
- $\mathbb{A}^n$ by acting on the affine coordinates
- $\mathbb{P}^n$ by acting on the homogeneous coordinates

- $\overline{K}[X]$ by acting on the coefficients.

One can identify

$$\mathbb{A}^n(K) = \{P \in \mathbb{A}^n : P^\sigma = P \text{ for all } \sigma \in G_{\overline{K}/K}\} \text{ and}$$
$$\mathbb{P}^n(K) = \{P \in \mathbb{P}^n : P^\sigma = P \text{ for all } \sigma \in G_{\overline{K}/K}\}.$$

For an (affine or projective) algebraic set $V$:

$$V(K) = \{P \in V : P^\sigma = P \text{ for all } \sigma \in G_{\overline{K}/K}\}$$

And for an (affine or projective) algebraic variety $V$:

$$K[V] = \{f \in \overline{K}[V] : f^\sigma = f \text{ for all } \sigma \in G_{\overline{K}/K}\}$$
$$K(V) = \{F \in \overline{K}(V) : F^\sigma = F \text{ for all } \sigma \in G_{\overline{K}/K}\}$$

## 1.5 Kähler Differentials

Let $R$ be a commutative ring (later we will set $R = \overline{K}$), and let $A$ be an algebra over $R$ (later we will set $A = \overline{K}(C)$ for a curve $C$). A *R-linear derivation on $A$* is a map $d : A \to M$, where $M$ is an $A$-module, s.t.

i) $d$ is $R$-linear, i.e.
   (a) $d(x+y) = d(x) + d(y)$ for all $x, y \in A$
   (b) $d(rx) = rd(x)$

ii) $d(xy) = d(x)y + xd(y)$ (*Leibniz Rule*)

Since $1 \in A$, one has $d(1) = d(1^2) = d(1) + d(1)$, hence $d(1) = 0$, and $d(r) = 0$ for all $r \in R$. The $A$-module generated by the symbols $dx$ for $x \in A$ subject to the relations above is called the module of *Kähler differentials*, denoted by $\Omega_{A/R}$. It is unique up to module-isomorphism and has the following universal property: For any *R-linear derivation $d' : A \to M$* on $A$ there is a unique $A$-module homomorphism $\phi : \Omega_{A/R} \to M$ such that $d' = d \circ \phi$. Setting $R = \overline{K}$ and $A = \overline{K}(C)$ for a curve $C$, then it is a fact that the module $\Omega_{A/R}$ is a 1-dimensional $\overline{K}$-vector space, denoted by $\Omega_C$. For $A = \overline{K}[X_0, \ldots, X_n]$ the vector space $\Omega_{A/K}$ has $K$-Basis $(dX_1, \ldots, dX_n)$.

## 1.6 Absolute Values and Valuations

A function $| \ | : K \to \mathbb{R}$ is called an *absolute value* on $K$ if for all $a, b \in K$

- $|a| = 0 \Leftrightarrow a = 0$
- $|a| \geq 0$
- $|ab| = |a||b|$
- $|a + b| \leq |a| + |b|$.

An absolute value defines a metric by $d(x,y) := |y - x|$. Two absolute values are called *equivalent*, if they define the same topology. An absolute value $| \ | : K \to \mathbb{R}$ is called *non-Archimedean*, if $|a+b| \leq \max\{|a|, |b|\}$ for all $a, b \in K$, otherwise it is called *Archimedean*.

A function $v : K \to \mathbb{R} \cup \infty$ is called *valuation* on $K$ if for all $a, b \in K$

- $v(a) = \infty \Leftrightarrow a = 0$
- $v(ab) = v(a) + v(b)$
- $v(a + b) \geq \min\{v(a), v(b)\}$

The valuation $v$ is called *discrete* if there exists $s \in \mathbb{Z}$ s.t. $v(K^*) = s\mathbb{Z}$. If this holds for $s = 1$, the valuation is called *normalized*.

Let $v$ be a valuation and $q \in \mathbb{R}_{>1}$, then one can construct a corresponding absolute value by

$$\begin{aligned} | \ |_q : K &\to \mathbb{R} \\ x &\mapsto \begin{cases} q^{-v(x)} & \text{if } x \neq 0 \\ 0 & \text{if } x = 0. \end{cases} \end{aligned}$$

On the other hand, one can construct a valuation from a given absolute value $| \ |$ in the following way: Let $s \in \mathbb{R}_{>0}$, then

$$\begin{aligned} v_s : K &\to \mathbb{R} \cup \{\infty\} \\ x &\mapsto \begin{cases} -s\log_s(|x|) & \text{if } x \neq 0 \\ \infty & \text{if } x = 0 \end{cases} \end{aligned}$$

is a valuation on $K$.

Let $| \ |$ be a non-Archimedean absolute value s.t. the corresponding valuation is discrete. Then the *ring of integers* $R_{||} := \{x \in K : |x| \leq 1\}$ is a *discrete valuation ring*, and the *valuation ideal* $\mathfrak{m}_{||} = \{x \in K : |x| < 1\}$ is generated by a single *uniformizing element*. Let $v$ be the valuation corresponding to $| \ |$, then $R_v := \{x \in K : v(x) \geq 0\} = R_{||}$ and $\mathfrak{m}_v := \{x \in K : v(x) > 0\} = \mathfrak{m}_{||}$. The *residue field* is given by $k_v := R_v/\mathfrak{m}_v = R_{||}/\mathfrak{m}_{||}$.

For a projective curve $C$ and a smooth point $P$ the local ring $\overline{K}[C]_P$ is a discrete valuation ring and we define the normalized valuation on the local ring

$$\begin{aligned} \text{ord}_P : \overline{K}[C]_P &\to \{0, 1, \ldots\} \cup \{\infty\} \\ f &\mapsto \sup\{d \in \mathbb{Z} | f \in M_P^d\}, \end{aligned}$$

where $M_P$ denotes the maximal ideal of $\overline{K}[C]_P$. We extend $ord_P$ to the rational function field in the following way

$$\begin{aligned} \text{ord}_P : \overline{K}(C) &\to \mathbb{Z} \cup \{\infty\} \\ f/g &\mapsto \text{ord}_P(f) - \text{ord}_P(g). \end{aligned}$$

Let $\omega \in \Omega_C$, then at any fixed smooth point $P$, $\omega$ can be written as $f \, dt_P$ where $t_P$ is a local uniformizing element with $\text{ord}_P(t_P) = 1$ and $f \in \overline{K}(C)$. The valuation $\text{ord}_P(f)$ only depends on $\omega$, and we denote it by $\text{ord}_P(\omega)$.

A field $K$ is called a *non-Archimedean local field*, if it is complete, i.e. every Cauchy sequence is convergent, with respect to a discrete non-Archimedean absolute value $| \ |$, and if its *residue field* is finite.

### 1.6.1 Absolute Values over rational Numbers $\mathbb{Q}$

In case of $K = \mathbb{Q}$, take a prime $p > 1$. Every $x \in \mathbb{Q}$ can be written as $x = p^n \frac{a}{b}$ with $n, a \in \mathbb{Z}$, $b \in \mathbb{N}$, $b \neq 0$, $\gcd(a,b) = 1$, $p \nmid a$, $p \nmid b$. The map

$$\begin{aligned} v_p : \mathbb{Q} &\to \mathbb{R} \cup \{\infty\} \\ x &\mapsto \begin{cases} n & \text{if } 0 \neq x = p^n \frac{a}{b} \\ \infty & \text{if } x = 0 \end{cases} \end{aligned}$$

is a discrete valuation and

$$\begin{aligned} |\ |_p : \mathbb{Q} &\to \mathbb{R} \\ x &\mapsto \begin{cases} p^{-v_p(x)} & \text{if } x \neq 0 \\ 0 & \text{if } x = 0 \end{cases} \end{aligned}$$

is the *normalized p-adic absolute value*. The *p-adic numbers* $\mathbb{Q}_p$ are defined to be the completion of $\mathbb{Q}$ with respect to $|\ |_p$.

According to *Ostrowski*'s Theorem (see [CASS] Chapter II Section 3) every non-trivial absolute value on $\mathbb{Q}$ is either non-archimedian and equivalent to a $p$-adic absolute value $|\ |_p$ for a prime $p$, or it is archimedian and equivalent to the normal absolute value on $\mathbb{R}$, denoted by $|\ |_\infty$.

### 1.6.2 The Set $M_K$ of inequivalent absolute Values, Decomposition and Inertia Group

Let $M_K$ be a complete set of inequivalent absolute values on $K$. Then the archimedian absolute values in $M_K$ are denoted by $M_K^\infty$ and the non-archimedian absolute values in $M_K$ are denoted by $M_K^0$.

Let $K$ be any field with a valuation $v$ on $K$, let $L/K$ be a finite Galois extension. Let $S_v$ be the set of equivalence classes of extensions of $v$ to $L$, i.e. for every valuation $v'$ on $L$ satisfying $v'(k) = v(k)$ for all $k \in K$ one has $[v'] \in S_v$. The Galois group $G_{L/K}$ acts on $S_v$ by composition:

$$\begin{aligned} G_{L/K} \times S_v &\to S_v \\ (g, [v']) &\mapsto [v' \circ \sigma] \end{aligned}$$

and one can prove that the action of $G_{L/K}$ on $S_v$ is independent of the choice of the representative $v'$ for the equivalence class $[v']$. For any extension $v'$ of $v$ on $L$ let

$$G_{v'} := \operatorname{stab}([v']) = \{\sigma \in G_{L/K} \mid \sigma([v']) = [v']\}$$

be the decomposition group of $v'$ in $S_v$. The *inertia group* of $v'$ is defined to be the subgroup

$$I_{v'} := \{\sigma \in G_{v'} \mid \sigma(x) = x \mod \mathfrak{m}_{v'} \text{ for all } x \in R_{v'}\}$$

and it depends on the choice of the extension $v'$ of $v$.

By definition there is a homomorphism $G_{v'} \to G_{\overline{K}_{v'}/K_{v'}}$. It is a fact that this is a surjection with kernel $I_{v'}$.

## 1.7 Algebraic Curves

(see also [SIL1], Chapter II, Sections 3 to 5)

An *algebraic curve* is a projective variety of dimension one.

### 1.7.1 Divisors

The *divisor group* of a curve $C$, denoted by $\text{Div}(C)$, is the free abelian group generated by the points of $C$. Hence, a divisor $D$ is a formal finite linear combination $\sum_{P \in C} n_P(P)$ with $n_p \in \mathbb{Z}$.

The *degree* of a divisor is defined to be $\deg D = \sum_{P \in C} n_P$. A divisor $D = \sum_{P \in C} n_P(P)$ is called *positive* or *effective*, if $n_P \geq 0$ for all $P \in C$. In this case, we write $D \geq 0$. For two divisors $D_1$ and $D_2$ we write $D_1 \geq D_2$ if $D_1 - D_2 \geq 0$.

Assume $C$ is a smooth curve. Any nonconstant $f \in \overline{K}(C)^*$ has a finite number of poles and of zeros only, and the divisor associated to $f$ is defined to be

$$\text{div}(f) = \sum_{P \in C} \text{ord}_P(f)(P).$$

A divisor of the form $D = div(f)$ for some $f \in \overline{K}(C)^*$ is called a *principal divisor*.

It is a fact that the degree of a principal divisor is zero.

Analogously, the *associated divisor* to $\omega \in \Omega_C$ is defined to be

$$\text{div}(\omega) = \sum_{P \in C} \text{ord}_P(\omega)(P).$$

Two divisors $D_1$ and $D_2$ are called *linearly equivalent*, denoted by $D_1 \sim D_2$, if $D_1 - D_2$ is a principal divisor.

### 1.7.2 The Picard Group

Let $\text{princ}(C) \subseteq Div(C)$ denote the group of principal divisors of $C$. The *divisor class group* or *Picard group* is defined to be the quotient $Pic(C) := Div(C)/princ(C)$.

The Galois group $G_{\overline{K}/K}$ of $\overline{K}/K$ acts on $Div(C)$. The subgroup of $Div(C)$ fixed by $G_{\overline{K}/K}$ is denoted by $Pic_K(C)$.

$\text{Pic}^0$ denotes the subgroup of $\text{Pic}(C)$ consisting of the divisor classes of degree 0. Analogously, $\text{Pic}^0_K$ is the subgroup of $\text{Pic}_K(C)$ consisting of the divisors of degree 0.

Let $\omega_1, \omega_2 \in \Omega_C$ be two differentials. Then there exists a function $f \in \overline{K}(C)^*$ s.t. $\omega_1 = f\omega_2$. This leads to $div(\omega_1) = div(f) + div(\omega_2)$. Thus $\omega_1$ and $\omega_2$ are linearly equivalent, hence contained in the same divisor class in $\text{Pic}(C)$. This class is called the *canonical divisor class on $C$*. The elements in $Div(C)$ are called *canonical divisors*.

### 1.7.3 The Genus of a Curve

Let $C$ be a nonsingular projective curve and let $P$ be a fixed point of $C$.

Let $\omega \in \Omega_C$, then we can write $\omega = f\, dt_P$ where $t_P$ is a local uniformizing element with $\mathrm{ord}_P(t_P) = 1$ and $f \in \overline{K}(C)$.

It is a fundamental fact that the only functions which are regular everywhere are constant functions. But there do exist nonconstant differential forms that are regular everywhere:

A differential $\omega \in \Omega_C$ is called *regular* at $P$, if $\mathrm{ord}_P(\omega) \geq 0$. Because of $\mathrm{ord}_P(\omega) = \mathrm{ord}_P(f)$, one has

$$\omega \text{ regular at } P \Leftrightarrow \mathrm{ord}_P(f) \geq 0.$$

So, we define

$$\Omega[K] := \{\omega \in \Omega_C | \omega \text{ regular at every point}\}.$$

This is a vectorspace and its dimension

$$g := \dim \Omega[K]$$

is the *genus* of $C$.

## 1.8 Elliptic Curves

([SIL1] Chapter 3, Section 3) An *elliptic curve* is a pair $(E, \mathcal{O})$, where $E$ is a smooth projective algebraic curve (i.e. a smooth onedimensional projective variety) of genus one together with a *point at infinity* $\mathcal{O} \in E$.

The elliptic curve $E$ is *defined over $K$*, denoted by $E/K$, if $E$ is defined over $K$ as a projective algebraic set and $O \in E(K)$.

### 1.8.1 Weierstrass Equations

An elliptic curve can be identified with an affine equation of the form

$$E: \quad y^2 + a_1 xy + a_3 y = x^3 + a_2 x^2 + a_4 x + a_6,$$

with $a_1, ..., a_6 \in \overline{K}$. $E/K \Leftrightarrow a_1, ..., a_6 \in K$.

This is called the *Weierstrass equation* of the elliptic curve, it is unique up to transformations $(x, y) \mapsto (x', y')$ where $(x, y) = (u^2 x' + r, u^3 y' + su^2 x' + t)$ with $u, r, s, t \in \overline{K}$ and $u \neq 0$. Transformations of this form are denoted by $T(u, r, s, t)$.

If $\mathrm{char} K \neq 2, 3$, then a Weierstrass equation can be reduced to a simple form

$$E: \quad y^2 = x^3 + Ax + B,$$

called *reduced Weierstrass equation*. The *discriminant* of a reduced Weierstrass equation is given by

$$\Delta = -16(4A^3 + 27B^2),$$

and the *j-invariant* is given by

$$j = -1728 \frac{(4A)^3}{\Delta}.$$

In the case of char$(K) = 2,3$ there are other formulas for $\Delta$ and $j$, see also [SIL1] Appendix A, Proposition 1.1. Setting $r = s = t = 0$ in the transformations above, one gets $(x,y) \mapsto (x',y')$ where $(x,y) = (u^2x', u^3y')$ with $u \in \overline{K}$ and $u \neq 0$. These are the only transformations preserving the simple form $E : y^2 = x^3 + Ax + B$ of Weierstrass equations. One gets

$$E: \quad y'^2 = x'^3 + A'x' + B',$$

with $u^4 A' = A$ and $u^6 B' = B$, and $u^{12} \Delta' = \Delta$.

Let $K$ be a non-Archimedean local field, complete with respect to a discrete, normalized valuation $v$ on $K$, and let $R := \{x \in K : v(x) \geq 0\}$ be the ring of integers of $K$. A Weierstrass equation with integer coefficients s.t. $v(\Delta)$ is minimal is called *local minimal Weierstrass equation* at $v$.

Local minimal Weierstrass equations are unique up to transformations $T(u,r,s,t)$ with $r,s,t \in R$ and $u \in R^*$. It is not always possible to get a local minimal Weierstrass equation in the reduced form $E: y^2 = x^3 + Ax + B$.

A Weierstrass equation over a number field is called *global minimal Weierstrass equation* if the discriminant is minimal at every non-archimedean valuation $v$. A Global minimal Weierstrass equation exists if the ring of integers in $K$ is a principal ideal domain.

In particular, if $K = \mathbb{Q}$, we can always find transformations leading to a global minimal Weierstrass equation. For further details see [SIL1] Chapter VIII Section 8.

### 1.8.2 Cremona Tables

Note that in the *Cremona database* [CREM] equations of elliptic curves over $\mathbb{Q}$ are given in a specific form: Let $E$ be an elliptic curve given by a Weierstrass equation of the form

$$E: \quad y^2 + a_1 xy + a_3 y = x^3 + a_2 x^2 + a_4 x + a_6,$$

with $a_1, ..., a_6 \in \mathbb{Q}$. Using transformations $(x,y) \mapsto (x',y')$ where $(x,y) = (u^2 x', u^3 y')$ with $u \in \mathbb{Q} \setminus \{0\}$ each coefficient $a_i$ will be divided by $u^i$. So, one can change the equation into an *integral form*, i.e. all coefficients $a_i \in \mathbb{Z}$. *Global minimal models* for $E$ over $\mathbb{Q}$ are defined to be integral models s.t. $|\Delta|$ is minimal. These minimal models are unique up to isomorphisms with $u = \pm 1$ and $r,s,t \in \mathbb{Z}$. One can choose $r,s$ and $t$ s.t. $a_1, a_3 \in \{0,1\}$ and $a_2 \in \{-1,0,1\}$. In the cremona database, models of this form are called *reduced*. The reduced minimal model is unique. It is used in the tables to identify the corresponding elliptic curve.

### 1.8.3 The Group Law

Let $E$ be an elliptic curve. There is a group law on $E$ that can be constructed geometrically. For more information and explicit formulas see [SIL1] Chapter III, Section 2.

In particular, $E(K)$ is a subgroup of $E$, called the *Mordell-Weil group*.

### 1.8.4 Torsion Points

([SIL1] III.4) Let $E/K$ be an elliptic curve and let $m \in \mathbb{Z}$. Because of the goup structure we can define a map

$$\begin{aligned} [m] : E &\to E \\ P &\mapsto [m]P := \underbrace{P + \ldots + P}_{m \text{ times}}. \end{aligned}$$

This map is called *multiplication-by-m isogeny*. I will give a definition of *isogenies* later in section 1.10.3.

For $m \in \mathbb{N}_{>0}$ the set

$$E[m] := \{P \in E : [m]P = O\},$$

consisting of the points of $E$ of order dividing $m$, is called the *m-torsion subgroup of E*.

The set

$$E_{tors} := \bigcup_{m=1}^{\infty} E[m],$$

consisting of all points of finite order is called *torsion subgroup of E*. If $E$ is defined over $K$, then the points of finite order in $E(K)$ are denoted by $E_{tors}(K)$.

### 1.8.5 The Mordell-Weil Theorem

Let $K$ be a number field. Let $E/K$ be an elliptic curve defined over the field $K$. The *Mordell-Weil theorem* states that the group $E(K)$ is finitely generated. One proves first the *weak-Mordell-Weil theorem*: Let $m \in \mathbb{N}_{\geq 2}$. Then the *weak Mordell-Weil group $E(K)/mE(K)$* is a finite group. After that, the Mordell-Weil Theorem is proved by using a descent procedure, and a theory of height functions on the points of elliptic curves. This is a major theorem, which I do not prove here. The case of $K = \mathbb{Q}$ is given in [SIL2]. The general case is proved in [SIL1] Chapter VIII. From the Mordell-Weil Theorem one has

$$E(K) = E(K)_{\text{tors}} \times \mathbb{Z}^r,$$

where the torsion subgroup $E(K)_{\text{tors}}$ is finite. The number $r \in \mathbb{N}$ is the *rank* of $E$.

## 1.9 Hyperelliptic Curves

Let $g > 1$ and let $(\phi_1, \ldots, \phi_g) \in \Omega[K]^g$ be a basis of regular differentials for a nonsingular projective curve. Since local rings are discrete valuation rings the elements $\phi_i$ can be represented locally by a product $\phi_i = f_i dt_p$, with $f_i \in \overline{K}(C)$ regular at $P$ and a local uniformizing element $t_p$. It is a fact that the $f_i(x)$ do not have a common zero. The functions $f_i$ are unique up to a factor, so we can define

$$\begin{aligned} \phi : C &\to \mathbb{P}^{g-1} \\ x &\mapsto [f_1(x), \ldots, f_g(x)]. \end{aligned}$$

This map either is an *embedding*, i.e. injective with non-singular image, or a 2-1 map. Those $C$, for which $\phi$ is 2-1, are called *hyperelliptic curves*.

Note that in case of $g = 2$, $\phi : C \to \mathbb{P}^1$ is always a 2-1 map.

Explicitly, the affine part of an hyperelliptic curve is given by an equation of the form
$$C_0 : y^2 + h(x)y = f(x),$$
where $f, h \in K[X]$ such that $\mathrm{disc}(f), \mathrm{disc}(h) \neq 0$, i.e. $C_0$ is smooth, and
$$\begin{aligned}\deg f &> 0 \text{ and} \\ \deg h &\leq g \text{ with} \\ g &= \begin{cases} \frac{\deg f - 1}{2} & \text{if } \deg f \text{ is odd} \\ \frac{\deg f - 2}{2} & \text{if } \deg f \text{ is even.} \end{cases}\end{aligned}$$

Note that the equation of $C_0$ is quadratic in $y$, hence $(x, y) \mapsto x$ gives a 2-1 map to affine space. One can show that homogenization would lead to a curve in $\mathbb{P}^2$, but it would not be smooth at infinity, whenever $\deg f \geq 4$.

On the other hand, the map
$$[1, x, x^2, \ldots x^{g+1}, y] : \quad C_0 \to \mathbb{P}^{g+2}$$
leads to a projective algebraic curve $C/K$ that is nonsingular at every point. Then $C$ is a *hyperelliptic curve* of genus $g$.

In the special case $g = 1$ one has *elliptic curves*.

For more details, see [SIL1] Section II, Example 2.5.1.

Let $[X_0, \ldots, X_{g-1}] = [1, x, x^2, \ldots x^{g-1}]$, then $C \cap \{X_0 \neq 0\}$ is isomorphic to $C_0$.

Note that in the case of $\mathrm{char} K \neq 2$, the equation of the affine curve $C_0$ can be reduced to
$$C_0 : y^2 = f(x),$$
where $f \in K[X]$ with $\deg f > 0$.

## 1.10 Maps

### 1.10.1 Morphisms between Algebraic Varieties

Let $V_1 \in \mathbb{P}^m$ and $V_2 \in \mathbb{P}^n$ be projective varieties. A *rational map* from $V_1$ to $V_2$, is a map of the form
$$\begin{aligned}\phi : V_1 &\to V_2 \\ P &\mapsto [f_0(P), ..., f_n(P)],\end{aligned}$$
where $f_0, ..., f_n \in \overline{K}(V_2)$, such that for all $P \in V_1$:
If $f_i \in \overline{K}[V_1]_P$ for all $i \in \{0, ..., n\}$ (i.e. all $f_i$ are *regular* or *defined* at $P$), then $\phi(P) \in V_2$.

We say the rational map $\phi = [f_0, ..., f_n]$ is *defined over* $K$, if there exists a $\lambda \in \overline{K}^*$ such that $\lambda f_0, \ldots, \lambda f_n \in K(V_1)$.

A rational map is called *regular* or *defined* at $P$ if there is a $g \in \overline{K}(V_1)$ with the properties

$$gf_i \in \overline{K}[V_1]_P \text{ for all } i \in \{0, ..., n\} \text{ and}$$
$$\text{there exists } i \in \{0, ..., n\} \text{ s.t. } gf_i(P) \neq 0.$$

A *morphism* from $V_1$ to $V_2$ is a rational map that is regular at every point.

If $f_0, \ldots f_n$ are homogeneous polynomials of degree $d$, then $\phi$ is called a morphism of *degree d*.

A morphism $\phi : V_1 \to V_2$ is called an *isomorphism* or $\overline{K}$-*isomorphism* between $V_1$ and $V_2$, if there is a morphism $\psi : V_2 \to V_1$ such that $\psi \circ \phi = id_{V_1}$ and $\phi \circ \psi = id_{V_2}$.

We call it a $K$-*isomorphism* if $\phi$ and $\psi$ can be defined over $K$.

### 1.10.2 Homomorphisms between Elliptic Curves

In general, a *homomorphism* is a *morphism* preserving algebraic structures. There is a group law on elliptic curves, so in our case we define homomorphisms between elliptic curves as follows:

Let $(E_1, O_1)$ and $(E_2, O_2)$ be elliptic curves over $K$. A *homomorphism* or $\overline{K}$-*homomorphism* between $E_1$ and $E_2$ is a $\overline{K}$-morphism $\phi : E_1 \to E_2$ between varieties, preserving the group structure of $E_1$.

Analogously, a homomorphism defined over $K$ is called $K$-homomorphism.

A $(K\text{-})$homomorphism $\phi : E_1 \to E_2$ is called a *(K-)isomorphism between elliptic curves*, if there is a $(K\text{-})$homomorphism $\psi : E_1 \to E_2$, such that $\psi \circ \phi = id_{E_1}$ and $\phi \circ \psi = id_{E_2}$. Note that $(K\text{-})$ isomorphisms between elliptic curves need to be $(K\text{-})$homomorphisms and not just $(K\text{-})$morphisms of varieties. It is a theorem, that any morphism $\phi : E_1 \to E_2$ satisfying $\phi(O_{E_1}) = O_{E_2}$ is always a $K$-homomorphism.

### 1.10.3 Isogenies

(see also [SIL1] Chapter 3, Section 4) Let $(E_1, O_{E_1})$ and $(E_2, O_{E_2})$ be elliptic curves. A morphism $\phi : E_1 \to E_2$ satisfying $\phi(O_{E_1}) = O_{E_2}$ is called an *isogeny* from $E_1$ to $E_2$. According to [SIL1] Chapter II, Theorem 4.8, every isogeny is also a homomorphism. Two curves $(E, O_{E_1})$ and $(E_2, O_{E_2})$ are called *isogeneous* if there exists an isogeny $\phi : E_1 \to E_2$ satisfying $\phi(E_1) \neq \{O_{E_2}\}$, i.e. $\phi$ is to be nonconstant.

We denote the kernel of $\phi$ by $E_1[\phi] := \{P \in E_1 : \phi(P) = O_{E_2}\}$.

For each nonconstant isogeny $\phi : E_1 \to E_2$ of degree $m$ there exists a unique isogeny $\hat{\phi} : E_2 \to E_1$, such that $\hat{\phi} \circ \phi = [m]$. The isogeny $\hat{\phi}$ is called the *dual isogeny* of $\phi$.

## 2 Twisting

**Notation 2.1** *Let $(E, O)$ be an elliptic curve over the field $K$ and let $\overline{K}$ be an algebraic closure of $K$. We use the following notations:*

$$\begin{aligned}
\text{Isom}(E) &:= \text{the automorhism group of } \overline{K}\text{-isomorphisms } \phi : E \to E \\
\text{Aut}(E) &:= \text{the subgroup of Isom}(E) \text{ consisting of those } \phi : E \to E \\
&\quad \text{mapping } O \text{ to } O \\
\text{Isom}_K(E) &:= \text{the subgroup of Isom}(E) \text{ consisting of those } \phi : E \to E \\
&\quad \text{that are also } K\text{-isomorphisms}
\end{aligned}$$

*Note the difference between* $\text{Isom}(E)$ *and* $\text{Aut}(E)$: *Elements in* $\text{Aut}(E)$ *are group isomorphisms.*

**Definition 2.2** (see also [SIL1] Chapter 10, Section 2 and Section 5) Let $(E, O)$ be an elliptic curve over $K$.

a) A *twist* of the *smooth projective curve* $C$ is a smooth curve $C'$ over $K$ that is isomorphic to $C$ over $\overline{K}$. We denote the set of $K$-isomorphism classes of twists by $\text{Twist}(C)$.

b) A twist of the *elliptic curve* $(E, O)$ is an elliptic curve $(E', O')$ that is isogeneous to $E$ over $\overline{K}$, i.e. there exists a $\overline{K}$-isomorphism $\phi : E \to E'$ over $\overline{K}$ mapping $O$ to $O'$. We denote the set of twists of $(E, O)$ modulo $K$-isomorphism by $\text{Twist}(E, O)$. If $E/K$, we write $\text{Twist}((E, O)/K)$.

**Theorem 2.3** *(see also [SIL1] Chapter 10, Section 2, Theorem 2.2 ) Let* $C/K$ *be a smooth projective curve. One can identify*

$$\text{Twist}(C/K) \cong H^1(G_{\overline{K}/K}, \text{Isom}(C))$$

*via the following natural bijection: Let* $C'/K$ *be a twist of* $C/K$ *with an* $\overline{K}$-*isomorphism* $\phi : C' \to C$. *Then the corresponding element in* $H^1(G_{\overline{K}/K}, Isom(C))$ *is represented by the 1-cocycle*

$$\begin{aligned}
\xi : G_{\overline{K}/K} &\to \text{Isom}(C) \\
\sigma &\mapsto \xi_\sigma := \phi^\sigma \phi^{-1}.
\end{aligned}$$

*Proof.* We show surjectivity. Let $\xi : G_{\overline{K}/K} \to \text{Isom}(C)$, $\sigma \mapsto \xi_\sigma$ be a cocycle. Let $\overline{K}(C)_\xi$ be an abstract field extension of $\overline{K}$ s.t. $\overline{K}(C)_\xi$ is isomorphic to $\overline{K}(C)$ via an isomorphism $Z : \overline{K}(C) \to \overline{K}(C)_\xi$ and s.t. the Galois action of $G_{\overline{K}/K}$ is "twisted by $\xi$", i.e. given by $Z(f)^\sigma = Z(f^\sigma \xi_\sigma)$ for all $f \in \overline{K}(C)$ and all $\sigma \in G_{\overline{K}/K}$.

Let $\mathcal{F}$ denote the subfield of $\overline{K}(C)_\xi$ consisting of all elements that are fixed by the Galois group.

a) $\mathcal{F} \cap \overline{K} = K$: Let $F \in \mathcal{F} \cap \overline{K}$.

Because of $F \in \mathcal{F} \subseteq \overline{K}(C)_\xi$ and because $Z$ is an isomorphism, there exists $f \in \overline{K}(C)$ s.t. $F = Z(f)$.

Because of $F \in \overline{K}$ we get $f \in \overline{K}$ since the restriction of $Z$ on $\overline{K}$ is the identity on $\overline{K}$.

Further, $f^\sigma \in \overline{K}$ is constant for every $\sigma \in G_{\overline{K}/K}$, thus uneffected by any element of $\text{Isom}(C)$, in particular

$$f^\sigma \xi_\sigma = f^\sigma \text{ for all } \sigma \in G_{\overline{K}/K}.$$

Because of $Z(f) \in \mathcal{F}$ one has
$$Z(f) = Z(f)^\sigma = Z(f^\sigma \xi_\sigma) \text{ for all } \sigma \in G_{\overline{K}/K}.$$
This leads to $Z(f) = Z(f^\sigma)$, and since $Z$ is an isomorphism one has $f = f^\sigma$ for all $\sigma \in G_{\overline{K}/K}$. Hence $f \in K$ and finally $F \in K$.

b) $\mathcal{F}$ is a finitely generated extension of $K$ of transcendence degree one:

$\overline{K}(C)_\xi$ is a $\overline{K}$-vectorspace, the Galois group acts continuously on $\overline{K}(C)_\xi$, and $\mathcal{F}$ is the subspace of $\overline{K}(C)_\xi$ that is fixed by the Galois Group. Then it is a general fact that $\overline{K}(C)_\xi$ has a finite basis of $\mathcal{F}$ vectors. In other words $\overline{K}\mathcal{F} = \overline{K}(C)_\xi$. This general fact follows directly from Hilbert's Theorem 90 applied to $GL_n(\overline{K})$, i.e. $H^1(G_{\overline{K}/K}, GL_n(\overline{K}) = 0$. For more details see [SIL1] Chp. II Lemma 5.8.1 and Ex. 2.12.

c) Field extensions fulfilling a) and b) together with field injections fixing $K$ are equivalent as a category to smooth projective curves over $K$ with nonconstant rational maps defined over $K$. See also [SIL1] Chp. II, Remark 2.2 and [HART] Chp. I, Corollary 6.12.

In contrast to [SIL1] I want to show that the equivalence of categrories leads to curves $C'$ and $C''$. In d) I show that $C' = C''$.

So, we can find a corresponding curve $C'/K$ s.t. $\mathcal{F} \cong K(C')$ and an isomorphism $\psi : C' \to C$ s.t. $Z|_{K(C)} = \psi^*$, where $\psi^* : K(C) \to K(C')$, $\psi^*(f) = f \circ \phi$ for every $f \in K(C)$.

Further, there is a smooth projective curve $C''/\overline{K}$ s.t. $\overline{K}(C)_\xi \cong \overline{K}(C'')$ and an isomorphism $\phi : C'' \to C$ s.t. $Z = \phi^*$, where $\phi^* : \overline{K}(C) \to \overline{K}(C'')$, $\phi^*(f) = f \circ \phi$ for every $f \in \overline{K}(C)$.

d) Because of
$$\overline{K}(C') = \overline{K}K(C') \stackrel{c)}{=} \overline{K}\mathcal{F} \stackrel{b)}{=} \overline{K}(C)_\xi \stackrel{Z}{\cong} \overline{K}(C)$$
the equivalence of categories leads to $C' = C''$.

e) For all $\sigma \in G_{\overline{K}/K}$ and all $f \in \overline{K}(C)$ we get
$$f^\sigma \phi^\sigma = (f\phi)^\sigma = \phi^*(f)^\sigma = Z(f)^\sigma = Z(f^\sigma \xi_\sigma) = \phi^*(f^\sigma \xi_\sigma) = f^\sigma \xi_\sigma \phi$$
and $\phi^\sigma = \xi_\sigma \phi$, i.e. $\xi_\sigma = \phi^\sigma \phi^{-1}$. ∎

**Example 2.4** (Constructing a quadratic Twist)
Let $\operatorname{char}(K) \neq 2$, let $K(\sqrt{d})/K$ be a quadratic extension, and for each $\sigma \in G_{\overline{K}/K}$ let
$$\chi(\sigma) := \begin{cases} +1 & \text{if } \sqrt{d}^\sigma = \sqrt{d} \\ -1 & \text{if } \sqrt{d}^\sigma = -\sqrt{d} \end{cases}$$
be the quadratic character giving the isomorphism $G_{K(\sqrt{d})/K} \cong \{1, -1\}$. Then
$$\begin{aligned} \xi : G_{\overline{K}/K} &\to \operatorname{Isom}(E) \\ \sigma &\mapsto \xi_\sigma := [\chi(\sigma)] \end{aligned}$$

is a 1-cocycle: For each $\sigma, \tau \in G_{\overline{K}/K}$ one has

$$
\begin{aligned}
\xi_{\sigma\tau} &= [\chi(\sigma\tau)] \\
&= \begin{cases} [1] & \text{if } \sqrt{d}^{\sigma\tau} = \sqrt{d} \\ [-1] & \text{if } \sqrt{d}^{\sigma\tau} = -\sqrt{d} \end{cases} \\
&= \begin{cases} [1] & \text{if } \left(\sqrt{d}^{\sigma} = \sqrt{d} \text{ and } \sqrt{d}^{\tau} = \sqrt{d}\right) \text{ or } \left(\sqrt{d}^{\sigma} = -\sqrt{d} \text{ and } \sqrt{d}^{\tau} = -\sqrt{d}\right) \\ [-1] & \text{if } \left(\sqrt{d}^{\sigma} = \sqrt{d} \text{ and } \sqrt{d}^{\tau} = -\sqrt{d}\right) \text{ or } \left(\sqrt{d}^{\sigma} = -\sqrt{d} \text{ and } (-\sqrt{d})^{\tau} = -\sqrt{d}\right) \end{cases} \\
&= \begin{cases} [1] & \text{if } (\xi_\sigma = [1] \text{ and } \xi_\tau = [1]) \text{ or } (\xi_\sigma = [-1] \text{ and } \xi_\tau = [-1]) \\ [-1] & \text{if } (\xi_\sigma = [1] \text{ and } \xi_\tau = [-1]) \text{ or } (\xi_\sigma = [1] \text{ and } \xi_\tau = [-1]) \end{cases} \\
&= \xi_\sigma \xi_\tau \\
&= (\xi_\sigma)^\tau \xi_\tau.
\end{aligned}
$$

The corresponding twist is given by a smooth curve $C/K$ that is isomorphic to $E$ via an $\overline{K}$-isomorphism $\phi: C \to E$. Let $\phi^*: \overline{K}(E) \to \overline{K}(C)$ denote the corresponding $\overline{K}$-isomorphism of function fields.

The Galois action on $\overline{K}(C)$ is given by

$$\phi^*(f)^\sigma = \phi^*(f^\sigma \xi_\sigma) \text{ for all } f \in \overline{K}(E) \text{ and } \sigma \in G_{\overline{K}/K}.$$

One has

$$
\begin{aligned}
\xi_\sigma: E &\to E \\
(x, y) &\mapsto [\chi(\sigma)](x, y) = (x, \chi(\sigma)y).
\end{aligned}
$$

For $h_x, h_y \in \overline{K}(E)$ with $h_x: (x, y) \mapsto x$ and $h_y: (x, y) \mapsto \frac{y}{\sqrt{d}}$ one has

$$
\begin{aligned}
h_x^\sigma: (x, y) &\mapsto \sigma(1)\, x = x = h_x(x, y) \\
h_y^\sigma: (x, y) &\mapsto \sigma\left(\frac{1}{\sqrt{d}}\right) y = \frac{\chi(\sigma)}{\sqrt{d}} y
\end{aligned}
$$

Further, $\phi^*(h_x) \in \overline{K}(C)$ and $\phi^*(h_y) \in \overline{K}(C)$ are fixed by the Galois group:

$$
\begin{aligned}
\phi^*(h_x)^\sigma &= \phi^*(h_x^\sigma \xi_\sigma) = \phi^*\big((x, y) \mapsto x\big) = \phi^*(h_x) \\
\phi^*(h_y)^\sigma &= \phi^*(h_y^\sigma \xi_\sigma) = \phi^*\left((x, y) \mapsto \frac{1}{\sqrt{d}} \chi(\sigma)^2 y\right) = \phi^*(h_y)
\end{aligned}
$$

Let $y^2 = f(x)$ be a Weierstrass equation of $E$. Then the values of $(h_x, h_y)$ represent an elliptic curve satisfying the Weierstrass equation

$$d \cdot h_y^2 = f(h_x).$$

This curve represents the twist $C$ of $E$ by $\chi$. It is isomorphic to $E$ over $K(\sqrt{d})$ by the identification $\phi^{-1}: E \to C$, $(x, y) \mapsto (x, \sqrt{d}y)$. One checks

$$\phi^\sigma \phi^{-1}(x, y) = \phi^\sigma(x, \sqrt{d}y) = \left(\sigma(1)x, \sigma\left(\frac{1}{\sqrt{d}}\right)\sqrt{d}y\right) = (x, \chi(\sigma)y) = \xi_\sigma(x, y).$$

∎

**Proposition 2.5** *(see also [SIL1] Chapter 10, Section 5, Theorem 5.3(a))*

$$\text{Twist}((E,O)/K) \cong H^1(G_{\overline{K}/K}, \text{Aut}(E))$$

*Proof.* This can be shown analogously as in Theorem 2.3: The difference lies in the fact that cocycles $\xi_\sigma : E \to E$ need to be elements in $Aut(E)$, hence they are group isomorphisms, and for each twist $(E', O)$ of $(E, O)$ the $\overline{K}$-isomorphism $\phi : E \to E'$ needs to map $O$ to $O'$. ∎

**Proposition 2.6** *Let* char $K \neq 2, 3$ *and let $E$ be an elliptic curve over $K$. Let*

$$n := \begin{cases} 2 & \text{if } j(E) \neq 0, 1728 \\ 4 & \text{if } j(E) = 1728 \\ 6 & \text{if } j(E) = 0. \end{cases}$$

*Let $\mu_n$ denote the group of $n^{th}$ roots of unity. Then as $G_{\overline{K}/K}$-modules*

$$\begin{aligned} \mu_n &\cong \text{Aut}(E) \\ \zeta &\mapsto [\zeta] : E \to E, (x,y) \mapsto (\zeta^2 x, \zeta^3 y). \end{aligned}$$

*In particular, $\text{Aut}(E)$ is finite and its order is given by $\text{ord}(\text{Aut}(E)) = n$.*

*Proof.* Using the Weierstrass equation to identify the elliptic curve E/K, automorphisms need to be of a specific form leading to $n^{th}$ roots. For more details see [SIL1] Chapter III Theorem 10.1 and Corollary 10.2. ∎

**Proposition 2.7** *a) Hilbert's Theorem 90:*

$$H^1(G_{\overline{K}/K}, \overline{K}^*) = 0.$$

*b) Kummer Duality: Let* $\text{char}(K) = 0$ *or* $\text{char}(K) \nmid n$. *Then*

$$H^1(G_{\overline{K}/K}, \mu_n) \cong K^*/(K^*)^n$$

*Proof.* a) [BSCH] Chapter 4 Section 8 Theorem 1 and Theorem 2.

b) This follows from a) and a long exact sequence of cohomology groups. For more details see [SIL1] Appendix B Prop 2.5(c). ∎

**Proposition 2.8** *(see also [SIL1] Chapter 10, Section 5, Theorem 5.4) Let* char $K \neq 2, 3$, *let $(E, O)$ be an elliptic curve over $K$ and let $n := ord(\text{Aut}(E))$. There is a canonical ismorphism*

$$\text{Twist}((E,O)/K) \cong K^*/(K^*)^n.$$

*Let $E : y^2 = x^3 + Ax + B$ be the reduced Weierstrass equation of $E$. Then for $D \in K^*/(K^*)^n$ the corresponding elliptic curve $E_D \in \text{Twist}((E,O),K)$ has Weierstrass equation*

$$E_D : y^2 = \begin{cases} x^3 + D^2 Ax + D^3 B & \text{if } j(E) \neq 0, 1728 \\ x^3 + DAx & \text{if } j(E) = 1728 \\ x^3 + DB & \text{if } j(E) = 0. \end{cases}$$

*Proof.* We have

$$\begin{aligned}
\text{Twist}((E,O)/K) &\cong H^1(G_{\overline{K}/K}, \text{Aut}(E)) &&\text{by Proposition 2.5}\\
&\cong H^1(G_{\overline{K}/K}, \mu_n) &&\text{by Proposition 2.6}\\
&\cong K^*/(K^*)^n &&\text{by Proposition 2.7b.}
\end{aligned}$$

We now compute equations for $E_D$ for $D \in (K^*)^n$.

- If $n = 2$: see Example 2.4
- If $n = 4$: Let $\delta \in \overline{K}$ s.t. $\delta^4 = D$. Then

$$\begin{aligned}
\xi : G_{\overline{K}/K} &\to \mu_4\\
\sigma &\mapsto \xi_\sigma := \sigma(\delta)/\delta
\end{aligned}$$

is a 1-cocycle: For each $\sigma, \tau \in G_{\overline{K}/K}$ one has

$$\begin{aligned}
(\xi_\sigma)^\tau \xi_\tau &= \tau\left(\frac{\sigma(\delta)}{\delta}\right)\frac{\tau(\delta)}{\delta}\\
&= \frac{\tau\sigma(\delta)}{\tau(\delta)}\frac{\tau(\delta)}{\delta}\\
&= \frac{\tau\sigma(\delta)}{\delta}\\
&= \xi_{\sigma\tau}.
\end{aligned}$$

Take the isomorphism

$$\begin{aligned}
\psi : \mu_4 &\to \text{Aut}(E)\\
\zeta &\mapsto \phi(\zeta) : E \to E, (x,y) \mapsto (\zeta^2 x, \zeta^3 x).
\end{aligned}$$

Then $\psi(\xi_\sigma) \in H^1(G_{\overline{K}/K}, \text{Aut}(E))$. Let $C/K$ be the corresponding twist according to 2.5.

The curve $C/K$ is isomorphic to $E$ via an $\overline{K}$-isomorphism $\phi : C \to E$. Let $\phi^* : \overline{K}(E) \to \overline{K}(C)$ denote the corresponding $\overline{K}$-isomorphism of function fields.

The Galois action on $\overline{K}(C)$ is given by

$$\phi^*(f)^\sigma = \phi^*(f^\sigma \xi_\sigma) \text{ for all } f \in \overline{K}(E) \text{ and } \sigma \in G_{\overline{K}/K}.$$

For $h_x, h_y \in \overline{K}(E)$ with $h_x : (x,y) \mapsto \delta^{-2} x$ and $h_y : (x,y) \mapsto \delta^{-3} y$ one has

$$\begin{aligned}
h_x^\sigma : (x,y) &\mapsto \sigma(\delta^{-2}) \, x\\
h_y^\sigma : (x,y) &\mapsto \sigma(\delta^{-3}) \, y
\end{aligned}$$

Further, $\phi^*(h_x) \in \overline{K}(C)$ and $\phi^*(h_y) \in \overline{K}(C)$ are fixed by the Galois group:

$$\begin{aligned}
\phi^*(h_x)^\sigma &= \phi^*(h_x^\sigma \psi(\xi_\sigma)) = \phi^*\big((x,y) \mapsto \sigma(\delta)^{-2}(\sigma(\delta)/\delta)^2 x\big) = \phi^*(h_x)\\
\phi^*(h_y)^\sigma &= \phi^*(h_y^\sigma \psi(\xi_\sigma)) = \phi^*\big((x,y) \mapsto \sigma(\delta)^{-3}(\sigma(\delta)/\delta)^3 x\big) = \phi^*(h_y)
\end{aligned}$$

Let $y^2 = x^3 + Ax$ be a Weierstrass equation of $E$. Then the values of $(h_x, h_y)$ satisfy

$$\begin{aligned} h_y(x,y)^2 &= (\delta^{-3} y)^2 \\ &= \delta^{-6}(x^3 + Ax) \\ &= \delta^{-6} x^3 + A\delta^{-6} x \\ &= \delta^{-6} x^3 + \delta^{-4} A \delta^{-2} x \\ &= (\delta^{-2} x)^3 + D^{-1} A \delta^{-2} x \\ &= f_x(x,y)^3 + D^{-1} A f_x(x,y). \end{aligned}$$

Setting $u := \delta^{-2} \in \overline{K}^*$, the substitution

$$\begin{aligned} (h_x, h_y) &= (u^2 f_y, u^3 f_y) \\ &= (\delta^{-4} f_x, \delta^{-6} f_y) \end{aligned}$$

leads to

$$\begin{aligned} (\delta^{-6} f_y(x,y))^2 &= (\delta^{-4} f_x(x,y))^3 + D^{-1} A \delta^{-4} f_x(x,y) \\ \Leftrightarrow \quad \delta^{-12} f_y(x,y)^2 &= \delta^{-12} f_x(x,y)^3 + \delta^{-8} A f_x(x,y) & | \cdot \delta^{12} = D^3 \\ \Leftrightarrow \quad f_y(x,y)^2 &= f_x(x,y)^3 + D f_x(x,y). \end{aligned}$$

- If $n = 6$: Let $\delta \in \overline{K}$ s.t. $\delta^6 = D$ and define $h_x$ and $h_y$ as in the previous case. Let $y^2 = x^3 + B$ be a Weierstrass equation of $E$. Then the values of $(h_x, h_y)$ satisfy

$$\begin{aligned} h_y(x,y)^2 &= (\delta^{-3} y)^2 \\ &= \delta^{-6}(x^3 + B) \\ &= \delta^{-6} x^3 + \delta^{-6} B \\ &= \delta^{-6} x^3 + D \delta^{-6} B \\ &= (\delta^{-2} x)^3 + \delta^{-6} B \\ &= f_x(x,y)^3 + \delta^{-6} B. \end{aligned}$$

Setting $u := \delta^{-5/3} \in \overline{K}^*$, the substitution

$$\begin{aligned} (h_x, h_y) &= (u^2 h_x, u^3 h_y) \\ &= ((\delta^{-5/3})^2 f_x, (\delta^{-5/3})^3 f_y) \\ &= (\delta^{-10/3} f_x, \delta^{-5} f_y) \end{aligned}$$

leads to

$$\begin{aligned} (\delta^{-5} f_y(x,y))^2 &= (\delta^{-10/3} f_x(x,y))^3 + \delta^{-6} B \\ \Leftrightarrow \quad \delta^{-10} f_y(x,y)^2 &= \delta^{-10} f_x(x,y)^3 + \delta^{-6} B & | \cdot \delta^{10} \\ \Leftrightarrow \quad f_y(x,y)^2 &= f_x(x,y)^3 + \delta^4 B \\ \Leftrightarrow \quad f_y(x,y)^2 &= f_x(x,y)^3 + DB. \end{aligned}$$ ∎

# 3 Homogeneous Spaces

**Definition 3.1** a) A *principal homogeneous space* or *G-Torsor* for an abelian Group $(G,+)$ consists of a pair $(S,\mu)$ where $S$ is a set and $\mu : S \times G \to S$ gives a simply transitive algebraic group action of $G$ on $S$.

b) Let $(E(\overline{K}),+)$ denote the group of the points on an elliptic curve $E/K$. If this group operates on the set $C(\overline{K})$ of points on a smooth curve $C/K$ via a simply transitive algebraic group action $\mu : C \times E \to C$ as a morphism of varieties defined over $K$, we say $C$ is a *homogeneous space* of $E$.

**Notation 3.2** *Instead of $\mu(p,P)$ we also write $p + P$. We define*

$$\begin{aligned} \nu : C \times C &\to E \\ (p,q) &\mapsto \text{ the unique } P \in E \text{ s.t. } \mu(p,P) = q \end{aligned}$$

*and write $p - q$ instead of $\nu(p,q)$.*

**Proposition 3.3** *Let $E/K$ be an elliptic curve and let $C/K$ be a homogeneous space for $E/K$. The curve $C/K$ is a twist of $E/K$.*

*Proof.* The isomorphism for the twist is given by the translation map

$$\begin{aligned} \Theta : E &\to C \\ P &\mapsto \mu(p_0, P) = p_0 + P, \end{aligned}$$

where $p_0 \in C$ is a fixed point. For further details see [SIL1] Chapter X, Proposition 3.2. ∎

## 3.1 Jacobian Varieties

In order to understand the meaning of homogeneous spaces, I want to identify them in a more abstract context with jacobian varieties. An *abelian variety* is defined to be a a nonsingular projective variety with a group structure; this is automatically abelian.

**Theorem 3.4** *Let $C/K$ be a curve of genus $g$. Then there exists an abelian variety $J(C)$ over $K$ of dimension $g$, such that*

$$\operatorname{Pic}_K^0(C) \cong J(C)(K) \text{ as abelian groups.}$$

*$J(C)$ is unique up to $K$-isomorphism and called the Jacobian of $C$.*

So, for a curve $C$ and an abelian variety $V$, we write $V = J(C)$, iff $V$ is $K$-isomorphic to $J(C)$ and we call $V$ *the* Jacobian of $C$.

*Proof (of Theorem 3.4).* See [CONF]. ∎

**Proposition 3.5** a) *Let $C/K$ be a curve of genus $1$. Then the Jacobian of $C$ is given by an elliptic curve $E/K$.*

b) *If $C/K$ is an elliptic curve, then $J(C) = C$.*

*Proof.* See [CONF]. ∎

**Theorem 3.6** *(see also [SUTH] Theorem 26.16.)*
*Let $E/K$ be an elliptic curve and let $C/K$ be a smooth curve. Then*

$$C \text{ is homogeneous space of } E \quad \Leftrightarrow \quad E = J(C).$$

*Proof.* See [CONF]. ∎

## 3.2 The Weil-Châtelet Group

**Definition 3.7** ([SIL1] X.3) Two homogeneous spaces $C/K$ and $C'/K$ are *equivalent* if there is an isomorphism $\Theta : C \to C'$, s.t.

$$\Theta(p + P) = \Theta(p) + P \text{ for all } p \in C \text{ and all } P \in E$$

**Definition 3.8** ([SIL1] X.3) The collection of equivalence classes of homogeneous spaces for $E/K$ is called the *Weil-Châtelet group* for $E/K$, denoted by $WC(E/K)$.

**Theorem 3.9** *Let $E/K$ be an elliptic curve. Then*

$$WC(E/K) \cong H^1(G_{\overline{K}/K}, E),$$

*where the cocycle corresponding to a homogeneous space $C/K$ of $E/K$ is given by*

$$\begin{aligned} G_{\overline{K}/K} &\to E \\ \sigma &\mapsto p_0^\sigma - p_0, \end{aligned}$$

*where $p_0$ is any point of $C$.*

*Proof.*
- One needs to show that the map given in Theorem 3.9 is well defined. I will not prove this in detail, referring to [SIL1] Chapter X, Theorem 3.6.:
  i) One can easily see that $\sigma \mapsto p_0^\sigma - p_0$ is a cocycle
  ii) Using [SIL1] Chapter X, Lemma 3.5 one can show that two equivalent homogeneous spaces are always mapped to the same cohomology class.

- Surjectivity: Let $\xi \in H^1(G_{\overline{K}/K}, E)$. We can embed $E$ into $\mathrm{Isom}(E)$ by identifying each point $P \in E$ with the translation $\tau_{-P} = \tau_P^{-1} \in \mathrm{Isom}(E)$. In particular, for $\xi_\sigma \in E$, we regard $\tau_{-\xi_\sigma}$. Since $\tau_{-\xi} : \sigma \mapsto \tau_{-\xi_\sigma}$ is contained in $H^1(G_{\overline{K}/K}, \mathrm{Isom}(E))$, we can use Theorem 2.3: There exists a curve $C/K$ and a $\overline{K}$-isomorphism $\phi : C \to E$ s.t.

$$\phi^\sigma \phi^{-1} = \tau_{\xi_\sigma}^{-1}.$$

We now need to prove that $C$ forms a homogeneous space of $E$. We define

$$\begin{aligned} \mu : C \times E &\to C \\ (p, P) &\mapsto \phi^{-1}(\tau_P(\phi(p))) = \phi^{-1}(\phi(p) + P) \end{aligned}$$

and check that this gives a simply transitive group action of $E$ on $C$ over $K$:

i) For every $p \in E$: $\mu(p, O) = \phi^{-1}(\phi(p) + O) = \phi^{-1}(\phi(p)) = p$

ii) For every $p \in C$ and $P, Q \in E$ one has

$$\begin{aligned} \mu(\mu(p, P), Q) &= \phi^{-1}(\phi(\mu(p, P)) + Q) \\ &= \phi^{-1}(\phi(\phi^{-1}(\phi(p) + P)) + Q) \\ &= \phi^{-1}((\phi(p) + P) + Q) \\ &= \phi^{-1}(\phi(p) + (P + Q)) \\ &= \mu(p, P + Q) \end{aligned}$$

iii) For every $p, q \in C$ one has

$$\begin{aligned} \mu(p, P) = q &\Leftrightarrow \phi^{-1}(\phi(p) + P) = q & | \phi() \\ &\Leftrightarrow \phi(p) + P = \phi(q) & | -\phi(P) \\ &\Leftrightarrow P = \phi(q) - \phi(p), \end{aligned}$$

which is unique for each pair of $p$ and $q$, according to [SIL1] Chapter X, Lemma 3.5.

iv) $\mu$ is defined over $K$, since it is fixed by the Galois group: From $\phi^\sigma \phi^{-1} = \tau_{\xi_\sigma}^{-1}$ one gets

$$\begin{aligned} \mu(p, P)^\sigma &= (\phi^{-1})^\sigma (\phi^\sigma(p^\sigma) + P^\sigma) & | \phi^\sigma = \tau_{\xi_\sigma}^{-1} \phi \\ &= (\phi^{-1})^\sigma ((\phi(p^\sigma) - \xi_\sigma) + P^\sigma) & | (\phi^{-1})^\sigma = \phi^{-1} \tau_{\xi_\sigma} \\ &= \phi^{-1}(((\phi(p^\sigma) - \xi_\sigma) + P^\sigma) + \xi_\sigma) \\ &= \phi^{-1}(\phi(p^\sigma) + P^\sigma) \\ &= \mu(p^\sigma, P^\sigma). \end{aligned}$$

We now show that the cocycle associated to the homogeneous space C/K is given by $\xi_\sigma$. Take any point $p_0 \in C$, e.g. $p_0 = \phi^{-1}(O)$ (This is a point of $C$ since $\phi : C \to E$). Then

$$\begin{aligned} p_0^\sigma - p_0 &= (\phi^\sigma)^{-1}(O) - \phi^{-1}(O) & | (\phi^\sigma)^{-1} = \phi^{-1} \tau_{\xi_\sigma} \\ &= \phi^{-1}(O + \xi_\sigma) - \phi^{-1}(O) & | \text{[SIL1] Chapter X, Lemma 3.5.} \\ &= O + \xi_\sigma - O \\ &= \xi_\sigma. \end{aligned}$$

- **Injectivity:** Assume $C/K$ and $C'/K$ are two homogeneous spaces. Let $\xi : \sigma \mapsto \xi_\sigma = p_0^\sigma - p_0$ and $\xi' : \sigma \mapsto \xi'_\sigma = p_0'^\sigma - p_0'$ be the corresponding cocycles, $p_0 \in C$, $p_0' \in C'$. Assume that $\xi$ and $\xi'$ are *cohomologous*, i.e. there exists $P_0 \in E$ s.t.

$$P_0^\sigma + \xi_\sigma = \xi'_\sigma + P_0.$$

Then

$$\begin{aligned} \Theta : C &\to C' \\ p &\mapsto p_0' + (p - p_0) + P_0. \end{aligned}$$

is an equivalence of $C$ and $C'$:

i) $\Theta$ is an isomorphism defined over $K$, since it is fixed by the Galois group, since for each $\sigma \in G_{\overline{K}/K}$:

$$\begin{aligned}
\Theta(p)^\sigma &= (p'_0 + (p - p_0) + P_0)^\sigma \\
&= p'^\sigma_0 + (p^\sigma - p^\sigma_0) + P^\sigma_0 \\
&= \underbrace{p'_0 + (p^\sigma - p_0) + P_0}_{=\Theta(p^\sigma)} + \underbrace{\left(P'^\sigma_0 + \underbrace{p'^\sigma_0 - p'_0}_{=\xi'_\sigma} - \underbrace{((p^\sigma_0 - p_0) + P_0)}_{=\xi_\sigma}\right)}_{=0} \\
&= \Theta(p^\sigma)
\end{aligned}$$

ii) $\Theta$ is compatible with the action of $E$ on $C$ and $C'$: For all $p \in C$ and all $P \in E$ one has

$$\begin{aligned}
\Theta(p + P) &= p'_0 + (p + P - p_0) + P_0) \\
&= \Theta(p) + P \quad \blacksquare
\end{aligned}$$

We summarize our results and get the following corollary:

**Corollary 3.10** *Let $E/K$ be an elliptic curve. Using the notations from above, Figure 1 shows our previous results and corresponding inclusions:*

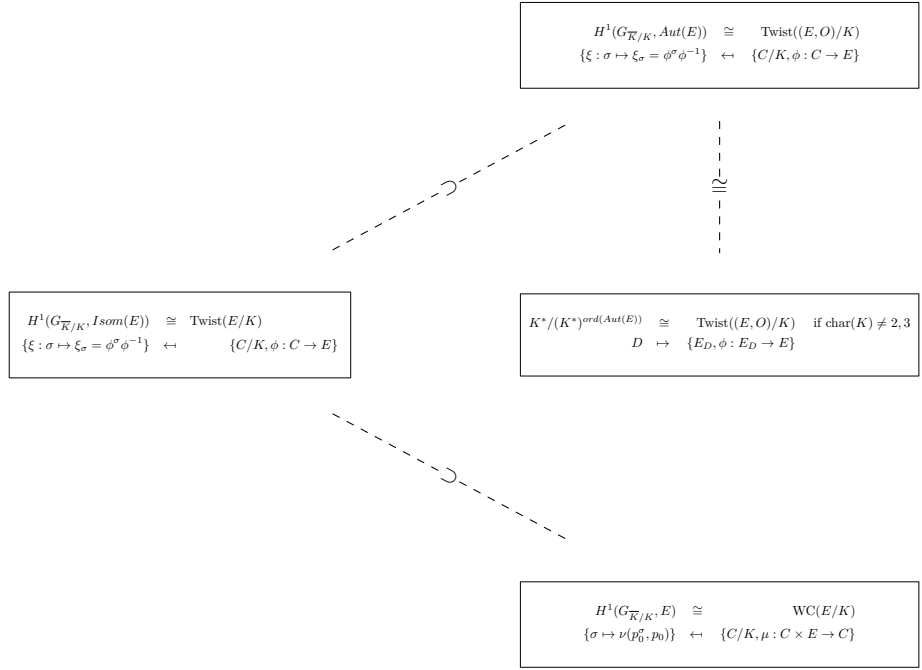

Figure 1: Summary of previous results and inclusions

### 3.3 Constructing a quadratic Homogeneous Space

**Proposition 3.11** *(See also [SIL1] Chapter X, Remark 3.7)*
*Assume* $\operatorname{char} K \neq 2$. *Let $K(\sqrt{d})/K$ be a quadratic extension and let $T \in E(K)$*

be a point of order 2. Let $E/K$ be an elliptic curve given by Weierstrass equation $E : y^2 = x^3 + ax^2 + bx$. The natural bijection defined in Theorem 3.9 maps the cocycle

$$\xi : G_{\overline{K}/K} \to E$$
$$\sigma \mapsto \begin{cases} O & \text{if } \sqrt{d}^\sigma = \sqrt{d} \\ T & \text{if } \sqrt{d}^\sigma = -\sqrt{d} \end{cases}$$

to a hyperellipic curve $C/K$, that can be represented by affine coordinates $(z, w)$ with equation

$$C : dw^2 = d^2 - 2adz^2 + (a^2 - 4b)z^4.$$

*Proof.* 1. We start constructing the twist $C$ of $E$ corresponding to $\xi$ by using the proof of Theorem 2.3 and Theorem 3.9: Each point $P$ of $E$ can be identified with a translation map $\tau_P : E \to E$, $p \mapsto p + P$, and $\xi$ can be regarded as a cocycle in $H^1(G_{\overline{K}/K}, Isom(E))$. Without loss of generality we can assume $E : y^2 = x^3 + ax^2 + bx$ and $T = (0,0)$ because of $\text{char}(K) \neq 2$ and $T \in E(K)$. Then we identify

$$\xi_\sigma = \begin{cases} \tau_O = id & \text{if } \sqrt{d}^\sigma = \sqrt{d} \\ \tau_T & \text{if } \sqrt{d}^\sigma = -\sqrt{d} \end{cases}$$

where

$$\tau_T : E \to E (x, y) \mapsto (x, y) + T = (x, y) + (0, 0) = \left(\frac{b}{x}, -\frac{by}{x^2}\right).$$

Let $C/K$ denote the twist of $E$ corresponding to $\xi$, and let $\phi : C \to E$ be the corresponding isomorphism, and let $\phi^* : \overline{K}(E) \to \overline{K}(C)$ denote the isomorphism of function fields.

Again, the Galois action on $\overline{K}(C)$ is given by

$$\phi^*(f)^\sigma = \phi^*(f^\sigma \xi_\sigma) \text{ for all } f \in \overline{K}(E) \text{ and } \sigma \in G_{\overline{K}/K}.$$

For $h_x, h_y \in \overline{K}(E)$ with $h_x : (x,y) \mapsto \frac{\sqrt{d}x}{y}$ and $h_y : (x,y) \mapsto \sqrt{d}\left(x - \frac{b}{x}\right)\left(\frac{x}{y}\right)^2$ one has

$$h_x^\sigma : (x, y) \mapsto \frac{\sqrt{d}^\sigma x}{y}$$
$$h_y^\sigma : (x, y) \mapsto \sqrt{d}^\sigma \left(x - \frac{b}{x}\right)\left(\frac{x}{y}\right)^2$$

Further, $\phi^*(h_1) \in \overline{K}(C)$ and $\phi^*(h_2) \in \overline{K}(C)$ are fixed by the Galois group:

- For $\sigma \in G_{\overline{K}/K}$ with $\sqrt{d}^\sigma = \sqrt{d}$ one has

$$\phi^*(h_x)^\sigma = \phi^*(h_x^\sigma \xi_\sigma) = \phi^*(h_x^\sigma \tau_O) = \phi^*(h_x^\sigma)$$
$$= \phi^*(h_x)$$

and
$$\phi^*(h_y)^\sigma = \phi^*(h_y^\sigma \xi_\sigma) = \phi^*(h_x^\sigma \tau_O) = \phi^*(h_x^\sigma)$$
$$= \phi^*(h_y).$$

- For $\sigma \in G_{\overline{K}/K}$ with $\sqrt{d}^\sigma = -\sqrt{d}$ one has
$$\phi^*(h_x)^\sigma = \phi^*(h_x^\sigma \xi_\sigma) = \phi^*(h_x^\sigma \tau_T)$$
$$= \phi^*\left((x,y) \mapsto \frac{\sqrt{d}^\sigma \frac{b}{x}}{\frac{-by}{x^2}}\right)$$
$$= \phi^*(h_x)$$

and
$$\phi^*(h_y)^\sigma = \phi^*(h_y^\sigma \xi_\sigma) = \phi^*(h_x^\sigma \tau_T)$$
$$= \phi^*\left((x,y) \mapsto \sqrt{d}^\sigma \left(\frac{b}{x} - \frac{b}{\frac{b}{x}}\right)\left(\frac{\frac{b}{x}}{\frac{-by}{x^2}}\right)^2\right)$$
$$= \phi^*(h_y).$$

Then the values of $(h_x, h_y)$ represent an hyperelliptic curve satisfying the equation

$$\begin{aligned}
d \cdot h_y(x,y)^2 &= d\left(\sqrt{d}\left(x - \frac{b}{x}\right)\left(\frac{x}{y}\right)^2\right)^2 & \mid h_x(x,y) = \frac{\sqrt{d}x}{y} \\
&= d\left(\sqrt{d}\left(x - \frac{b}{x}\right)\left(\frac{h_x(x,y)}{\sqrt{d}}\right)^2\right)^2 \\
&= \left(x - \frac{b}{x}\right)^2 h_x(x,y)^4 \\
&= \left(\left(x + \frac{b}{x}\right)^2 - 4b\right) h_x(x,y)^4 \\
&= \left(\left(\left(x + a + \frac{b}{x}\right) - a\right)^2 - 4b\right) h_x(x,y)^4 & \mid y^2 = x^3 + ax^2 + bx \\
&= \left(\left(\frac{y^2}{x^2} - a\right)^2 - 4b\right) h_x(x,y)^4 & \mid h_x(x,y) = \frac{\sqrt{d}x}{y} \\
&= \left(\left(\frac{d}{h_x(x,y)} - a\right)^2 - 4b\right) h_x(x,y)^4 \\
&= \left(\left(\frac{d}{h_x(x,y)}\right)^2 - 2\frac{ad}{h_x(x,y)^2} + a^2 - 4b\right) h_x(x,y)^4 \\
&= d^2 - 2adh_x(x,y)^2 + (a^2 - 4b) h_x(x,y)^4.
\end{aligned}$$

This curve represents the twist $C$ of $E$ corresponding to $\xi$ and the isomorphism is given by

$$\begin{aligned}
\phi : E &\to C \\
(x,y) &\mapsto (h_x(x,y), h_y(x,y)).
\end{aligned}$$

2. The next step is to verify that the natural bijection given in Theorem 3.9 maps $C$ to $\xi$. We write

$$\begin{aligned}
\phi : E &\to C \\
(x,y) &\mapsto (h_x(x,y), h_y(x,y)) \\
&= \left(\frac{\sqrt{d}x}{y}, \sqrt{d}\left(x - \frac{b}{x}\right)\left(\frac{x}{y}\right)^2\right) \\
&= \left(\frac{\sqrt{d}xy}{y^2}, \sqrt{d}\left(x - \frac{b}{x}\right)\frac{x^2}{y^2}\right) \qquad | \; y^2 = x^3 + ax^2 + bx \\
&= \left(\frac{\sqrt{d}y}{x^2 + ax + b}, \sqrt{d}\left(x - \frac{b}{x}\right)\frac{x}{x^2 + ax + b}\right) \\
&= \left(\frac{\sqrt{d}y}{x^2 + ax + b}, \frac{\sqrt{d}\left(x^2 - b\right)}{x^2 + ax + b}\right)
\end{aligned}$$

and can compute $\phi(0,0) = (0, -\sqrt{d})$ and $\phi(O) = (0, \sqrt{d})$. We choose $p := (0, \sqrt{d}) \in C$ and check

$$\begin{aligned}
\sigma \mapsto p^\sigma - p &= (0,\sqrt{d})^\sigma - (0,\sqrt{d}) \\
&= \begin{cases} (0,\sqrt{d}) - (0,\sqrt{d}) & \text{if } \sqrt{d}^\sigma = \sqrt{d} \\ (0,-\sqrt{d}) - (0,\sqrt{d}) & \text{if } \sqrt{d}^\sigma = -\sqrt{d} \end{cases} \\
&= \begin{cases} O & \text{if } \sqrt{d}^\sigma = \sqrt{d} \\ \phi^{-1}(0,-\sqrt{d}) - \phi^{-1}(0,\sqrt{d}) & \text{if } \sqrt{d}^\sigma = -\sqrt{d} \end{cases} \\
&= \begin{cases} O & \text{if } \sqrt{d}^\sigma = \sqrt{d} \\ (0,0) = T & \text{if } \sqrt{d}^\sigma = -\sqrt{d} \end{cases} \\
&\cong \xi_\sigma.
\end{aligned}$$
∎

# 4 The Selmer and Shafarevich-Tate Groups

According to Section 1.6.2, let $M_K$ denote a complete set of inequivalent absolute values on $K$, and for each $v \in M_K$ let $K_v \supset \overline{K}$ be the completion of $K$ at $v$ and $G_v \subset G_{\overline{K}/K}$ the corresponding decomposition group. For example, if $K = \mathbb{Q}$ let $\mathbb{Q}_p$ denote the $p$-adic numbers for $p$ prime.

**Proposition 4.1** *The following diagram is commutative:*

$$\begin{array}{ccccccccc}
0 & \longrightarrow & E'(K)/\phi(E(K)) & \longrightarrow & H^1(G_{\overline{K}/K}, E[\phi]) & \longrightarrow & WC(E/K)[\phi] & \longrightarrow & 0 \\
& & \downarrow & & \downarrow & & \downarrow & & \\
0 & \to & \prod_{v \in M_K} E'(K_v)/\phi(E(K_v)) & \to & \prod_{v \in M_K} H^1(G_v, E[\phi]) & \to & \prod_{v \in M_K} WC(E/K_v)[\phi] & \to & 0
\end{array}$$

*Proof.* We start with the short exact sequence

$$0 \longrightarrow E[\phi] \longrightarrow E \longrightarrow E' \longrightarrow 0$$
.

This gives a long exact sequence

$$0 \longrightarrow H^0(G_{\overline{K}/K}, E[\phi]) \longrightarrow H^0(G_{\overline{K}/K}, E) \longrightarrow H^0(G_{\overline{K}/K}, E') \xrightarrow{\delta}$$
$$H^1(G_{\overline{K}/K}, E[\phi]) \longrightarrow H^1(G_{\overline{K}/K}, E) \longrightarrow H^1(G_{\overline{K}/K}, E').$$

In section 1.4 we noted that for any variety $V$ its rational points can be expressed as follows

$$V(K) = \{P \in V : P^\sigma = P \text{ for all } \sigma \in G_{\overline{K}/K}\}.$$

We know from cohomology theory that the 0th cohomology group is always given by the elements fixed by the galois group. This yields

$$0 \longrightarrow E(K)[\phi] \longrightarrow E(K) \longrightarrow E'(K) \xrightarrow{\delta}$$
$$H^1(G_{\overline{K}/K}, E[\phi]) \longrightarrow H^1(G_{\overline{K}/K}, E) \longrightarrow H^1(G_{\overline{K}/K}, E').$$

This gives the fundamental short exact sequence

$$0 \longrightarrow E'(K)/\phi(E(K)) \xrightarrow{\delta} H^1(G_{\overline{K}/K}, E[\phi]) \longrightarrow H^1(G_{\overline{K}/K}, E)[\phi] \longrightarrow 0.$$

For each $v \in M_K$ the group $G_v$ acts on $E(\overline{K}_v)$ and $E'(\overline{K}_v)$. So, for each $v \in M_K$ we get an exact sequence

$$0 \longrightarrow E'(K_v)/\phi(E(K_v)) \xrightarrow{\delta} H^1(G_v, E[\phi]) \longrightarrow H^1(G_v, E)[\phi] \longrightarrow 0.$$

Using Theorem 3.9, one can replace the first cohomology groups for the module $E$ by the corresponding Weil-Châtelet groups, which gives the result. ∎

**Definition 4.2** (See also [SIL1] Chapter X, Section 4) The Selmer and Shafarevich-Tate groups are defined to be

$$S^\phi(E/K) := \ker\{H^1(G_{\overline{K}/K}, E[\phi]) \to \prod_{v \in M_K} WC(E/K_v)\}$$

$$\text{Ш}(E/K) := \ker\{WC(E/K) \to \prod_{v \in M_K} WC(E/K_v)\}.$$

## 4.1 Finiteness of the Selmer Group

**Theorem 4.3** *(See [SIL1] Chapter X, Theorem 4.2 b)*
*Let $E/K$, $E'/K$ be elliptic curves and let $\phi : E/K \to E'/K$ be an isogeny defined over $K$.*

*The Selmer Group $S^{(\phi)}(E/K)$ is finite.*

Figure 2 shows how this theorem is proved using the following propositions.

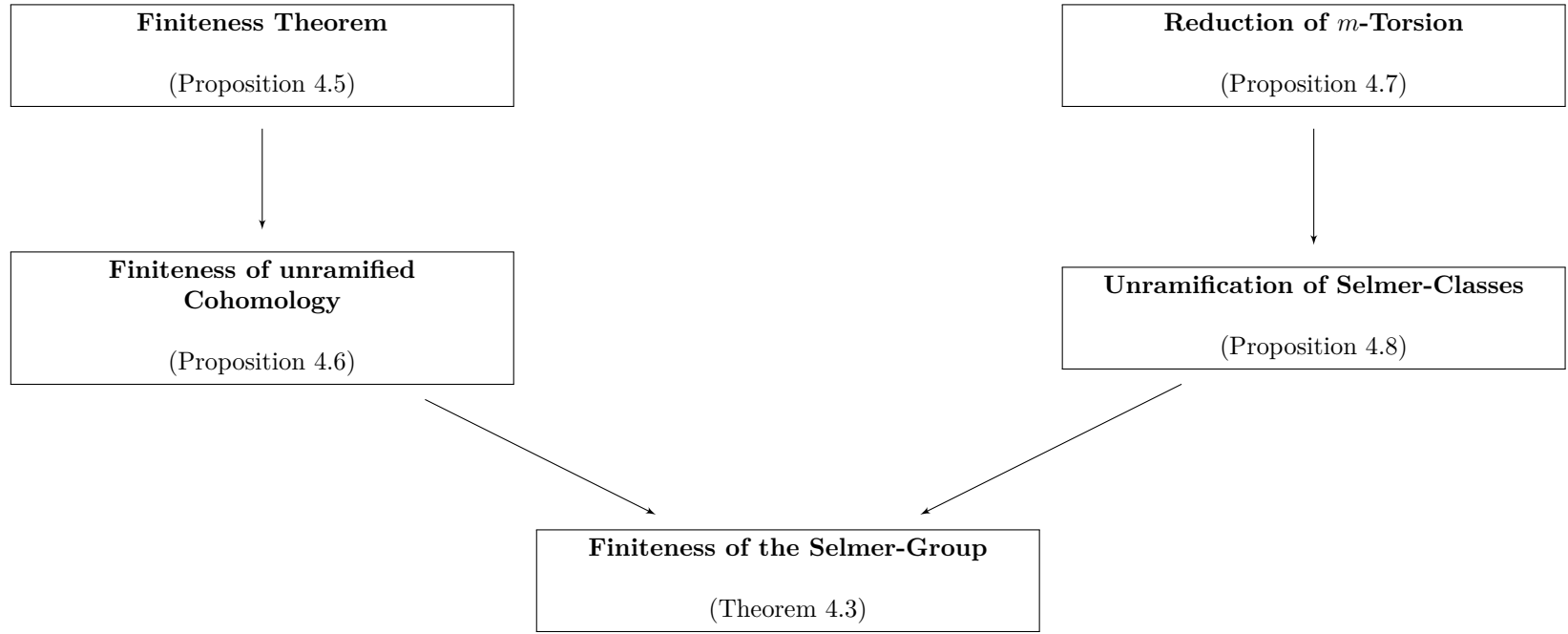

Figure 2: Structure of the proof of Theorem 4.3

**Definition 4.4** (See also [SIL1] Chapter VII Section 4)
Let $M$ be a $G_{\overline{K}/K}$-module and let $v \in M_K^0$.
An extension $L/K$ is called *unramified* at $v$, if the inertia subgroup $I_v$ of $G_{L/K}$ is trivial.
A cohomology class $\xi \in H^r(G_{\overline{K}/K}, M)$

$$\begin{aligned} \xi : G_{\overline{K}/K} &\to M \\ \sigma &\mapsto \xi_\sigma \end{aligned}$$

is called *unramified* at $v$ if

$$\xi_\sigma = 1_M \text{ for all } \sigma \in I_v,$$

where $1_M$ is the neutral element in the group $M$.

**Proposition 4.5** *(see also [SIL1] Chapter VIII, Proposition 1.6) Let $K$ be a number field, $S$ be a finite set of places s.t. $M_K^\infty \subset S \subset M_K$ and let Let $L/K$ be the maximal abelian Galois extension of exponent $m \geq 2$ that is unramified at $v$ for every $v \notin S$. Then $L/K$ is a finite extension.*

*Proof.* The proof uses fundamental facts of algebraic number theory. I will not prove this here, referrring to [SIL1] Chapter VIII, Proposition 1.6. ∎

**Proposition 4.6** *(See also [SIL1] Chapter X, Lemma 4.3) Let $K$ be a number field and let $\overline{K}$ be its algebraic closure. Let $M$ be a finite $G_{\overline{K}/K}$-module. Let $S$ be a finite set of places s.t. $M_K^\infty \subset S \subset M_K$. Then*

$$H^1(G_{\overline{K}/K}, M; S) = \{\xi \in H^1(G_{\overline{K}/K}, M) \mid \xi \text{ unramified at } v \text{ for } v \notin S\}$$

*is finite.*

*Proof.*  i) Reduction to the case of trivial action: Because $G_{\overline{K}/K}$ acts continuously on $M$ and $M$ is finite, there exists $H \subset G_{\overline{K}/K}$ of finite index s.t. $H$ acts trivially on $M$. Then

$$K' := Fix(H)$$

gives a finite extension of $K$ and

$$G_{K'/K} = G_{\overline{K}/K}/H.$$

$G_{\overline{K}/K'} = H$ acts trivially on $M$. There is an exact *Inflation-Restriction Sequence* (see also [SIL1] Appendix B, Proposition 2.4)

$$0 \to \underbrace{H^1(G_{K'/K}, M^H)}_{=H^1(G_{\overline{K}/K}/H, M^H)} \overset{\text{Inf}}{\to} H^1(G_{\overline{K}/K}, M) \overset{\text{Res}}{\to} \underbrace{H^1(G_{\overline{K}/K'}, M)}_{=H^1(H,M)}.$$

$H^1(G_{\overline{K}/K}/H, M^H)$ is finite, so one gets

$$H^1(H, M) \text{ finite} \Rightarrow H^1(G_{\overline{K}/K}, M) \text{ finite}.$$

ii) Now assume that the action of $G_{\overline{K}/K}$ is trivial. Then we know from cohomology theory that the first cohomology group is given by

$$H^1(G_{\overline{K}/K}, M; S) = \text{Hom}(G_{\overline{K}/K}, M; S).$$

$M$ is finite, so there exists $m \geq 2$ s.t. $m \cdot M = \{0\}$. So, let $L/K$ be the maximal abelian Galois extension of exponent $m \geq 2$ that is unramified at $v$ for every $v \notin S$. One has

$$\text{Hom}(G_{\overline{K}/K}, M; S) \stackrel{iii)}{=} \text{Hom}(G_{L/K}, M).$$

According to 4.5, $L/K$ is a finite extension and one gets the result.

iii) Proof of $\text{Hom}(G_{\overline{K}/K}, m; S) = \text{Hom}(G_{L/K}, M)$:

$\xi : G_{\overline{K}/K} \to M$ is a homomorphism to $M$ abelian of exponent $m$
$\Leftrightarrow \quad \xi$ is trivial on $H \subset G_{\overline{K}/K}$ of index $m$ with $G/H$ abelian.

So, $\xi : G/H \to M$ is a homomorphism. But since $G/H = G_{L/K}$ and

$L/K$ is unramified at $v$
$\Leftrightarrow \quad I_v = \{1\},$

$\xi : G/H \to M$ is unramified at $v$. ∎

**Proposition 4.7** *(See also [SIL1] Chapter VII Proposition 3.1 and Restatement in Chapter VIII Proposition 1.4) Let $K$ be a local field, i.e. complete with respect to a discrete valuation $v \in M_K^0$ with finite residue field $k$. Let $E/K$ be an elliptic curve with good reduction at $v$. Let $m \geq 1$ be prime to $\text{char}(k)$. Then the reduction*

$$E(K)[m] \to \tilde{E}_v(k_v)$$

*is an injection.*

*Proof.* Without proof here, see [SIL1] Chapter VII Proposition 3.1. ∎

**Proposition 4.8** *(See also [SIL1] Chapter X, Corollary 4.4) Let $E$ and $E'$ be elliptic curves over $K$ and let $\phi : E/K \to E'$ be an isogeny defined over $K$. Let $S$ be a finite set of places s.t. all infinite places $M_K^\infty$, the set of finite places $v \in K_K^0$ of bad reduction of $E$, and the set of finite places dividing $\deg(\phi)$ are contained in $S$.*
*Then one has*

$$S^\phi(E/K) \subset H^1(G_{\overline{K}/K}, E[\phi]; S).$$

*Proof.* Let $\xi \in S^\phi(E/K)$. Because of

$$S^\phi(E/K) \quad := \quad \ker\{H^1(G_{\overline{K}/K}, E[\phi]) \to \prod_{v \in M_K} WC(E/K_v)\}$$

$\xi$ is trivial in $WC(E/K_v) \cong H^1(G_v, E)$, i.e. there exists $P \in E(\overline{K}_v)$ s.t. $\xi_\sigma = P^\sigma - P$ for $\sigma \in G_v$. The action of $\sigma \in I_v$ is trivial on $\overline{K}_v$. So, reduction $E \to \tilde{E}_v$ gives

$$\widetilde{P^\sigma - P} = \tilde{P}^\sigma - \tilde{P} = O.$$

But $P^\sigma - P \in E[\phi] \subset E[m]$. According to 4.7 one gets that $E(K)[m]$ injects into $\tilde{E}_v$, since we have good reduction at $v \in S$. So, $\xi_\sigma = P^\sigma - P = O$ for all $\sigma \in I_v$, i.e. $\xi$ is unramified at $v \notin S$. ∎

Now we are able to prove the finiteness of the Selmer group:

*Proof (of Theorem 4.3).* According to Proposition 4.8 one has

$$S^\phi(E/K) \subset H^1(G_{\overline{K}/K}, E[\phi]; S),$$

for a set of places $M_K^\infty \subset S \subset M_K$. From 4.6 we know that $H^1(G_{\overline{K}/K}, E[\phi]; S)$ is finite, this gives the result. ∎

## 4.2 Finiteness of the Shafarevich-Tate Group

The finiteness of the Shafarevich-Tate group has not been proven in general, but only in very few cases.

**Conjecture 4.9** *Let $E/K$ be an elliptic curve.*

*The Shafarevich-Tate group $Ш(E/K)$ is finite.*

A general proof of the finiteness of $Ш$ does not exist yet, but there is the following proposition:

**Proposition 4.10** *(See also [SIL1] Chapter X, Theorem 4.14, Cassels-Tate-Pairing)*

*a) Let $E/K$ be an elliptic curve. There exists an alternating bilinear pairing*

$$\Gamma : Ш(E/K) \times Ш(E/K) \to \mathbb{Q}/\mathbb{Z},$$

*such that its kernel on each side equals the divisible elements in $Ш(E/K)$.*

*b) If $Ш(E/K)$ is finite then*

$$\#Ш(E/K) \text{ is a perfect square.}$$

*In this case $\dim_2(Ш(E/K))$ is even.*

*Proof.* b) follows from a). Without proof here. ∎

# 5 Descent via Two-Isogeny

## 5.1 Computing the Selmer Group

(See also [SIL1] Chapter X Proposition 4.9, Chapter X Proposition 4.8, and Chapter III Example 4.5)

Now we give a procedure to compute the Selmer group:

**Proposition 5.1** *Let*

$$E : y^2 = x^3 + ax^2 + bx$$
$$E' : y^2 = x^3 + a'x^2 + b'x$$

*be elliptic curves over $K$ s.t.*

$$a' = -2a$$
$$b' = a^2 - b.$$

*a) The map*

$$\phi : E \to E'$$
$$(x, y) \mapsto \left(\frac{y^2}{x^2}, \frac{y(b-x^2)}{x^2}\right)$$

*is an isogeny of degree 2 with kernel*

$$E[\phi] = \{O, (0,0)\}.$$

*b) For each $d \in K^*$ the hyperelliptic curve*

$$C_d : dw^2 = d^2 - 2adz^2 + (a^2 - 4b)z^4$$

*is a homogeneous space of $E$.*

*c) Define*

$$S := M_K^\infty \cup \{v \in M_K^0 \mid v(2) \neq 0\} \cup \{v \in M_K^0 \mid v(b) \neq 0\} \cup \{v \in M_K^0 \mid v(b') \neq 0\}$$

*and*

$$K(S, 2) := \{b \in K^*/(K^*)^2 : \mathrm{ord}_v(b) \equiv 0 \mod 2 \text{ for all } v \notin S\}.$$

*Then, there is an exact sequence*

$$0 \to E'(K)/\phi(E(K)) \xrightarrow{\delta} K(S, 2) \xrightarrow{\vartheta} \mathrm{WC}(E/K)[\phi],$$

*where*

$$\vartheta : K(S, 2) \to \mathrm{WC}(E/K)[\phi]$$
$$d \mapsto \{C_d/K\}$$

*and*

$$\delta : E'(K)/\phi(E(K)) \to K(S, 2)$$
$$\begin{cases} (X, Y) \mapsto X \\ O \mapsto 1 \\ (0, 0) \mapsto a^2 - 4b. \end{cases}$$

*d) The selmer group is given by*

$$S^{(\phi)}(E/K) \cong \{d \in K(S, 2) : C_d(K_v) \neq \emptyset \text{ for all } v \in S\}.$$

So, for each $d \in K(S,2)$ one tries to find out if $C_d(K_v)$ is empty or not.

Before we prove this proposition, we need to introduce the following generalized pairings:

**Proposition 5.2** *(Generalized Weil-Pairing, see also [SIL1] Exercise 3.15) Let $E/K$ and $E'/K$ be elliptic curves, and let $\phi : E \to E'$ be an isogeny of degree $m$. Assume that $\mathrm{char} K = 0$ or $m$ is prime to $\mathrm{char} K$.*
*There exists a bilinear, nondegenerate and Galois invariant pairing*

$$e_\phi : \ker\phi \times \ker\hat\phi \to \mu_m,$$

*where $\mu_m$ denotes the group of $m^{th}$ roots of unity.*

*Proof.* The pairing $e_\phi$ is a generalization of the standard Weil pairing on the m-torsion points $e_m : E[m] \times E[m] \to \mu_m$ constructed in [SIL1], Chapter III, Section 8. The generalized Weil pairing is constructed analogously, so I do not give the whole proof here, but the following ideas: We use the fact that a divisor $D = \sum_{P \in E} n_p(P) \in \mathrm{Div}(E)$ is principal if and only if its degree is 0 and $\sum_{P \in E}[n_P]P = O$ on $E$, see also [SIL1] Chapter III, Corollary 3.5. So, one can show that for $T \in \ker\hat\phi$, there are functions $f_T \in \overline{K}(E')$ and $g_T \in \overline{K}(E)$ such that

$$\mathrm{div}(f_T) = m(T) - m(O) \text{ and } f_T \circ \phi = g_T^m.$$

Let $S \in \ker\phi$. For any $X \in E$ one has

$$g_T(X+S)^m = f_T(\phi X + \phi S) = f_T(\phi X) = g_T(X)^m.$$

So, if we take $\frac{g_T(X+S)}{g_T(X)}$ as a function of $X$, then for each $X$ the value is in $\mu_m$. In particular, as a function of $X$, it is not surjective. Using the fact that a morphism of curves is either constant or surjective, one can define

$$e_\phi(S,T) = \frac{g_T(X+S)}{g_T(X)}.$$

One can show that $e_\phi$ is bilinear, nondegenerate and Galois invariant analogously as for the standard Weil pairing $e_m$ on the m-torsion points, see also [SIL1] Chapter III, Proposition 8.1. ∎

So, in the following assume $\mathrm{char} K = 0$ or $m$ is prime to $\mathrm{char} K$. In particular, the Weil pairing is defined and according to [SIL1] Chapter 3, Corollary 6.4(b) we can identify

$$E[m] \quad = \quad \mathbb{Z}/m\mathbb{Z} \times \mathbb{Z}/m\mathbb{Z}.$$

We need the following generalization of the Kummer pairing. This generalization is not given in [SIL1], but necessary for the next proposition.

**Proposition 5.3** *(Generalized Kummer Pairing) Let $E/K$ and $E'/K$ be elliptic curves and let $\phi : E \to E'$ be a nontrivial isogeny. Assume $E[\phi] \subset K(E)$. Then there exists a bilinear pairing*

$$\kappa : E'(K) \times G_{\overline{K}/K} \to E[\phi]$$

*Its kernel on the left is $\phi E(K)$ and its kernel on the right is $G_{\overline{K}/L}$ with*

$$L = K(\bigcup_{\substack{Q \in E(\overline{K}) \\ \phi Q \in E'(K)}} Q).$$

*We write $K(\phi^{-1}(E))$ for $L$.*

*Proof.* This is a generalization of the standard Kummer pairing

$$E(K) \times G_{\overline{K}/K} \to E[m]$$

for the isogeny $[m] : E \to E$.

For $P \in E'(K)$ one can choose any point $Q \in E$ s.t. $\phi Q = P$ and can show that the definition $\kappa(P, \sigma) = Q^\sigma - Q$ is well defined and leads to a bilinear pairing with the given properties. We will not do this here in detail, referring to [SIL1] Chapter VIII, Proposition 1.2. But I would like to note the following: Analogously as in the proof of Proposition 4.1 we take the exact sequence

$$0 \to E[\phi] \overset{\text{id}}{\to} E \overset{\phi}{\to} E' \to 0$$

We have already seen that this gives the fundamental short exact sequence

$$0 \longrightarrow E'(K)/\phi(E(K)) \overset{\delta}{\longrightarrow} H^1(G_{\overline{K}/K}, E[\phi]) \longrightarrow H^1(G_{\overline{K}/K}, E)[\phi] \longrightarrow 0$$

.

Under the assumtion that $E[\phi] \subset K(E)$, the Galois action is trivial on $E[\phi]$, so

$$H^1(G_{\overline{K}/K}, E[\phi]) = \text{Hom}(G_{\overline{K}/K}, E[\phi])$$

such that $\delta$ yields the generalized Kummer pairing. ∎

**Proposition 5.4** *(See also [SIL1] Exercise 10.1) Let $E/K$ and $E'/K$ be elliptic curves, and let $\phi : E \to E'$ be an isogeny of degree $m$. Assume $E[\phi] \subset E(K)$ and $E'[\hat{\phi}] \subset E'(K)$.*

*a) There exists a bilinear pairing*

$$b : E'(K)/\phi(E(K)) \times E'[\hat{\phi}] \to K(S, m)$$

*defined by*

$$e_\phi(\delta_\phi(P), T) = \delta_K(b(P, T)),$$

*where*

$$e_\phi : E[\phi] \times E'[\hat{\phi}] \to \mu_m,$$

*is the generalized Weil pairing given in Proposition 5.2, and*

$$\delta_\phi : E'(K) \to H^1(G_{\overline{K}/K}, E[\phi])$$
$$\delta_K : K^*/(K^*)^m \to H^1(G_{\overline{K}/K}, \mu_m)$$

*are the usual connecting homomorphisms.*

b) *The pairing b is nondegenerate on the left.*

c) *For $T \in E'[\hat{\phi}]$, there are $f_T \in K(E')$ and $g_T \in K(E)$ s.t.*

$$\operatorname{div}(f_T) = m(T) - m(O) \text{ and } f_T \circ \phi = g_T^m.$$

*Then, for all $P \neq O, T$, the value $b(P,T)$ can be expressed as follows*

$$b(P,T) = f_T(P) \mod (K^*)^m.$$

d) *If $m = 2$, i.e. $E'[\hat{\phi}] = \{O, T\}$, one has*

$$b(P,T) = x(P) - x(T) \mod (K^*)^2.$$

*Proof.* This Proposition is a generalization of [SIL1] Chapter X, Theorem 1.1 and can be proved analogously using the generalized Weil and generalized Kummer pairing instead. So, I would like to give the following ideas showing the difference between the original and the generalized proof:

a) 
- First, note that under the assumption $E[\phi] \subset E(K)$ and $E'[\hat{\phi}] \subset E'(K)$, the properties of the generalized Weil pairing yield $\mu_m \in K^*$. This is analogous to [SIL1] Chapter III, Corollary 8.1.1.
- The pairing is well-defined: We need to show that the image is contained in $K(S, m)$. The right kernel of the generalized Kummer pairing is $G_{\overline{K}/L}$, where $L/K$ is a Galois extension given in Proposition 5.3.

  Analogously as in [SIL1] Chapter VIII, Proposition 1.5 b, one can show that the extension $L/K$ is unramified outside $S$. The proof is analogous, but one needs to use the fact that $E[\phi]$ injects into $E[m]$.

  Let $P \in E(K)/\phi E(K)$, $T \in E'[\hat{\phi}]$ and $\beta \in \overline{K}^*$ such that $b(P,T)^m = b$. I want to prove that the extension $K(\beta)$ of $K$ is contained in $L$. The following argument is not given in the proof of [SIL1] Chapter X, Theorem 1.1 c: According to the definition we have

  $$e_\phi(\delta_E(P), T) = \delta_K(b(P,T)),$$

  i.e. for all $\sigma \in G_{\overline{K}/K}$:

  $$e_\phi(\delta_E(P)(\sigma), T) = \delta_K(b(P,T))(\sigma).$$

  According to Hilbert 90:

  $$\delta_K(b(P,T))(\sigma) = \frac{\beta^\sigma}{\beta} \in \mu_m.$$

  But, on the other hand, one has

  $$\delta_E(P)(\sigma) = \kappa(P)(\sigma) = Q^\sigma - Q \in E[\phi],$$

  where $Q \in E(\overline{K})$ s.t. $\phi Q = P$. So, one has

  $$\frac{e_\phi(Q,T)^\sigma}{e_\phi(Q,T)} = e_\phi(Q^\sigma - Q, T) = \frac{\beta^\sigma}{\beta} \in \mu_m.$$

But then, since $\mu_m \in K^*$:

$$K(\beta) \subset K(\bigcup_{\substack{Q' \in E(\overline{K}) \\ \phi Q' \in E'(K)}} Q') = L.$$

Now let $v \in M_K$ be a finite place and let $v(m) = 0$. Then, according to Definition 4.4, one can show

$$K(\beta)/K \text{ is unramified at } v \Leftrightarrow \text{ord}_v(\beta^m) \equiv 0 \mod m$$

which gives the result.

- The pairing is bilinear: Using the definition $e_\phi(\delta_\phi(P), T) = \delta_K(b(P,T))$, the bilinearity of $b$ follows from the following facts:
  - The generalized Weil pairing $e_\phi$ is bilinear according to Proposition 5.2.
  - The connecting homomorphism $\delta_\phi$ is bilinear because it can be identified with the generalized Kummer pairing $\kappa : E'(K) \times G_{\overline{K}/K} \to E[\phi]$ which is bilinear.
  - The connecting homomorphism $\delta_K$ is an isomorphism because of the Kummer duality $K^*/(K^*)^m \cong H^1(G_{\overline{K}/K}, \mu_m)$.

b) We need to show that if $b(P,T) = 1$ for all $T \in E[\phi]$, then $P$ is trivial. If $b(P,T) = 1$ for all $T \in E[\phi]$, then $e_\phi(\kappa(P,\sigma), T) = 1$ for all $T \in E[\phi]$ and $\sigma \in G_{\overline{K}/K}$. The left nondegeneracy of the generalized Weil pairing leads to $\kappa(P\sigma) = 0$ for all $\sigma$. Then, since the left kernel of the generalized Kummer pairing equals $\phi E(K)$, we get $P \in \phi E(K)$, so $P \in E'(K)/\phi E(K)$ is trivial.

c) This is analogous to [SIL1] Chapter X, Theorem 1.1d. Let $\beta = b(P,T)^{1/m}$ and take some $Q \in E(\overline{K})$ s.t. $P = \phi Q \in E'(\overline{K})$. Then

$$e_\phi(\delta_\phi(P), T) = \delta_K(b(P,T)).$$

Applying Hilbert 90, the isomorphism $\delta_K$ can be written as

$$\delta_K : K^*/(K^*)^m \to \text{Hom}(G_{\overline{K}/K}, \mu_m)$$
$$\delta_K(b)(\sigma) = \frac{\beta^\sigma}{\beta}.$$

This gives

$$e_\phi(\delta_\phi(P), T) = \beta^\sigma/\beta.$$

The left side of the equation equals

$$e_\phi(\delta_\phi(P), T) = e_\phi(Q^\sigma - Q, T)$$
$$= \frac{g_T(X + Q^\sigma - Q)}{g_T(X)} \text{ for any } X \in E.$$

By setting X=Q we get

$$e_\phi(\delta_\phi(P), T) = \frac{g_T(Q)^\sigma}{g_T(Q)}.$$

Finally, we have

$$\frac{g_T(Q)^\sigma}{g_T(Q)} = \frac{\beta^\sigma}{\beta}.$$

Using the fact that the map $\delta_K$ is an isomorphism, we conclude

$$g_T(Q)^m \equiv \beta^m \mod (K^*)^m,$$

and because of

$$f_T(P) = f_T \circ \phi(Q) = g_T(Q)^m$$

one gets

$$\begin{aligned}f_T(P) &\equiv \beta^m \mod (K^*)^m \\ &\equiv b(P, T) \mod (K^*)^m.\end{aligned}$$

d) It is sufficient to show this for the curves $E$ and $E'$ and the map $\phi$ given in Proposition 5.1. One can show that $\phi$ is an isogeny of degree 2 with kernel $((0,0), O)$. Setting $T = (0,0)$, we see that $f_T(x, y) = x$ satisfies

$$\operatorname{div}(f_T) = 2(T) - 2(O)$$

and

$$f_T \circ \phi = \frac{y^2}{x^2}$$

is a square. ∎

Now we can give the proof of Proposition 5.1, showing how to compute the Selmer group:

*Proof (Proof of Proposition 5.1).* a) The statement can be proved straightforwardly, without computation here.

b) This is already proved by Proposition 3.11.

c) For c) note the following: One can identify $E[\phi]$ with $\mu_2$ as $G_{\overline{K}/K}$-modules. This gives

$$K^*/(K^*)^2 \cong H^1(G_{\overline{K}/K}, E[\phi]) \text{ by Proposition 2.7 (Kummer)}$$
$$\text{and } K(S, 2) \cong H^1(G_{\overline{K}/K}, E[\phi]; S).$$

Using Proposition 4.8 we identify $S^{(\phi)} \subset K(S, 2)$. According to the definition of $S^{(\phi)}(E/K)$ and $\text{III}(E/K)$ and the the exact sequences given in Chapter 4 one has

$$0 \to E'(K)/\phi(E(K)) \to S^{(\phi)}(E/K) \to \text{III}(E/K)[\phi] \to 0.$$

This gives the exact sequence in c). Now we want to compute the connecting homomorphism $\delta$.

It is obvious that $\delta : O \mapsto 1$.

Let $T \in E'[\hat{\phi}]$. According to Proposition 5.4, we get

$$b(P,T) = x(P) - x(T) \mod (K^*)^2.$$

As we did in the proof of Proposition 3.11, we can assume $T = (0,0)$, such that the connecting homomorphism $\delta$ maps $P$ to $x(P) - x((0,0)) = x(P) - 0 = x(P)$, or with other notation $\delta : P = (X, Y) \mapsto X$.

Now I want to show that $\delta$ maps $(0,0)$ to $a^2 - 4b$. The Weierstrass equation of $E$ leads to

$$\begin{aligned} E : y^2 &= x^3 + ax^2 + bx \\ &= x(x^2 + ax + b) \\ &= (x - e_1)(x - e_2)(x - e_3) \end{aligned}$$

with

$$\begin{aligned} e_1 &= 0 \\ e_{2,3} &\in K(\sqrt{a^2 - 4b}) \end{aligned}$$

The following argument is not given in [SIL1]. According to [SIL1] Chapter 3, Corollary 6.4(b) we can identify the 2-torsion group with the Klein four-group:

$$E[2] = \mathbb{Z}/2\mathbb{Z} \times \mathbb{Z}/2\mathbb{Z}.$$

In particular, let $E[2]$ consist of the points $O$ and $(0,0)$ and two other points $s_1$ and $s_2$. Then, in the Klein four-group, one has $s_1 s_2 = (0,0)$ and $\delta$ maps $(0,0)$ to $a^2 - 4b$.

d) Analogously to Theorem 3.9, one can see that an element $C/K$ of $WC(E/K)$ is trivial in $WC(E/K_v)$ if $C$ has any point $p_0 \in K_v$. Then, the cocycle

$$\begin{aligned} \xi : G_v &\to E \\ \sigma &\mapsto p_0^\sigma - p_0 \end{aligned}$$

is trivial. Since $S^\phi \subset H^1(G_{\overline{K}/K}, E[\phi]; S)$ by Proposition 4.8, we get the identification in d). ∎

**Proposition 5.5** *One has*

$$P \in C_d(K) \quad \Rightarrow \quad \delta(\psi(P)) \equiv d \mod (K^*)^2,$$

*where*

$$\begin{aligned} \psi : C_d &\to E' \\ (z, w) &\mapsto \left( \frac{d}{z^2}, -\frac{dw}{z^3} \right). \end{aligned}$$

*Proof.* Note that the proof of Proposition 3.11 gives a natural map

$$\begin{aligned} \Theta : E &\to C_d \\ (x, y) &\mapsto \left( \frac{\sqrt{d}x}{y}, \sqrt{d}\left(x - \frac{b}{x}\right)\left(\frac{x}{y}\right)^2 \right). \end{aligned}$$

This is an isomorphism and its inverse is ([SIL1] Chp. X, Remark 3.7)

$$\Theta^{-1}: C_d \to E$$
$$(z,w) \mapsto \left(\frac{\sqrt{d}w - az^2 + d}{2z^2}, \frac{dw - a\sqrt{d}z^2 + d\sqrt{d}}{2z^3}\right).$$

Using the isogeny $\phi: E \to E'$, $(x,y) \mapsto \left(\frac{y^2}{x^2}, \frac{y(b-x^2)}{x^2}\right)$ one gets

$$\psi := \Theta^{-1} \circ \phi : C_d \to E'$$
$$(z,w) \mapsto \left(\frac{d}{z^2}, -\frac{dw}{z^3}\right).$$

Since $\delta$ takes the first coordinate and $\frac{d}{z^2} \equiv d \mod \mathbb{Q}^2$ one gets the result. ∎

## 5.2 Computing the weak Mordell-Weil group

Having computed the Selmer group, we can give the following procedure to compute the weak Mordell-Weil group:

Let $E$ and $E'$ be elliptic curves over $K$, let $\phi : E \to E'$ be an isogeny of degree $m$, and let $\hat{\phi} : E' \to E$ be its dual isogeny, in particular $\hat{\phi} \circ \phi = [m]$.

a) (See also [SIL1] Chapter X, Theorem 4.2a) The sequences

$$0 \to E'(K)/\phi(E(K)) \to S^{(\phi)}(E/K) \to \text{III}(E/K)[\phi] \to 0$$
$$0 \to E(K)/\hat{\phi}(E'(K)) \to S^{(\hat{\phi})}(E'/K) \to \text{III}(E'/K)[\hat{\phi}] \to 0$$

are exact.

b) (See also [SIL1] Chapter 10, Section 4.7) The sequences

$$0 \to \frac{E'(K)[\hat{\phi}]}{\phi(E(K)[m])} \to \frac{E'(K)}{\phi(E(\mathbb{Q}))} \xrightarrow{\hat{\phi}} \frac{E(K)}{mE(K)} \to \frac{E(K)}{\hat{\phi}(E'(K))} \to 0$$

and

$$0 \to \frac{E(K)[\phi]}{\hat{\phi}(E'(K)[m])} \to \frac{E(K)}{\hat{\phi}(E'(K))} \xrightarrow{\phi} \frac{E'(K)}{mE'(\mathbb{Q})} \to \frac{E'(K)}{\phi(E(K))} \to 0.$$

are exact.

This follows directly from the definitions of the Selmer and Shafarevich-Tate groups and Proposition 4.1 and from the facts that $m = \deg(\phi)$ and $\hat{\phi}$ is dual to $\phi$, i.e. $\hat{\phi} \circ \phi = [m]$.

For example, using these properties, the first exact sequence in b) is obvious since

$$\text{im}\left(\frac{E'(K)[\hat{\phi}]}{\phi(E(K))[m]} \to \frac{E'(K)}{\phi(E(K))}\right) = \frac{E'(K)[\hat{\phi}]}{\phi(E(K))} = \ker\left(\frac{E'(K)}{\phi(E(K))} \xrightarrow{\hat{\phi}} \frac{E(K)}{mE(K)}\right)$$

$$\text{im}\left(\frac{E'(K)}{\phi(E(K))} \xrightarrow{\hat{\phi}} \frac{E(K)}{mE(K)}\right) = \frac{\hat{\phi}(E'(K))}{m(E(K))} = \ker\left(\frac{E(K)}{mE(K)} \to \frac{E(K)}{\hat{\phi}(E'(K))}\right).$$

The other exact sequence follows by changing the roles of $\phi$ and $\hat{\phi}$.

So, if one was able to compute the Selmer Groups and if the corresponding Shafarevich-Tate groups are trivial, one can determine $E(K)/\phi(E'(K))$ and $E'(K)/\phi(E(K))$ by a). Using the exact sequence in b), one can determine the weak Mordell-Weil Group.

## 5.3 Example: $E : y^2 = x^3 + 6x^2 + x$

(See [SIL1] Chapter X, Exercise 10.14(a), analogously to [SIL1] Chapter X, Example 4.10) The discriminant is $\Delta = 2^9$. So, we have $S = \{2, \infty\}$, and $\mathbb{Q}(S, 2) = \{\pm 1, \pm 2\}$.
For $a := 6$, $b := 1$, we determine $-2 \cdot a = -12$ and $a^2 - 4 \cdot b = 36 - 4 = 32$, s.t.

$$E' = x^3 - 12 \cdot x^2 + 32 \cdot x.$$

For each value of $d \in \mathbb{Q}(S, 2)$, the homogeneous space for $E$ is given by

$$C_d : dw^2 = d^2 - 2adz^2 + (a^2 - 4b)z^4,$$

so we have

$$C_d : dw^2 = d^2 - 12dz^2 + 32z^4.$$

And for $E'$ with $a' := -12$ and $b' := 32$, we have $-2a' = 24$ and $a'^2 - 4b' = 16$: So, the homogeneous space for $E'$ is given by the equation

$$C'_d : dw^2 = d^2 + 24dz^2 + 16z^4.$$

The *Selmer group* $S^{(\phi)}(E/\mathbb{Q})$ is a subgroup of $\mathbb{Q}(S, 2)$ and represented by those $d \in \mathbb{Q}(S, 2)$, s.t. $C_d$ is solvable over $\mathbb{R}$ and every field of $p$-adic numbers $\mathbb{R}$ and $\mathbb{Q}_p$ for all primes $p$.

It is $E'(\mathbb{Q})/\phi(E(\mathbb{Q}))$ a subgroup of $S^{(\phi)}(E/\mathbb{Q})$ and represented by those $d \in \mathbb{Q}(S, 2)$, s.t. $C_d$ is solvable over $\mathbb{Q}$.

Finally, the *Shafarevich-Tate group* $\text{III}(E/\mathbb{Q})$ is the quotient of $S^{(\phi)}(E/\mathbb{Q})$ by $E'(\mathbb{Q})/\phi(E(\mathbb{Q}))$. So, the Shafarevich-Tate group is represented by those homogeneous spaces, that are solvable over $\mathbb{R}$ and over the $p$-adic numbers $\mathbb{Q}_p$ for all primes $p$, but that have no rational solution.

For each $d \in \mathbb{Q}(S, 2)$, we try to find a rational point of $C_d$ or a prime number $p$ s.t. $C_d(\mathbb{Q}_p) = \emptyset$.

In the first case, we determine $d$ (representing $C_d$) to be contained in $E'(\mathbb{Q})/\phi(E(\mathbb{Q}))$.

In the second case, we determine $d$ to be *not* contained in the Selmer group at all.

If none of this both cases holds, this cannot be determined in a finite amount of time. Then, $d$ represents an element of the Shafarevich-Tate group.

- The case $d = 1$ does not need to be regarded, because the neutral element 1 is always contained in $S^{(\phi)}(E/\mathbb{Q})$ and in each subgroup. Remember also $\delta(\mathcal{O}) \mapsto 1$.

- If $d = 2$, the homogeneous space for $E$ is given by

$$C_2 : 2w^2 = 2^2 - 12 \cdot 2z^2 + 32z^4.$$

Dividing by 2, we get $w^2 = 2 - 12z^2 + 16z^4$. Substitution $Z = 2z$ gives $w^2 = 2 - 3Z^2 + Z^4$. This has a solution $(w, Z) = (0, 1)$, because $0^2 = 2 - 3 + 1$.
Further, we know $\delta(0, 0) = 32 \equiv 2 \ (mod \ \mathbb{Q}^{*2})$.
So, we have

$$2 \in E'(\mathbb{Q})/\phi(E(\mathbb{Q})) \subset S^{(\phi)}(E/\mathbb{Q}).$$

- If $d = -2$, the homogeneous space for $E$ is given by

$$C_{-2} : -2w^2 = 2^2 - 12 \cdot (-2)z^2 + 32z^4.$$

Dividing by $-2$, we get $w^2 = -2 - 12z^2 - 16z^4$. Regarding this equation modulo $2^3$, we have $w^2 = 6 + 4z^2$. This has no integer solution modulo 8: For $w^2$ and $z^2$ there are just the possibilities $0^2 \equiv 0$, $1^2 \equiv 1$, $2^2 \equiv 4$, $3^2 \equiv 1$, $4^2 \equiv 0$, $5^2 \equiv 4$, $6^2 \equiv 4$, $7^2 \equiv 1 \ mod \ 8$. So, the set of squares of integers is given by $\{0, 1, 4\}$. Because of $4z^2 \in \{0, 4\}$ and $6 + 4z^2 \in \{2, 6\}$, there is no integer solution modulo 8. This leads to $C_{-2}(\mathbb{Q}_2) = \emptyset$, and

$$-2 \notin S^{(\phi)}(E/\mathbb{Q}).$$

- If $d = -1$: Since $S^{(\Phi)}(E/\mathbb{Q})$ is a subgroup of $\mathbb{Q}(S, 2)$, we also have

$$-1 \notin S^{(\phi)}(E/\mathbb{Q}).$$

because otherwise it would be $2 \cdot (-1) = -2 \in S^{(\Phi)}(E/\mathbb{Q})$ a contradiction.

Finally, this leads to $S^{(\Phi)}(E/\mathbb{Q}) = E'(\mathbb{Q})/\phi(E(\mathbb{Q})) = \{1, 2\}$ and $\text{Ш}(E/\mathbb{Q}) = \{1\}$.

Now we do the same for the homogeneous space of $E'$:

- If $d = \pm 2$, the homogeneous space for $E'$ is given by

$$C'_{\pm 2} : \pm 2w^2 = 2^2 + 24 \cdot (\pm 2)z^2 + 16z^4.$$

Reducing modulo $2^3$, we would get $\pm 2w^2 = 2^2$. But this is a contradiction, since $ord_2(2w^2)$ is odd, and $ord_2(2^2)$ is even.
So, it is $C'_2(\mathbb{Q}_2) = C'_2(\mathbb{Q}_{-2}) = \emptyset$ and

$$\pm 2 \notin S^{(\hat{\phi})}(E'/\mathbb{Q}).$$

- If $d = -1$, the homogeneous space for $E'$ is given by

$$C'_{-1} : -w^2 = 2^2 - 24z^2 + 16z^4.$$

This has a solution $(z, w) = (1, 2)$, because $-2^2 = 4 - 8 = 2^2 - 24 \cdot 1^2 + 16 \cdot 1^4$.
So, we have

$$-1 \in E(\mathbb{Q})/\hat{\phi}(E'(\mathbb{Q})) \subset S^{(\hat{\phi})}(E'/\mathbb{Q}).$$

This gives $S^{(\hat{\Phi})}(E'/\mathbb{Q}) = E(\mathbb{Q})/\hat{\phi}(E'(\mathbb{Q})) = \{\pm 1\}$ and $\text{III}(E'/\mathbb{Q}) = \{1\}$.

So, now we know that

$$E'(\mathbb{Q})/\phi(E(\mathbb{Q})) \cong \mathbb{Z}/2\mathbb{Z} \text{ and } E(\mathbb{Q})/\hat{\phi}(E'(\mathbb{Q})) \cong \mathbb{Z}/2\mathbb{Z}.$$

Using the exact sequence in [SIL1] Chapter 10, Section 4.7, we have

$$0 \to \frac{E'(\mathbb{Q})[\hat{\phi}]}{\phi(E(\mathbb{Q})[2])} \to \frac{E'(\mathbb{Q})}{\phi(E(\mathbb{Q}))} \xrightarrow{\hat{\phi}} \frac{E(\mathbb{Q})}{2E(\mathbb{Q})} \to \frac{E(\mathbb{Q})}{\hat{\phi}(E'(\mathbb{Q}))} \to 0$$

and

$$0 \to \frac{E(\mathbb{Q})[\phi]}{\hat{\phi}(E'(\mathbb{Q})[2])} \to \frac{E(\mathbb{Q})}{\hat{\phi}(E'(\mathbb{Q}))} \xrightarrow{\phi} \frac{E'(\mathbb{Q})}{2E'(\mathbb{Q})} \to \frac{E'(\mathbb{Q})}{\phi(E(\mathbb{Q}))} \to 0.$$

This leads to

$$E(\mathbb{Q})/2(E(\mathbb{Q})) \cong \mathbb{Z}/2\mathbb{Z} \text{ and } E'(\mathbb{Q})/2(E'(\mathbb{Q})) \cong \mathbb{Z}/2\mathbb{Z}.$$

# Part II
# Elliptic Curves over $K = \mathbb{Q}$ with $j$-Invariant $1728$

## 6 The Curve $E_D : Y^2 = X^3 + DX$

### 6.1 The $j$-Invariant $1728$

**Corollary 6.1** *Let $E$ be an elliptic curve over $\mathbb{Q}$. Then*

$$j(E) = 1728 \Leftrightarrow \exists D \in \mathbb{Q}^* \text{ s.t. } E_D : y^2 = x^3 + Dx \text{ is isomorphic to } E.$$

*Proof.* "$\Rightarrow$": This follows directly from Proposition 2.8.
"$\Leftarrow$": For $E_D : y^2 = x^3 + Ax + B$ with $A = D$ and $B = 0$ we have

$$\begin{aligned} j(E_D) &= -1728 \frac{(4A)^3}{\Delta} \\ &= -1728 \frac{(4A)^3}{-16(4A^3 + 27B^2)} \quad | B = 0 \\ &= 1728. \end{aligned}$$

Two elliptic curves have the same $j$-invariant if and only if they are isomorphic over $\overline{K}$, see also [SIL1] Chapter III, Proposition 1.4.b. It follows $j(E) = j(E_D) = 1728$. ∎

### 6.2 Computing $E_{D,tors}(\mathbb{Q})$

Let $E : y^2 = f(x) = x^3 + Dx$ be an elliptic curve with $D \in \mathbb{Q}^*/(\mathbb{Q}^*)^4$.

In the following, we identify $D$ with a fourth-power-free integer, such that $D$ is uniquely defined by $E$.

**Proposition 6.2** *The torsion group is given by*

$$E_{D,\text{tors}}(\mathbb{Q}) \cong \begin{cases} \mathbb{Z}/4\mathbb{Z} & \text{if } D = 4, \\ \mathbb{Z}/2\mathbb{Z} \times \mathbb{Z}/2\mathbb{Z} & \text{if } -D \text{ is a perfect square,} \\ \mathbb{Z}/2\mathbb{Z} & \text{otherwise.} \end{cases}$$

*Proof.* Let $q > 2$ be an odd prime and $\mathbb{F}_q$ a finite field of order $q$.

i) Assume $q \equiv 3 \mod 4$. We prove $\#\widetilde{E}_D(\mathbb{F}_q) = q+1$ (See also [SIL1] Exercise 10.17):

Since $2 \mid (q-1)$ there is a unique subgroup of $F_q^*$ of order $\frac{q-1}{2}$ consisting of the square roots in $F_q^*$, i.e. the equation $v = u^2$ has $\frac{q-1}{2}$ solutions in $\mathbb{F}_q^* \times \mathbb{F}_q^*$. Since $q \equiv 3 \mod 4$ the map $u^2 \mapsto u^4$ is an automorphism of $(\mathbb{F}_q^*)^2$. And $v \mapsto v^2$ is a 2-1 map in $\mathbb{F}_q$. So, the equation $v^2 = u^4 + 16D$ has exactly $2 \cdot \frac{q-1}{2} = q - 1$ solutions in $\mathbb{F}_q^* \times \mathbb{F}_q^*$.

Let $\widetilde{E_D}$ be the elliptic curve reduced modulo $q$. The map

$$\begin{aligned} \phi : C &\to \widetilde{E_q} \\ (u, v) &\mapsto \left( \frac{u^2}{4}, \frac{uv}{8} \right) \end{aligned}$$

is a well-defined isogeny: Let $(u, v) \in C$, then $D = \frac{v^2 - u^4}{16}$ and

$$\left(\frac{u^2}{4}\right)^3 + D\left(\frac{u^2}{4}\right) = \frac{u^6}{64} + \frac{v^2 - u^4}{16}\frac{u^2}{4}$$
$$= \frac{u^6}{64} + \frac{v^2 u^2 - u^6}{64}$$
$$= \frac{u^2 v^2}{64}$$
$$= \left(\frac{uv}{8}\right)^2.$$

So, $\phi(u, v) \in E$. Taking into account that $(0, 0)$ and $O$ are solutions for $E$, one gets $\#E(\mathbb{F}_p) = (q - 1) + 2 = q + 1$.

ii) According to Proposition 4.7 $E_{D_{\text{tors}}}(\mathbb{Q})$ injects into $\widetilde{E}_D(\mathbb{F}_q)$. For all but finitely many primes $q \equiv 3 \mod 4$ one has $\#E_{D_{\text{tors}}}(\mathbb{Q}) \mid (p + 1)$ such that $\#E_{D_{\text{tors}}}(\mathbb{Q}) \mid 4$, hence $\#E_{D,\text{tors}}(\mathbb{Q}) \in \{1, 2, 4\}$. So, one has

$$E_{D,\text{tors}}(\mathbb{Q}) \in \{\mathbb{Z}/2\mathbb{Z}, (\mathbb{Z}/2\mathbb{Z})^2, \mathbb{Z}/4\mathbb{Z}\}.$$

iii) We have

$$E[2] \subset E(\mathbb{Q}) \quad \Leftrightarrow \quad f(x) = x^3 + Dx \text{ factors completely over } \mathbb{Q}$$
$$\Leftrightarrow \quad -D \text{ is a perfect square.}$$

Since $E[2] \cong \mathbb{Z}/2\mathbb{Z} \times \mathbb{Z}/2\mathbb{Z}$ one gets

$$E_{D,\text{tors}}(\mathbb{Q}) = \mathbb{Z}/2\mathbb{Z} \times \mathbb{Z}/2\mathbb{Z} \quad \Leftrightarrow \quad -D \text{ is a perfect square.}$$

iv) One has

$$E_{D,\text{tors}}(\mathbb{Q}) = \mathbb{Z}/4\mathbb{Z} \quad \Leftrightarrow \quad E(\mathbb{Q}) \text{ has a point of order } 4$$
$$\Leftrightarrow \quad (0, 0) \in 2E(\mathbb{Q})$$

Let $P := (D^{1/2}, (4D^3)^{1/4})$. According to the duplication formula one has $[2]P = (0, 0)$. Since we assumed $D$ to be a fourth-power-free integer, one gets

$$(0, 0) \in 2E(\mathbb{Q}) \quad \Leftrightarrow \quad D = 4,$$

hence

$$E_{D,\text{tors}}(\mathbb{Q}) = \mathbb{Z}/4\mathbb{Z} \quad \Leftrightarrow \quad D = 4 \qquad \blacksquare$$

## 6.3 An upper bound for the rank of $E_D(\mathbb{Q})$

**Proposition 6.3** *Let $\dim_2$ denote the dimension of an $\mathbb{F}_2$ vectorspace. Then, have the following relation between $E(\mathbb{Q})/2E(\mathbb{Q})$ and $E'(\mathbb{Q})/\phi E(\mathbb{Q})$ and $E(\mathbb{Q})/\hat{\phi} E'(\mathbb{Q})$:*

$$\dim_2 \frac{E(\mathbb{Q})}{2E(\mathbb{Q})} = \dim_2 \frac{E'(\mathbb{Q})}{\phi E(\mathbb{Q})} + \dim_2 \frac{E(\mathbb{Q})}{\hat{\phi}(E'(\mathbb{Q}))} + \dim_2 \phi(E(\mathbb{Q}))[2] - \dim_2 E'(\mathbb{Q})[\hat{\phi}].$$

*Proof.* From Section 5.2 b) we have the exact sequence

$$0 \to \frac{E'(\mathbb{Q})[\hat{\phi}]}{\phi(E(\mathbb{Q})[2])} \to \frac{E'(\mathbb{Q})}{\phi(E(\mathbb{Q}))} \xrightarrow{\hat{\phi}} \frac{E(\mathbb{Q})}{2E(\mathbb{Q})} \to \frac{E(\mathbb{Q})}{\hat{\phi}(E'(\mathbb{Q}))} \to 0.$$

Then

$$\dim_2 \frac{E(\mathbb{Q})}{2E(\mathbb{Q})} + \dim_2 \frac{E'(\mathbb{Q})[\hat{\phi}]}{\phi(E(\mathbb{Q}))[2]} = \dim_2 \frac{E'(\mathbb{Q})}{\phi(E(\mathbb{Q}))} + \dim_2 \frac{E(\mathbb{Q})}{\hat{\phi}(E'(\mathbb{Q}))}$$

Since $\dim_2 \frac{E'(\mathbb{Q})[\hat{\phi}]}{\phi(E(\mathbb{Q}))[2]} = \dim_2 E'(\mathbb{Q})[\hat{\phi}] - \dim_2 \phi(E(\mathbb{Q}))[2]$ one gets the result. ∎

**Proposition 6.4** *(See also [SIL1] Chapter X, Proposition 6.2(b)) Let $\nu(2D)$ be the number of distinct primes dividing $2D$. Then there is an upper bound for the rank of $E(\mathbb{Q})$:*

$$\mathrm{rank} E(\mathbb{Q}) < 2\nu(2D).$$

*Proof.* This is an application of the Two-Descent described in Section 5: One has

$$E' : Y^2 = X^3 - 4DX$$

and an isogeny

$$\begin{aligned} \phi : E &\to E' \\ \phi(x,y) &= \left( \frac{y^2}{x^2}, \frac{y\left(y(D-x^2)\right)}{x^2} \right) \end{aligned}$$

One easily sees that $\phi$ and its dual isogeny $\hat{\phi} : E' \to E$ have degree 2. Further,

$$\begin{aligned} S &= \{p \in M_{\mathbb{Q}} : p \text{ divides } 2D\} \cup \{\infty\} \\ \mathbb{Q}(S,2) &= \{b \in \mathbb{Q}^*/(\mathbb{Q}^*)^2 : \mathrm{ord}_p \equiv 0 \mod 2 \text{ for all } p \notin S\}. \end{aligned}$$

For each $d \in \mathbb{Q}(S,2)$ one gets homogeneous spaces $C_d/\mathbb{Q} \in \mathrm{WC}(E/\mathbb{Q})$ and $C'_d/\mathbb{Q} \in \mathrm{WC}(E'/\mathbb{Q})$ with equations

$$\begin{aligned} C_d : \quad dw^2 &= d^2 - 4Dz^4 \\ C'_d : \quad dW^2 &= d^2 + 16DZ^4, \end{aligned}$$

and substitution $Z \to Z/2$ leads to $C'_d : dW^2 = d^2 + DZ^4$. According to Proposition 6.3 one has

$$\dim_2 \frac{E(\mathbb{Q})}{2E(\mathbb{Q})} = \dim_2 \frac{E'(\mathbb{Q})}{\phi E(\mathbb{Q})} + \dim_2 \frac{E(\mathbb{Q})}{\hat{\phi}(E'(\mathbb{Q}))} + \dim_2 \phi(E(\mathbb{Q}))[2] - \dim_2 E'(\mathbb{Q})[\hat{\phi}].$$

Because $E'(\mathbb{Q})/\phi E(\mathbb{Q})$ and $E(\mathbb{Q})/\hat{\phi}E'(\mathbb{Q})$ inject into the group $\mathbb{Q}(S,2)$ one has

$$\dim_2 \frac{E(\mathbb{Q})}{2E(\mathbb{Q})} \leq \#\mathbb{Q}(S,2) + \#\mathbb{Q}(S,2) + \dim_2 \phi(E(\mathbb{Q}))[2] - \dim_2 E'(\mathbb{Q})[\hat{\phi}].$$

Since the group $\mathbb{Q}(S,2)$ is generated by $-1$ and the primes dividing $2D$, we get $\#\mathbb{Q}(S,2) = 2(1+\nu(2D))$ and

$$\dim_2 E(\mathbb{Q})/2E(\mathbb{Q}) \leq 2\left(1+\nu(2D)\right) - \dim_2 E'(\mathbb{Q})[\hat{\phi}] + \dim_2 \phi(E(\mathbb{Q})[2])$$

Because of $E'(\mathbb{Q})[\hat{\phi}] \cong \mathbb{Z}/2\mathbb{Z}$, one has (*)

$$\dim_2 E(\mathbb{Q})/2E(\mathbb{Q}) \leq 2\nu(2D) + \dim_2 \phi(E(\mathbb{Q})[2])$$

From the Mordell-Weil theorem we know

$$E(\mathbb{Q}) \cong E(\mathbb{Q})_{\text{tors}} \times \mathbb{Z}^{\text{rank}E(\mathbb{Q})}$$

- If $E(\mathbb{Q})[2] \cong \mathbb{Z}/2\mathbb{Z}$: One has

$$\dim_2 E(\mathbb{Q})/2E(\mathbb{Q}) = 1 + \text{rank}E(\mathbb{Q}).$$

Since $\phi$ has kernel $\{(0,0), O\}$, we have

$$\phi(E(\mathbb{Q})[2]) \cong 0.$$

and $\dim_2 \phi(E(\mathbb{Q})[2]) = 1$. Substituting this results into (*) one gets

$$\dim_2 E(\mathbb{Q})/2E(\mathbb{Q}) \leq 2\nu(2D) + \dim_2 \phi(E(\mathbb{Q})[2])$$
$$\Leftrightarrow \quad 1 + \text{rank}E(\mathbb{Q}) \leq 2\nu(2D) + 1$$
$$\Leftrightarrow \quad \text{rank}E(\mathbb{Q}) \leq 2\nu(2D).$$

- If $E(\mathbb{Q})[2] \cong \mathbb{Z}/2\mathbb{Z} \times \mathbb{Z}/2\mathbb{Z}$: One has

$$\dim_2 E(\mathbb{Q})/2E(\mathbb{Q}) = 2 + \text{rank}E(\mathbb{Q}).$$

and

$$\phi(E(\mathbb{Q})[2]) \cong \mathbb{Z}/2\mathbb{Z},$$

and $\dim_2 \phi(E(\mathbb{Q})[2]) = 2$. Then according to (*)

$$\dim_2 E(\mathbb{Q})/2E(\mathbb{Q}) \leq 2\nu(2D) + \dim_2 \phi(E(\mathbb{Q})[2])$$
$$\Leftrightarrow \quad 2 + \text{rank}E(\mathbb{Q}) \leq 2\nu(2D) + 2$$
$$\Leftrightarrow \quad \text{rank}E(\mathbb{Q}) \leq 2\nu(2D).$$

So, in both cases we have $\text{rank}E(\mathbb{Q}) \leq 2\nu(2D)$. One can increase this upper bound by 1, if one takes into account the following argument: If $d < 0$ then at least one homogeneous space $C_d$ or $C'_d$ has no real points, because one side of the equation is negative, whereas the other side is positive. So, $C_d(\mathbb{R}) = \emptyset$ or $C'_d(\mathbb{R}) = \emptyset$, such that

$$\text{rank}E(\mathbb{Q}) \leq 2\nu(2D) - 1 \qquad \blacksquare$$

# 7 The Case $D = p$ for an Odd Prime $p$

Let $p \neq 2$ be prime, and let

$$E_p : y^2 = x^3 + px$$

be the Weierstrass equation for the elliptic curve $E_p$. Using the notations from above let $\phi : E_p \to E_p'$ be an isogeny of degree 2 with kernel $E_p[\phi] = \{O, (0,0)\}$.

**Corollary 7.1**

$$E_{p,\text{tors}}(\mathbb{Q}) \;\;=\;\; \mathbb{Z}/2\mathbb{Z}$$

*Proof.* This follows directly from Proposition 6.2. ∎

**Definition 7.2** Let $a \in \mathbb{Z}$ and let $p$ be prime. Let

$$\left(\frac{a}{p}\right) := \begin{cases} 1 & \text{if } a \text{ is a quadratic residue of } p \\ -1 & \text{if } a \text{ is a quadratic non-residue of } p \\ 0 & \text{if } a \equiv 0 \mod p \end{cases}$$

denote the *Legendre* symbol.

**Proposition 7.3** *(Quadratic Reciprocity Law) Let $p, q > 2$ be prime. Then*

$$\left(\frac{p}{q}\right) = (-1)^{\frac{p-1}{2}} \left(\frac{p}{q}\right).$$

*Further, we have*

$$\left(\frac{-1}{p}\right) = (-1)^{\frac{p-1}{2}}$$

*so, $-1$ is a quadratic residue modulo $p$ if and only if $p \equiv 1 \mod 4$, and*

$$\left(\frac{2}{p}\right) = (-1)^{\frac{p^2-1}{8}}$$

*so, 2 is a quadratic residue modulo $p$ if and only if $p \equiv \pm 1 \mod 8$.*

*Proof.* This is a fundamental fact of algebraic number theory, see also [SCPI] Paragraph 3, Theorem 5 and Theorem 6, [SCHM] Section 2.2 or [SCPI] Section 8.3. ∎

**Proposition 7.4** *(Euler) Let $p > 2$ be prime and $a \in \mathbb{Z}$ s.t. $p \nmid a$. Then*

$$\left(\frac{a}{p}\right) = a^{\frac{p-1}{2}} \mod p.$$

*Proof.* This is a fundamental fact of algebraic number theory, see also [SCPI] Section 8.3 ∎

**Proposition 7.5**

$$S^{(\hat\phi)}(E'_p/\mathbb{Q}) = \{1, p\}$$

$$S^{(\phi)}(E_p/\mathbb{Q}) = \begin{cases} \{1, -p\} & \text{if } p \equiv 7, 11 \mod 16 \\ \{1, -p, -1, -p\} & \text{if } p \equiv 5 \mod 16 \\ \{1, -p, 2, -2p\} & \text{if } p \equiv 15 \mod 16 \\ \{1, -p, -2, 2p\} & \text{if } p \equiv 3 \mod 16 \\ \{\pm 1, \pm p\} & \text{if } p \equiv 13 \mod 16 \\ \{\pm 1, \pm p, \pm 2, \pm 2p\} & \text{if } p \equiv 1, 9 \mod 16. \end{cases}$$

**Corollary 7.6** *(See also [SIL1] Chapter X Proposition 6.2b)*

$$S^{(\hat\phi)}(E'_p/\mathbb{Q}) \cong \mathbb{Z}/2\mathbb{Z}$$

$$S^{(\phi)}(E_p/\mathbb{Q}) \cong \begin{cases} \mathbb{Z}/2\mathbb{Z} & \text{if } p \equiv 7, 11 \mod 16 \\ (\mathbb{Z}/2\mathbb{Z})^2 & \text{if } p \equiv 3, 5, 13, 15 \mod 16 \\ (\mathbb{Z}/2\mathbb{Z})^3 & \text{if } p \equiv 1, 9 \mod 16. \end{cases}$$

*Proof (Proposition 7.5).* Since we execute a 2-descent, we identify squares in $\mathbb{Q}(S, 2)$. And since 2 does not divide $p$, we get the following representatives

$$\mathbb{Q}(S, 2) = \{\pm 1, \pm 2, \pm p, \pm 2p\}.$$

Let $\delta : E'(\mathbb{Q})/\phi(E(\mathbb{Q})) \to K(S, 2)$ denote the map from Section 5.1 (b). Set $a := 0$ and $b := p$ from the Weierstrass equation of $E$, then

$$\delta((0, 0)) = a^2 - 4b = -4p \equiv -p.$$

So we find $p$ to be a representative in $S^{(\phi)}(E'_p/\mathbb{Q})$. We have

$$E'_p : y^2 = x^3 - 4p$$

and set $a' := 0$ and $b' := -4p$. Let $\delta'$ denote the map from Section 5.1 (b) with roles of $E$ and $E'$ reversed. One gets

$$\delta'((0, 0)) = a'^2 - 4b' = -4(-4p) \equiv p.$$

So, $p$ is a representative in $S^{(\hat\phi)}(E_p/\mathbb{Q})$. Consider the homogeneous spaces

$$D_d : \quad dw^2 = d^2 - 4pz^4$$
$$C'_d : \quad dW^2 = d^2 + pZ^4.$$

1. Computing $S^{(\hat\phi)}(E_p/\mathbb{Q})$:

   (a) We show that $S^{(\hat\phi)}(E_p/\mathbb{Q})$ does not contain even numbers: Take the homogeneous space

   $$C'_2 : 2W^2 = 4 + pZ^4$$

and reduce modulo 2:

$$2W^2 \equiv 4 + pZ^4 \mod 2$$
$$\Leftrightarrow \quad 0 \equiv 0 + pZ^4 \mod 2 \qquad | \, p \equiv 1 \mod 2$$
$$\Leftrightarrow \quad Z^4 \equiv \phantom{0} \mod 2$$
$$\Leftrightarrow \quad Z \equiv 0 \mod 2,$$

then $Z \equiv 0 \mod 4$. Reduction modulo 4 gives

$$2W^2 \equiv 4 + pZ^4 \mod 4$$
$$\Leftrightarrow \quad 2W^2 \equiv 0 \mod 4$$
$$\Leftrightarrow \quad W^2 \equiv 0 \mod 4$$
$$\Leftrightarrow \quad W \equiv 0 \mod 2.$$

But then we find that $C'_2$ has no solutions in $\mathbb{Q}_8$:

$$2W^2 \equiv 4 + pZ^4 \mod 8$$
$$\Leftrightarrow \quad 0 \equiv 4 + 0 \mod 8.$$

(b) Let $d < 0$. $S^{(\hat{\phi})}(E_p/\mathbb{Q})$ does not contain negative numbers, since

$$C'_d : \underbrace{dW^2}_{<0} = \underbrace{d^2}_{>0} + \underbrace{pZ^4}_{>0}.$$

So, $S^{(\hat{\phi})}(E_p/\mathbb{Q}) = \{1, p\} \cong \mathbb{Z}/2\mathbb{Z}$.

2. Computing $S^{(\phi)}(E_p/\mathbb{Q})$: Since $p$ is prime and 2 does not divide $p$, we can use the fact that

$$d \in S^{(\phi)}(E_p/\mathbb{Q}) \Leftrightarrow C_d \text{ has no solutions in } \mathbb{Q}_p \text{ and } \mathbb{Q}_2.$$

(a) $d = -1$: Take the homogeneous space

$$C_{-1} : -w^2 = 1 - 4pz^4.$$

i. Reducing modulo $p$:

$$-w^2 \equiv 1 - 4pz^4 \mod p$$
$$\Leftrightarrow \quad w^2 \equiv -1 \mod p.$$

It is a fact that the $p$-adic numbers contain all $(p-1)^{\text{th}}$ roots of unity. Let $\zeta$ be a generator for the $(p-1)^{\text{th}}$ roots of unity. Since $p$ is odd, 2 divides $p-1$. If $p \equiv 1 \mod 4$ then one can take $w = \zeta^{\frac{p-1}{4}}$ s.t. $w^2 = \zeta^{\frac{p-1}{2}} \equiv -1 \mod p$. So, there is a solution for $C_{-1}$ in $\mathbb{Q}_p$ if and only if $p \equiv 1 \mod 4$. This is also given in Proposition 7.3.

ii. Reducing modulo 2:
My following argument differs from the one given in [SIL1]:
$(w, z) = (1, 0)$ is a solution for

$$-w^2 = 1 - 4pz^4 \mod 2$$

since $-1^2 = 1 - 0 \mod 2 \Leftrightarrow -1 \equiv 1 \mod 2$. In particular, this still holds in case of $p \equiv 1 \mod 4$.

(b) $d = 2$: Take the homogeneous space
$$C_2 : 2w^2 = 4 - 4pz^4.$$
One can divide by 2 s.t.
$$C_2 : w^2 = 2 - 2pz^4.$$

i. Reducing modulo $p$: This gives
$$w^2 \equiv 2 \mod p.$$
According to Proposition 7.3, 2 is a quadratic residue modulo $p$ if and only if $p \equiv 1, 7 \mod 8$

ii. Reducing modulo 2: This gives
$$w^2 \equiv 0 \mod 2$$
so $w \equiv 0 \mod 2$ and we can substitute $(z, w) \to (Z, 2W)$, s.t.
$$C_2 : 2W^2 = 1 - pZ^4.$$
We can assume $p \equiv 1, 7 \mod 8$, so we only consider the following cases for $p$, and can compute:

| $p$ | Solution $(Z, W)$ of reduced equation |
|---|---|
| 1 mod 8 | no solutions |
| 7 mod 16 | no solutions |
| 15 mod 32 | $(3, 1)$ |
| 31 mod 32 | $(1, 1)$ |

Finally, we get that $-2$ is contained in the Selmer group if and only if $p \equiv 1, 9, 15 \mod 16$

(c) $d = -2$: Take the homogeneous space
$$C_{-2} : -2w^2 = 4 - 4pz^4.$$
One can divide by 2 s.t.
$$C_{-2} : -w^2 = 2 - 2pz^4.$$

i. Reducing modulo $p$: This gives
$$w^2 \equiv -2 \mod p.$$
Using Proposition 7.3 and Proposition 7.4 one gets
$$\left(\frac{-2}{p}\right) = \left(\frac{2}{p}\right)\left(\frac{-1}{p}\right)$$
$$= \begin{cases} 1 \text{ if } \left(\frac{2}{p}\right) = \left(\frac{-1}{p}\right) \\ -1 \text{ else.} \end{cases}$$
$$= \begin{cases} 1 & \text{if } (p \equiv 1, 7 \mod 8 \text{ and } p \equiv 1 \mod 4) \\ & \text{or } (p \equiv 3, 5 \mod 8 \text{ and } p \equiv 3 \mod 4) \\ -1 & \text{else.} \end{cases}$$
$$= \begin{cases} 1 \text{ if } p \equiv 1, 3 \mod 8 \\ -1 \text{ else.} \end{cases}$$

ii. Reducing modulo 2: This gives

$$w^2 \equiv 0 \mod 2$$

so $w \equiv 0 \mod 2$ and we can substitute $(z, w) \to (Z, 2W)$, s.t.

$$C_{-2} : -2W^2 = 1 - pZ^4.$$

We can assume $p \equiv 1, 3 \mod 8$, so we consider the following cases for $p$, and can compute:

| $p$ | Solution $(Z, W)$ of reduced equation |
|---|---|
| 1 mod 32 | $(1, 0)$ |
| 3 mod 32 | $(3, 11)$ |
| 9 mod 32 | $(1, 2)$ |
| 11 mod 16 | no solution |
| 17 mod 32 | $(3, 0)$ |
| 19 mod 32 | $(1, 3)$ |
| 25 mod 32 | $(3, 2)$ |

Finally, we get that $-2$ is contained in the Selmer group if and only if $p \equiv 1, 3, 9 \mod 16$ ∎

**Proposition 7.7**

$$\mathrm{rank} E_p(\mathbb{Q}) + \dim_2 \mathrm{III}(E_p/\mathbb{Q})[2] = \begin{cases} 0 & \text{if } p \equiv 7, 11 \mod 16 \\ 1 & \text{if } p \equiv 3, 5, 13, 15 \mod 16 \\ 2 & \text{if } p \equiv 1, 9 \mod 16. \end{cases}$$

*Proof.* In the following we write $E$ for $E_p$ and $E'$ for $E'_p$. From Section 5.2 b) we have the exact sequence

$$0 \to \frac{E'(\mathbb{Q})[\hat{\phi}]}{\phi(E(\mathbb{Q})[2])} \to \frac{E'(\mathbb{Q})}{\phi(E(\mathbb{Q}))} \xrightarrow{\hat{\phi}} \frac{E(\mathbb{Q})}{2E(\mathbb{Q})} \to \frac{E(\mathbb{Q})}{\hat{\phi}(E'(\mathbb{Q}))} \to 0.$$

This yields (*)

$$\dim_2 E'(\mathbb{Q})[\hat{\phi}]/\phi(E(\mathbb{Q})[2]) + \dim_2 E(\mathbb{Q})/2E(\mathbb{Q})$$
$$= \dim_2 E'(\mathbb{Q})/\phi(E(\mathbb{Q})) + \dim_2 E(\mathbb{Q})/\hat{\phi}(E(\mathbb{Q})).$$

Analogously, using the exact sequences from Section 5.2 a)

$$0 \to E'(\mathbb{Q})/\phi(E(\mathbb{Q})) \to S^{(\phi)}(E/\mathbb{Q}) \to \mathrm{III}(E/\mathbb{Q})[\phi] \to 0$$
$$0 \to E(\mathbb{Q})/\hat{\phi}(E'(\mathbb{Q})) \to S^{(\hat{\phi})}(E'/\mathbb{Q}) \to \mathrm{III}(E'/\mathbb{Q})[\hat{\phi}] \to 0$$

one has (**)

$$\dim_2 E'(\mathbb{Q})/\phi(E(\mathbb{Q})) = \dim_2(S^{(\phi)}(E/\mathbb{Q})) - \dim_2(\mathrm{III}(E/\mathbb{Q})[\phi])$$
$$\dim_2 E(\mathbb{Q})/\hat{\phi}(E(\mathbb{Q})) = \dim_2(S^{(\hat{\phi})}(E'/\mathbb{Q})) - \dim_2(\mathrm{III}(E'/\mathbb{Q})[\hat{\phi}]).$$

So, (*) and (**) give (***)

$$\dim_2 E'(\mathbb{Q})[\hat{\phi}]/\phi(E(\mathbb{Q})[2]) + \dim_2 E(\mathbb{Q})/2E(\mathbb{Q})$$
$$= \dim_2(S^{(\phi)}(E/\mathbb{Q})) - \dim_2(\text{III}(E/\mathbb{Q})[\phi])$$
$$+ \dim_2(S^{(\hat{\phi})}(E'/\mathbb{Q})) - \dim_2(\text{III}(E'/\mathbb{Q})[\hat{\phi}]).$$

i) From the Mordell-Weil theorem we know

$$E(\mathbb{Q}) \cong E(\mathbb{Q})_{\text{tors}} \times \mathbb{Z}^{\text{rank} E(\mathbb{Q})}.$$

This gives

$$E(\mathbb{Q})/2E(\mathbb{Q}) \cong (\mathbb{Z}/2\mathbb{Z})^{1+\text{rank} E(\mathbb{Q})}.$$

ii) Because $\phi$ has kernel $\{(0,0), O\}$, one has $\phi(E(\mathbb{Q})[2]) \cong 0$. Since $\hat{\phi}$ is its dual isogeny with the same kernel, one gets

$$E'(\mathbb{Q})[\hat{\phi}]/\phi(E(\mathbb{Q})[2]) \cong \mathbb{Z}/2\mathbb{Z}.$$

iii) We know

$$E(\mathbb{Q})/\hat{\phi}(E'(\mathbb{Q})) \subseteq S^{(\hat{\phi})}(E'/\mathbb{Q}).$$

Because of $S^{(\hat{\phi})}(E'/\mathbb{Q}) = \mathbb{Z}/2\mathbb{Z}$ by Proposition 7.6, one gets

$$E(\mathbb{Q})/\hat{\phi}(E'(\mathbb{Q})) \cong S^{(\hat{\phi})}(E'/\mathbb{Q}),$$

which we prove by contradiction: Otherwise the left side would need to be trivial, i.e. $\hat{\phi}(E'(\mathbb{Q})) = E(\mathbb{Q})$. Using $\phi\hat{\phi} = [2]$ we would have $\phi(E\mathbb{Q}) = 2E(\mathbb{Q})$, which is not true.

iv) The equality $E(\mathbb{Q})/\hat{\phi}(E'(\mathbb{Q})) \cong S^{(\hat{\phi})}(E'/\mathbb{Q})$ and (**) lead to

$$\text{III}(E'/\mathbb{Q})[\hat{\phi}] \cong 0.$$

The sequence

$$0 \to \text{III}(E/\mathbb{Q})[\phi] \to \text{III}(E/\mathbb{Q})[2] \xrightarrow{\phi} \underbrace{\text{III}(E'/\mathbb{Q})[\hat{\phi}]}_{=0}$$

is exact, s.t.

$$\dim_2 \text{III}(E/\mathbb{Q})[2] = \dim_2 \text{III}(E/\mathbb{Q})[\phi].$$

Substitution in (***) yields

$$1 + 1 + \text{rank} E(\mathbb{Q}) = \dim_2(S^{(\phi)}(E/\mathbb{Q})) + 2 - \dim_2 \text{III}(E/\mathbb{Q})[2]$$
$$\Leftrightarrow \text{rank} E(\mathbb{Q}) + \dim_2 \text{III}(E/\mathbb{Q})[2] = \dim_2(S^{(\phi)}(E/\mathbb{Q})).$$

The values for $S^{(\phi)}(E/\mathbb{Q})$ are also given in Proposition 7.6. This gives the result. ∎

**Corollary 7.8** *a) If $p \equiv 7, 11 \mod 16$:*

$$\mathrm{rank} E_p(\mathbb{Q}) \ = \ 0.$$

*b) If $\mathrm{III}(E_p/\mathbb{Q})$ is finite and $p \equiv 3, 5, 13, 15 \mod 16$:*

$$\mathrm{rank} E_p(\mathbb{Q}) \ = \ 1.$$

*c) If $\mathrm{III}(E_p/\mathbb{Q})$ is finite and $p \equiv 1, 9 \mod 16$:*

$$\mathrm{rank} E_p(\mathbb{Q}) \ \in \ \{0, 2\}.$$

*Proof.* This follows directly from Proposition 7.7 and Proposition 4.10 b). ∎

In the following chapter, we show that if $p \equiv 1 \mod 8$, and if 2 is not a quartic residue for $p$, one has $\mathrm{rank} E_p(\mathbb{Q}) = 0$.

# 8 Prime Numbers $p \equiv 1 \mod 8$

**Corollary 8.1** *Let $p$ be an odd prime. Consider the elliptic curve*

$$E_p : y^2 = x^3 + px.$$

*Then $\mathbb{Q}(S, 2) = \{\pm 1, \pm 2, \pm p, \pm 2p\}$ and*

$$S^\phi(E_p/\mathbb{Q}) = \mathbb{Q}(S, 2) \ \Leftrightarrow \ p \equiv 1 \mod 8.$$

*Proof.* This follows directly from Proposition 7.5. ∎

**Proposition 8.2** *a) (See also [SIL1] Chapter X, Remark 6.5.1)*
*Let $p \equiv 1 \mod 8$ be prime. Then there are integers $A, B \in \mathbb{Z}$ s.t.*

$$p \ = \ A^2 + B^2$$

*and*

$$AB \ \equiv \ 0 \mod 4.$$

*b) (See also [SIL1] Chapter X, Case II in the proof of Proposition 6.4.)*
*Let $p \equiv 1 \mod 4$ be prime. Then there are integers $A, B \in \mathbb{Z}$ s.t.*

$$p \ = \ A^2 + B^2$$

*and*

$$A \ \equiv \ 1 \mod 2$$
$$B \ \equiv \ 0 \mod 2.$$

*Proof.* This is a fundamental fact of algebraic number theory, referring to Fermat's theorem on sums of two squares. For more details see [DCOX].

**Proposition 8.3** *(Quartic Reciprocity) Let $p \equiv 1 \mod 8$ be prime. Let $A, B \in \mathbb{Z}$ be integers from Proposition 8.2 s.t.*

$$p = A^2 + B^2$$
$$AB \equiv 0 \mod 4.$$

*For any $a \in \mathbb{Z}$ let*

$$\left(\frac{a}{p}\right)_4 := \begin{cases} 1 & \text{if } a \text{ is a quartic residue of } p \\ -1 & \text{if } a \text{ is a quartic non-residue of } p \\ 0 & \text{if } a \equiv 0 \mod p \end{cases}$$

*Then*

$$\left(\frac{2}{p}\right)_4 = (-1)^{AB/4}.$$

*In particular*

$$2 \text{ is a quartic residue modulo } p \quad \Leftrightarrow \quad AB \equiv 0 \mod 8.$$

*Proof.* For a fancy proof see [SIL1] Chapter X, Proposition 6.6. ∎

**Proposition 8.4** *(See also [SIL1] Chapter X, Proposition 6.5) Let $p \equiv 1 \mod 8$, and assume that 2 is not a quartic residue for $p$. Consider the elliptic curve*

$$E_p : y^2 = x^3 + px.$$

*The homogeneous spaces*

$$\begin{aligned} C_{-1}: & \quad w^2 = -1 + 4pz^4 \\ C_{-2}: & \quad w^2 = -2 + 2pz^4 \\ C_2: & \quad w^2 = 2 - 2pz^4 \end{aligned}$$

*have no points defined over $\mathbb{Q}$.*

*Proof.* We prove this by contradiction:

1. $C_{-1} : w^2 = -1 + 4pz^4$: Let $(z, w)$ be a rational solution of $C_{-1}$. Now we prove that we can write $(z, w)$ in the form $(\frac{r}{2t}, \frac{s}{2t^2})$ with $r, s, t \in \mathbb{Z}$ s.t. $\gcd(r, t) = 1$:

    - Let $z = \frac{a}{b} \in \mathbb{Q}$ with $a, b \in \mathbb{Z}$ and $\gcd(a, b) = 1$.
        - In case of $2 \mid b$: Then, there exists $r, t \in \mathbb{Z}$ s.t. $b = 2t$ and $a = r$ leads to $z = \frac{r}{2t}$, $\gct(r, t) = 1$.
        - In case of $2 \nmid b$: Then, write $z = \frac{2a}{2b}$, set $r = 2a$ and $t = b$ s.t. $z = \frac{r}{2t}$, $\gct(r, t) = 1$.
    - Since we assume $(z, w)$ to be a rational solution of $C_{-1}$, one has

    $$\begin{aligned} & w^2 = -1 + 4pz^4 & \left| z = \frac{r}{2t} \right. \\ \Leftrightarrow \quad & w^2 = -1 + 4p\frac{r^4}{16t^4} & \left| \cdot (16t^4) \right. \\ \Leftrightarrow \quad & 16t^4 w^2 = -16t^4 + 4pr^4 & \left| \cdot \frac{1}{4} \right. \\ \Leftrightarrow \quad & 4t^4 w^2 = -4t^4 + pr^4 \end{aligned}$$

We substitute $s := 2t^2 w$. Now we show $s \in \mathbb{Z}$:

Assume $s = \frac{x}{y}$ with $x, y \in \mathbb{Z}$, $\gcd(x,y) = 1$. If $y \neq 1$, then there exists a prime $n \geq 2$ s.t. $n \mid y$. But $s^2 = pr^4 - 4t^4$ would lead to $x^2 = y^2(pr^4 - 4t^4)$, hence $n \mid x$, which is a contradiction to $\gcd(x,y) = 1$.

Now one has

$$w^2 = -1 + 4pz^4$$
$$\Leftrightarrow s^2 = pr^4 - 4t^4.$$

By Proposition 8.2 b) there are integers $A, B \in \mathbb{Z}$ s.t. $A \equiv 1 \mod 4$, $B \equiv 0 \mod 4$ and $p = A^2 + B^2$.
Then

$$p\left(Br^2 + 2t^2\right)^2 + A^2 s^2 \qquad \mid p = A^2 + B^2$$
$$= p\left(Br^2 + 2t^2\right)^2 + \left(p - B^2\right) s^2 \qquad \mid s^2 = pr^4 - 4t^4$$
$$= p\left(Br^2 + 2t^2\right)^2 + \left(p - B^2\right)\left(pr^4 - 4t^4\right)$$
$$= p\left(B^2 r^4 + \underbrace{4Br^2 t^2}_{2(2Bt^2)pr^2} + 4t^4\right) + \underbrace{p^2 r^4}_{(pr^2)^2} - 4pt^4 - pB^2 r^4 + \underbrace{4B^2 t^4}_{(2Bt^2)^2}$$
$$= (pr^2 + 2Bt^2)^2$$

and

$$\left(pr^2 + 2Bt^2 + As\right)\left(pr^2 + 2Bt^2 - As\right) = p\left(Br^2 + 2t^2\right)^2.$$

By Lemma 8.5, one can check that $gcd(pr^2 + 2Bt^2 + As, pr^2 + 2Bt^2 - As)$ is a square or twice a square. Then, there exist $u, v \in \mathbb{Z}$ s.t.

$$\begin{cases} pr^2 + 2Bt^2 \pm As = pu^2 \\ pr^2 + 2Bt^2 \mp As = v^2 \\ Br^2 + 2t^2 = uv \end{cases} \text{ or } \begin{cases} pr^2 + 2Bt^2 \pm As = 2pu^2 \\ pr^2 + 2Bt^2 \mp As = 2v^2 \\ Br^2 + 2t^2 = 2uv \end{cases}.$$

One can eliminate $s$ and gets

$$\begin{cases} 2pr^2 + 4Bt^2 = pu^2 + v^2 \\ Br^2 + 2t^2 = uv \end{cases} \text{ or } \begin{cases} pr^2 + 2Bt^2 = pu^2 + v^2 \\ Br^2 + 2t^2 = 2uv \end{cases}.$$

Since $p \equiv 1 \mod 8$ we can apply Proposition 8.3 as follows: Because of $\left(\frac{2}{p}\right)_4 = 1$, one gets $AB \equiv 0 \mod 8$. Because $A \equiv 1 \mod 2$, this gives us $B \equiv 0 \mod 8$.

Reducing the systems of equations modulo 8, we get

$$\begin{cases} 2pr^2 \equiv pu^2 + v^2 \mod 8 \\ 2t^2 = uv \mod 8 \end{cases} \text{ or } \begin{cases} pr^2 \equiv pu^2 + v^2 \mod 8 \\ t^2 = uv \mod 8 \end{cases}.$$

- In the first case, the second equation $2t^2 = uv \mod 8$ would imply that either $u$ or $v$ is even. Remember that we assume $p$ to be odd. If one of $u$ or $v$ is even and the other one is odd, then $pu^2 + v^2$ is odd, which is a contradiction to the first equation.

- In the second case, the second equation $t^2 = uv \mod 8$ implies that $u$ or $v$ are both even or both odd. If they are both even, then the first equation implies that $r$ is even, and the second equation implies that $t$ is even. This is a contradiction to $\gcd(r,t) = 1$.
  If $u$ and $v$ both are odd, then $pu^2 + v^2$ is even. The first equation implies $4 \equiv pu^2 + v^2 \mod 8$, which is a contradiction to the assumption that $u$ and $v$ are odd.

2. $C_{\pm 2} : w^2 = \pm 2 \mp 2pz^4$: Assume $(z, w)$ is a rational solution of $C_{\pm 2}$. Since $z \in \mathbb{Q}$ we can write $z = \frac{r}{t}$ with $r, t \in \mathbb{Z}$ s.t. $\gcd(r, t) = 1$. We get

$$w^2 = \pm 2 + \mp 2pz^4 \qquad \mid z = \frac{r}{t}$$

$$\Leftrightarrow \quad w^2 = \pm 2 \mp 2p\frac{r^4}{t^4} \qquad \mid \cdot \frac{t^4}{2}$$

$$\Leftrightarrow \quad 2t^4 w^2 = \pm 1 \mp pr^4 \qquad \mid \cdot \frac{1}{4}$$

$$\Leftrightarrow \quad 2(\frac{wt^2}{2})^2 = \pm 1 \mp pr^4$$

We substitute $s = \frac{wt^2}{2}$. Now we show $s \in \mathbb{Z}$:
Assume $s = \frac{x}{y}$ with $x, y \in \mathbb{Z}$, $\gcd(x, y) = 1$. If $y \neq 1$, then there exists a prime $n \geq 2$ s.t. $n \mid y$. But $\pm 2s^2 = t^4 - pr^4$ would lead to $\pm 2x^2 = y^2(t^4 - pr^4)$, hence $n \mid x$, which is a contradiction to $\gcd(x, y) = 1$.
So, $(z, w) = (\frac{r}{t}, \frac{2s}{t^2})$ with $r, s, t \in \mathbb{Z}$ s.t. $\gcd(r, s, t) = 1$ and $\pm 2s^2 = t^4 - pr^4$.
Then for any odd prime $q$ dividing $s$ one has

$$\left(\frac{p}{q}\right) = 1 \qquad \mid \text{Proposition 7.3}$$

$$\Leftrightarrow \quad \left(\frac{q}{p}\right) = 1 \qquad \mid \left(\frac{2}{p}\right) = 1$$

$$\Leftrightarrow \quad \left(\frac{s}{p}\right) = 1$$

$$\Leftrightarrow \quad \left(\frac{s^2}{p}\right)_4 = 1.$$

But then, $2s^2 = t^4 - pr^4$ implies $\left(\frac{\pm 2}{p}\right)_4 = 1$, which is a contradiction to the fact $\left(\frac{-1}{p}\right)_4 = 1$ for $p \equiv 1 \mod 8$ and the assumption $\left(\frac{2}{p}\right)_4 \neq 1$. ∎

I proved the following Lemma, which is not given in [SIL1]:

**Lemma 8.5** *(Referring to [SIL1], Proof of Proposition X.6.5, Case II) Let*

- *$p$ prime, $p \equiv 1 \mod 8$*

- $r, t \in \mathbb{Z}$: $\gcd(r, t) = 1$
- $s \in \mathbb{Z}$: $s^2 p = r^4 - 4t^4$
- $A, B \in \mathbb{Z} : A \equiv 1 \mod 4,\ B \equiv 0 \mod 4,\ p = A^2 + B^2$

and assume the following factorization
$$(pr^2 + 2Bt^2 + As)(pr^2 + 2Bt^2 - As) = p(Br^2 + 2t^2)^2.$$

Then $g := \gcd(pr^2 + 2Bt^2 + As, pr^2 + 2Bt^2 - As)$ is a square or twice a square.

*Proof.* One has
$$g \mid \big((pr^2 + 2Bt^2 + As) \pm (pr^2 + 2Bt^2 - As)\big)$$
$$\Rightarrow\ g \mid 2(pr^2 + 2Bt^2) \text{ and } g \mid 2As$$
$$\Rightarrow\ g \mid 2\gcd(pr^2 + 2Bt^2, As).$$

Now set
$$g^* := \gcd(pr^2 + 2Bt^2, As). \tag{1}$$

Because of $\gcd(pr^2 + 2Bt^2, As) \mid \gcd(pr^2 + 2Bt^2 + As, pr^2 + 2Bt^2 - As)$, one gets
$$g^* \mid g \mid 2g^*,$$
and it is sufficient to show that
$$g^* := \gcd(pr^2 + 2Bt^2, As) \text{ is a square or twice a square.}$$

Because of
$$\begin{aligned} pr^2 + 2Bt^2 &= (A^2 + B^2)r^2 + 2Bt^2 \\ &= A^2 r^2 + B^2 r^2 + 2Bt^2 \\ &= A^2 r^2 + B(Br^2 + 2t^2), \end{aligned}$$

equation (1) gives
$$g^* = \gcd\left(A^2 r^2 + B(Br^2 + 2t^2), As\right). \tag{2}$$

The factorization $(pr^2 + 2Bt^2 + As)(pr^2 + 2Bt^2 - As) = p(Br^2 + 2t^2)^2$ implies $(g^*)^2 \mid (Br^2 + 2t^2)^2$ since $p$ is prime. This leads to
$$g^* \mid (Br^2 + 2t^2). \tag{3}$$

Then, (2) and (3) imply
$$g^* \mid A^2 r^2. \tag{4}$$

One has
$$s^2 = pr^4 - 4t^4 \tag{5}$$
$$= (A^2 + B^2)r^4 - 4t^4$$
$$= A^2 r^4 + B^2 r^4 - 4t^4$$
$$= A^2 r^4 + (Br^2 - 2t^2)(Br^2 + 2t^2). \tag{6}$$

Let $q^k \mid g^*$, $q$ prime, $k \in \mathbb{N}$, $k$ odd. Then, show $q = 2$ or $q^{k+1} \mid g^*$:

i) Begin with
$$(4) \Rightarrow q^k \mid (Ar)^2 \Rightarrow q^{\frac{k+1}{2}} \mid Ar \Rightarrow q^{k+1} \mid A^2 r^2.$$

ii) One has
$$\left.\begin{array}{l}(3) \Rightarrow q^k \mid (Br^2 + 2t^2) \\ (4) \Rightarrow q^k \mid A^2 r^2 \Rightarrow q^k \mid A^2 r^4\end{array}\right\} \stackrel{(6)}{\Rightarrow} q^k \mid s^2 \Rightarrow q^{\frac{k+1}{2}} \mid s \Rightarrow q^{k+1} \mid s^2.$$

Further, according to (i)
$$q^{k+1} \mid A^2 r^2 \Rightarrow q^{k+1} \mid A^2 r^4.$$

So, using equation (6), $q^{k+1} \mid s^2$ and $q^{k+1} \mid A^2 r^4$ lead to
$$q^{k+1} \mid (Br^2 - 2t^2)(Br^2 + 2t^2).$$

Because of $(3) \Rightarrow q^k \mid (Br^2 + 2t^2)$, it follows
$$q^{k+1} \mid (Br^2 + 2t^2) \quad \text{or} \quad q \mid (Br^2 - 2t^2).$$

If $q^{k+1} \mid (Br^2 + 2t^2)$ one has $q^{k+1} \mid B(Br^2 + 2t^2)$.
If $q \mid (Br^2 - 2t^2)$ one has
$$q \mid \big((Br^2 + 2t^2) \pm (Br^2 - 2t^2)\big) \Rightarrow q \mid 2Br^2 \text{ and } q \mid 4t^2.$$

Because of $\gcd(r, t) = 1$ and $q$ prime, it follows $q = 2$ or $q \mid B$.
If $q \mid B$, then $q \cdot q^k \mid B(Br^2 + 2t^2) \Rightarrow q^{k+1} \mid B(Br^2 + 2t^2)$.

iii) One has
$$\left.\begin{array}{l}(ii) \Rightarrow q^{\frac{k+1}{2}} \mid s \\ \text{If } q \nmid r : (i) \Rightarrow q^{\frac{k+1}{2}} \mid Ar \Rightarrow q^{\frac{k+1}{2}} \mid A\end{array}\right\} \Rightarrow q^{k+1} \mid As$$

If $q \mid r$ : Then $q \mid pr^4$, such that (5) and $q \mid s^2$ give $q \mid 4t^2$. Because of $\gcd(r, t) = 1$, it follows $q \nmid t$ and $q = 2$.

According to (2) one has
$$g^* = \gcd(\underbrace{A^2 r^2 + B(Br^2 + 2t^2)}_{(i)\phantom{xx}(ii)}, \underbrace{As}_{(iii)}),$$

such that
$$\left.\begin{array}{ll}(i) & \Rightarrow \phantom{q = 2 \text{ or }} q^{k+1} \mid A^2 r^2 \\ (ii) & \Rightarrow q = 2 \text{ or } q^{k+1} \mid B(Br^2 + 2t^2) \\ (iii) & \Rightarrow q = 2 \text{ or } q^{k+1} \mid As\end{array}\right\} \Rightarrow q = 2 \text{ or } q^{k+1} \mid g^*.$$

∎

**Corollary 8.6** *Let $p \equiv 1 \mod 8$ and let 2 be no quartic residue for $p$.*

$$\mathrm{rank} E_p(\mathbb{Q}) = 0 \text{ and } \mathrm{III}(E_p(\mathbb{Q}))[2] = (\mathbb{Z}/2\mathbb{Z})^2.$$

*Proof.* From Corollary 8.1 we know that $S^\phi(E_p/\mathbb{Q}) = \mathbb{Q}(S,2) = \{\pm 1, \pm 2, \pm p, \pm 2p\}$.

$\mathrm{III}(E_p(\mathbb{Q}))[2] = (\mathbb{Z}/2\mathbb{Z})^2$ follows from the facts that $C_{-1}, C_2$ and $C_{-2}$ have no rational points and $-p$ is contained in the image of $(0,0)$. Then, Proposition 7.7 leads to $\mathrm{rank} E_p(\mathbb{Q}) = 0$. ∎

So, if $p \equiv 1 \mod 8$ and 2 is not a quartic residue for $p$, then there do exist homogeneous spaces containing $\mathbb{Q}_v$-rational points for every completion $\mathbb{Q}_v$ of $\mathbb{Q}$ on the one hand, but no rational points on the other hand. This is called a failure of the *Hasse principle*. In this case, our descent procedure does not terminate.

**Corollary 8.7** *A failure of the Hasse principle does occur.*

*Proof.* This follows from Proposition 8.1 and Proposition 8.4. ∎

# 9 Prime Number $p = 17$

### 9.0.1 Computation in `SageMath`

I wrote the following script `example.sage` in SageMath:

```
print '____________________'
print 'Defining the curve:'
E = EllipticCurve([0,0,0,17,0])
print 'E = '+`E`
d=E.discriminant()
print 'discriminant = ' + `d` + ' = ' + `d.factor()`
print 'E minimal? ' + `E.is_minimal()`
print '____________________'
print 'Cremona Data:'
try: print E.cremona_label()
except (RuntimeError) as msg: print('RuntimeError during call of E.
    cremona_label(): ' + `msg`)
#print E.cremona_label()
c = E.conductor()
print 'conductor = ' + `c` + ' = ' + `c.factor()`
print '____________________'
print 'Local Data:'
print E.local_data()
print '____________________'
print 'Torsion:'
T = E.torsion_subgroup()
print T
print 'Torsion Points: ' + `T.points()`
print '____________________'
print 'Rank:'
print 'rank = '+`E.rank()`
print '____________________'
print 'two_descent_by_two_isogeny:'
from sage.schemes.elliptic_curves.descent_two_isogeny import
   two_descent_by_two_isogeny
print two_descent_by_two_isogeny(E, verbosity=3)
print '____________________'
print 'simon_two_descent:'
print E.simon_two_descent(verbose=1)
```

We get the following result by typing `load("example.sage")`:

```
sage: load("example.sage")
____________________
Defining the elliptic curve:
E = EllipticCurve([0,0,0,17,0])
E = Elliptic Curve defined by y^2 = x^3 + 17*x over Rational Field
discriminant = -314432 = -1 * 2^6 * 17^3
E minimal? True
____________________
Cremona Data:
RuntimeError during call of E.cremona_label(): RuntimeError('Cremona label
   not known for Elliptic Curve defined by y^2 = x^3 + 17*x over Rational
   Field.',)
conductor = 18496 = 2^6 * 17^2
____________________
Local Data:
[Local data at Principal ideal (2) of Integer Ring:
Reduction type: bad additive
Local minimal model: Elliptic Curve defined by y^2 = x^3 + 17*x over Rational
    Field
Minimal discriminant valuation: 6
Conductor exponent: 6
Kodaira Symbol: II
Tamagawa Number: 1, Local data at Principal ideal (17) of Integer Ring:
Reduction type: bad additive
Local minimal model: Elliptic Curve defined by y^2 = x^3 + 17*x over Rational
    Field
Minimal discriminant valuation: 3
Conductor exponent: 2
Kodaira Symbol: III
Tamagawa Number: 2]
____________________
```

```
28 Torsion:
29 Torsion Subgroup isomorphic to Z/2 associated to the Elliptic Curve defined
      by y^2 = x^3 + 17*x over Rational Field
30 Torsion Points: [(0 : 1 : 0), (0 : 0 : 1)]
31 --------------------
32 Rank:
33 rank = 0
34 --------------------
35 two_descent_by_two_isogeny:
36
37 2-isogeny
38
39 changing coordinates
40
41 new curve is y^2 == x( x^2 + (0)x + (17) )
42 new isogenous curve is y^2 == x( x^2 + (0)x + (-68) )
43 Found small global point, quartic (1,0,0,0,17)
44 Found small global point, quartic (17,0,0,0,1)
45 Found small global point, quartic (1,0,0,0,-68)
46 ELS without small global points, quartic (2,0,0,0,-34)
47 ELS without small global points, quartic (17,0,0,0,-4)
48 ELS without small global points, quartic (34,0,0,0,-2)
49 ELS without small global points, quartic (-1,0,0,0,68)
50 ELS without small global points, quartic (-2,0,0,0,34)
51 Found small global point, quartic (-17,0,0,0,4)
52 ELS without small global points, quartic (-34,0,0,0,2)
53
54 Results:
55 2 <= #E(Q)/phi'(E'(Q)) <= 2
56 2 <= #E'(Q)/phi(E(Q)) <= 8
57 #Sel^(phi')(E'/Q) = 2
58 #Sel^(phi)(E/Q) = 8
59 1 <= #Sha(E'/Q)[phi'] <= 1
60 1 <= #Sha(E/Q)[phi] <= 4
61 1 <= #Sha(E/Q)[2], #Sha(E'/Q)[2] <= 4
62 0 <= rank of E(Q) = rank of E'(Q) <= 2
63 (2, 2, 2, 8)
64 --------------------
65 simon_two_descent:
66  Elliptic curve: Y^2 = x^3 + 17*x
67  E[2] = [[0], [0, 0]]
68  Elliptic curve: Y^2 = x^3 + 17*x
69  trivial points on E(Q) = [[0, 0], [1, 1, 0], [0, 0], [0, 0]]
70
71  points on E(Q) = [[0, 0]]
72
73 [E(Q):phi'(E'(Q))]   = 2
74 #S^(phi')(E'/Q)      = 2
75 #III(E'/Q)[phi']     = 1
76
77  trivial points on E'(Q) = [[0, 0], [1, 1, 0], [0, 0], [0, 0]]
78
79  points on E'(Q) = [[0, 0]]
80  points on E(Q) = [[0, 0]]
81
82 [E'(Q):phi(E(Q))]    >= 2
83 #S^(phi)(E/Q)        = 8
84 #III(E/Q)[phi]       <= 4
85
86 #III(E/Q)[2]         <= 4
87 #E(Q)[2]             = 2
88 #E(Q)/2E(Q)          >= 2
89
90 0 <= rank            <= 2
91
92 points = [[0, 0]]
93 (0, 3, [])
```

So, the Cremona-Label cannot be found by `SageMath`.

### 9.0.2  Cremona Tables

We take a look at the Cremona database available at

```
http://johncremona.github.io/ecdata/
```

and select the conductor range `10000-19999`. We get the following result:
$E$ has the Cremona-Label `18496 k1`.

- **Summary Table:**
  `https://raw.githubusercontent.com/JohnCremona/ecdata/master/count/count.10000-19999`

  ```
  18496 21
  ```

- **Table One:**
  `https://raw.githubusercontent.com/JohnCremona/ecdata/master/allcurves/allcurves.10000-19999`:

  ```
  18496 21
  ```

- **Table Two:**
  `https://raw.githubusercontent.com/JohnCremona/ecdata/master/allgens/allgens.10000-19999`:

  ```
  18496 k 1 [0,0,0,17,0] 0 [2] [0:0:1]
  ```

- **Table Three:**
  `https://raw.githubusercontent.com/JohnCremona/ecdata/master/aplist/aplist.10000-19999`:

  ```
  18496 k - 0 4 0 0 6 + 0 0 4 0 12 -8 0 0 14 0 12 0 0 -16 0 0 -10 8
  ```

# Part III
# Outlook on further Research

## 10 The Curve $E_p : Y^2 = X^3 + pX$

As we have seen in Corollary 7.6, the Selmer group is maximal for $p \equiv 1 \mod 8$ and we saw that the failure of the Hasse principle does occur in this case if 2 is not a quartic residue for $p$.

As a further research one could try to find out what happens if 2 is a quartic residue for $p$.

If we assume that the Shafarevich-Tate Group is finite, the rank could be 0 or 2. One migth try to distinguish the cases for which rank equals 0 or 2.

We have already determined the rank for all other cases $p \not\equiv 1 \mod 8$ in Corollary 7.8 under the Assumption that Ш is finite. As a further research cases $p \not\equiv 1 \mod 8$ could be analyzed without using the Shafarevich conjecture.

### 10.1 Examples

I wrote the following code in `SageMath` to determine the rank if 2 is a quartic residue for $p \equiv 1 \mod 8$:

```
pmax=input('pmax = ');
for p in [1..pmax]:
  Q = quadratic_residues(p);
  if ((is_prime(p)) & (2 in Q) & (mod(p,8)==1) ):
    E = EllipticCurve([0,0,0,p,0])
    print E
    print 'p=', p
    print 'rank=', E.rank()
    print '____________________'
#Finding quadratic residues, see also:
#http://doc.sagemath.org/html/en/constructions/number_theory.html
```

We get the following result, showing that rank 0 and 2 occur:

```
sage: load("quad-res-2.sage")
pmax = 200
Elliptic Curve defined by y^2 = x^3 + 17*x over Rational Field
p= 17
rank= 0
____________________
Elliptic Curve defined by y^2 = x^3 + 41*x over Rational Field
p= 41
rank= 0
____________________
Elliptic Curve defined by y^2 = x^3 + 73*x over Rational Field
p= 73
rank= 2
____________________
Elliptic Curve defined by y^2 = x^3 + 89*x over Rational Field
p= 89
rank= 2
____________________
Elliptic Curve defined by y^2 = x^3 + 97*x over Rational Field
p= 97
rank= 0
____________________
Elliptic Curve defined by y^2 = x^3 + 113*x over Rational Field
p= 113
rank= 2
____________________
Elliptic Curve defined by y^2 = x^3 + 137*x over Rational Field
p= 137
```

```
29 rank= 0
30 --------------------
31 Elliptic Curve defined by y^2 = x^3 + 193*x over Rational Field
32 p= 193
33 rank= 0
34 --------------------
```

So, for $p \leq 200$ the rank is 0 if $p = 17, 41, 97, 137$ or $193$ and one has rank 2 in case of $p = 73, 89$ and $113$.

# 11 Elliptic Curves over $\mathbb{Q}$ with $j$-Invariant $0$

## 11.1 The Curve $E_D : Y^2 = X^3 + D$

Let $E$ be an elliptic curve over $\mathbb{Q}$. Then, analogously to Part II, one has

$$j(E) = 0 \Leftrightarrow \exists D \in \mathbb{Q}^* \text{ s.t. } E_D : y^2 = x^3 + D \text{ is isomorphic to } E.$$

Let $E_D : y^2 = f(x) = x^3 + D$ be an elliptic curve.

Again, we identify $D$ with a sixth-power-free integer, such that $D$ is uniquely defined by $E$. According to [SIL1] Exercise 10.19, the torsion subgroup of $E$ is given by

$$E_{D,\text{tors}}(\mathbb{Q}) \cong \begin{cases} \mathbb{Z}/6\mathbb{Z} & \text{if } D = 1, \\ \mathbb{Z}/3\mathbb{Z} & \text{if } D \neq 1 \text{ is a square, or if } D = -432 \\ \mathbb{Z}/2\mathbb{Z} & \text{if } D \neq 1 \text{ is a cube} \\ 1 & \text{otherwise} \end{cases}$$

So, if $D = 1$ or if $D$ is a cube, then $E[2] \subset E_{D,\text{tors}}(\mathbb{Q})$ and one can try to execute a two-descent via two-isogeny. If $D$ is a square, or if $D = -432$ one could try a 3-desecent. As a further research, one could try to give an upper bound for the rank of $E$ and try to determine the Selmer group analogously as in the previous part. And one might analyze the question whether there are cases in which the Hasse principle fails or not.

## 11.2 Examples

### 11.2.1 The case $D = 1$

We use `SageMath` to compute the two-descent for

$$E : y^2 = x^3 + 1.$$

We use `SageMath` to compute the two-descent. Again, we execute the script `example.sage`. Setting $E$ as above, we get the following result:

```
1 sage: load("example.sage")
2 --------------------
3 Defining the elliptic curve:
4 E = EllipticCurve([0,0,0,0,1])
5 E = Elliptic Curve defined by y^2 = x^3 + 1 over Rational Field
6 discriminant = -432 = -1 * 2^4 * 3^3
7 E minimal? True
8 --------------------
9 Cremona Data:
10 36a1
11 conductor = 36 = 2^2 * 3^2
```

```
--------------------
Local Data:
[Local data at Principal ideal (2) of Integer Ring:
Reduction type: bad additive
Local minimal model: Elliptic Curve defined by y^2 = x^3 + 1 over Rational
    Field
Minimal discriminant valuation: 4
Conductor exponent: 2
Kodaira Symbol: IV
Tamagawa Number: 3, Local data at Principal ideal (3) of Integer Ring:
Reduction type: bad additive
Local minimal model: Elliptic Curve defined by y^2 = x^3 + 1 over Rational
    Field
Minimal discriminant valuation: 3
Conductor exponent: 2
Kodaira Symbol: III
Tamagawa Number: 2]
--------------------
Torsion:
Torsion Subgroup isomorphic to Z/6 associated to the Elliptic Curve defined
    by y^2 = x^3 + 1 over Rational Field
Torsion Points: [(0 : 1 : 0), (2 : 3 : 1), (0 : 1 : 1), (-1 : 0 : 1), (0 : -1
    : 1), (2 : -3 : 1)]
--------------------
Rank:
rank = 0
--------------------
two_descent_by_two_isogeny:

2-isogeny

changing coordinates

new curve is y^2 == x( x^2 + (-3)x + (3) )
new isogenous curve is y^2 == x( x^2 + (6)x + (-3) )
Found small global point, quartic (1,0,-3,0,3)
Found small global point, quartic (3,0,-3,0,1)
Found small global point, quartic (1,0,6,0,-3)
Found small global point, quartic (-3,0,6,0,1)

Results:
2 <= #E(Q)/phi'(E'(Q)) <= 2
2 <= #E'(Q)/phi(E(Q)) <= 2
#Sel^(phi')(E'/Q) = 2
#Sel^(phi)(E/Q) = 2
1 <= #Sha(E'/Q)[phi'] <= 1
1 <= #Sha(E/Q)[phi] <= 1
1 <= #Sha(E/Q)[2], #Sha(E'/Q)[2] <= 1
0 <= rank of E(Q) = rank of E'(Q) <= 0
(2, 2, 2, 2)
--------------------
simon_two_descent:
 Elliptic curve: Y^2 = x^3 + 1
 E[2] = [[0], [-1, 0]]
 Elliptic curve: Y^2 = x^3 - 3*x^2 + 3*x
 trivial points on E(Q) = [[0, 0], [1, 1, 0], [0, 0], [0, 0], [1, 1], [3, 3]]

 points on E(Q) = [[3, 3]]

[E(Q):phi'(E'(Q))]  = 2
#S^(phi')(E'/Q)     = 2
#III(E'/Q)[phi']    = 1

 trivial points on E'(Q) = [[0, 0], [1, 1, 0], [0, 0], [0, 0], [1, 2], [-3,
    6]]

 points on E'(Q) = [[-3, 6]]
 points on E(Q) = [[1, -1]]

[E'(Q):phi(E(Q))]   = 2
#S^(phi)(E/Q)       = 2
#III(E/Q)[phi]      = 1

#III(E/Q)[2]        = 1
```

```
81 #E(Q)[2]            = 2
82 #E(Q)/2E(Q)         = 2
83 rank                = 0
84
85 points = [[0, 0]]
86 Rank determined successfully, saturating...
87 (0, 1, [])
```

### 11.2.2 The case $D \neq 1$ cubic

Let $D \neq 1$ be cubic, and let

$$E_p : y^2 = x^3 + D$$

be the Weierstrass equation for the elliptic curve $E_D$. Using the notations from above let $\phi : E_D \to E'_D$ be an isogeny of degree 2 with kernel $E_p[\phi] = \{O, (0,0)\}$.

$D = 8$:

```
 1 sage: load("example.sage")
 2 --------------------
 3 Defining the elliptic curve:
 4 E = EllipticCurve([0,0,0,0,8])
 5 E = Elliptic Curve defined by y^2 = x^3 + 8 over Rational Field
 6 discriminant = -27648 = -1 * 2^10 * 3^3
 7 E minimal? True
 8 --------------------
 9 Cremona Data:
10 576a1
11 conductor = 576 = 2^6 * 3^2
12 --------------------
13 Local Data:
14 [Local data at Principal ideal (2) of Integer Ring:
15 Reduction type: bad additive
16 Local minimal model: Elliptic Curve defined by y^2 = x^3 + 8 over Rational
    Field
17 Minimal discriminant valuation: 10
18 Conductor exponent: 6
19 Kodaira Symbol: I0*
20 Tamagawa Number: 2, Local data at Principal ideal (3) of Integer Ring:
21 Reduction type: bad additive
22 Local minimal model: Elliptic Curve defined by y^2 = x^3 + 8 over Rational
    Field
23 Minimal discriminant valuation: 3
24 Conductor exponent: 2
25 Kodaira Symbol: III
26 Tamagawa Number: 2]
27 --------------------
28 Torsion:
29 Torsion Subgroup isomorphic to Z/2 associated to the Elliptic Curve defined
    by y^2 = x^3 + 8 over Rational Field
30 Torsion Points: [(0 : 1 : 0), (-2 : 0 : 1)]
31 --------------------
32 Rank:
33 rank = 1
34 --------------------
35 two_descent_by_two_isogeny:
36
37 2-isogeny
38
39 changing coordinates
40
41 new curve is y^2 == x( x^2 + (-6)x + (12) )
42 new isogenous curve is y^2 == x( x^2 + (12)x + (-12) )
43 Found small global point, quartic (1,0,-6,0,12)
44 Found small global point, quartic (3,0,-6,0,4)
45 Found small global point, quartic (1,0,12,0,-12)
46 Found small global point, quartic (6,0,12,0,-2)
47 Found small global point, quartic (-2,0,12,0,6)
48 Found small global point, quartic (-3,0,12,0,4)
49
```

```
50 Results:
51  2 <= #E(Q)/phi'(E'(Q)) <= 2
52  4 <= #E'(Q)/phi(E(Q)) <= 4
53  #Sel^(phi')(E'/Q) = 2
54  #Sel^(phi)(E/Q) = 4
55  1 <= #Sha(E'/Q)[phi'] <= 1
56  1 <= #Sha(E/Q)[phi] <= 1
57  1 <= #Sha(E/Q)[2], #Sha(E'/Q)[2] <= 1
58  1 <= rank of E(Q) = rank of E'(Q) <= 1
59  (2, 2, 4, 4)
60  --------------------
61 simon_two_descent:
62   Elliptic curve: Y^2 = x^3 + 8
63   E[2] = [[0], [-2, 0]]
64   Elliptic curve: Y^2 = x^3 - 6*x^2 + 12*x
65   trivial points on E(Q) = [[0, 0], [1, 1, 0], [0, 0], [0, 0], [1/4, 13/8],
       [3, 3]]
66
67   points on E(Q) = [[3, 3]]
68
69  [E(Q):phi'(E'(Q))]  = 2
70  #S^(phi')(E'/Q)     = 2
71  #III(E'/Q)[phi']    = 1
72
73   trivial points on E'(Q) = [[0, 0], [1, 1, 0], [0, 0], [0, 0], [1, 1], [-2,
       8]]
74
75   points on E'(Q) = [[-2, 8], [0, 0]]
76   points on E(Q) = [[4, -4], [0, 0]]
77
78  [E'(Q):phi(E(Q))]   = 4
79  #S^(phi)(E/Q)       = 4
80  #III(E/Q)[phi]      = 1
81
82  #III(E/Q)[2]        = 1
83  #E(Q)[2]            = 2
84  #E(Q)/2E(Q)         = 4
85  rank                = 1
86
87 points = [[0, 0], [3, 3]]
88 Rank determined successfully, saturating...
89 (1, 2, [(1 : 3 : 1)])
```

$D = 27$:

```
1 sage: load("example.sage")
2 --------------------
3 Defining the elliptic curve:
4 E = EllipticCurve([0,0,0,0,27])
5 E = Elliptic Curve defined by y^2 = x^3 + 27 over Rational Field
6 discriminant = -314928 = -1 * 2^4 * 3^9
7 E minimal? True
8 --------------------
9 Cremona Data:
10 144a3
11 conductor = 144 = 2^4 * 3^2
12 --------------------
13 Local Data:
14 [Local data at Principal ideal (2) of Integer Ring:
15 Reduction type: bad additive
16 Local minimal model: Elliptic Curve defined by y^2 = x^3 + 27 over Rational
    Field
17 Minimal discriminant valuation: 4
18 Conductor exponent: 4
19 Kodaira Symbol: II
20 Tamagawa Number: 1, Local data at Principal ideal (3) of Integer Ring:
21 Reduction type: bad additive
22 Local minimal model: Elliptic Curve defined by y^2 = x^3 + 27 over Rational
    Field
23 Minimal discriminant valuation: 9
24 Conductor exponent: 2
25 Kodaira Symbol: III*
26 Tamagawa Number: 2]
27 --------------------
```

```
28 Torsion:
29 Torsion Subgroup isomorphic to Z/2 associated to the Elliptic Curve defined
     by y^2 = x^3 + 27 over Rational Field
30 Torsion Points: [(0 : 1 : 0), (-3 : 0 : 1)]
31 --------------------
32 Rank:
33 rank = 0
34 --------------------
35 two_descent_by_two_isogeny:
36
37 2-isogeny
38
39 changing coordinates
40
41 new curve is y^2 == x( x^2 + (-9)x + (27) )
42 new isogenous curve is y^2 == x( x^2 + (18)x + (-27) )
43 Found small global point, quartic (1,0,-9,0,27)
44 Found small global point, quartic (3,0,-9,0,9)
45 Found small global point, quartic (1,0,18,0,-27)
46 Found small global point, quartic (-3,0,18,0,9)
47
48 Results:
49 2 <= #E(Q)/phi'(E'(Q)) <= 2
50 2 <= #E'(Q)/phi(E(Q)) <= 2
51 #Sel^(phi')(E'/Q) = 2
52 #Sel^(phi)(E/Q) = 2
53 1 <= #Sha(E'/Q)[phi'] <= 1
54 1 <= #Sha(E/Q)[phi] <= 1
55 1 <= #Sha(E/Q)[2], #Sha(E'/Q)[2] <= 1
56 0 <= rank of E(Q) = rank of E'(Q) <= 0
57 (2, 2, 2, 2)
58 --------------------
59 simon_two_descent:
60  Elliptic curve: Y^2 = x^3 + 27
61  E[2] = [[0], [-3, 0]]
62  Elliptic curve: Y^2 = x^3 - 9*x^2 + 27*x
63  trivial points on E(Q) = [[0, 0], [1, 1, 0], [0, 0], [0, 0]]
64
65  points on E(Q) = [[0, 0]]
66
67 [E(Q):phi'(E'(Q))]  = 2
68 #S^(phi')(E'/Q)     = 2
69 #III(E'/Q)[phi']    = 1
70
71  trivial points on E'(Q) = [[0, 0], [1, 1, 0], [0, 0], [0, 0]]
72
73  points on E'(Q) = [[0, 0]]
74  points on E(Q) = [[0, 0]]
75
76 [E'(Q):phi(E(Q))]   = 2
77 #S^(phi)(E/Q)       = 2
78 #III(E/Q)[phi]      = 1
79
80 #III(E/Q)[2]        = 1
81 #E(Q)[2]            = 2
82 #E(Q)/2E(Q)         = 2
83 rank                = 0
84
85 points = [[0, 0]]
86 Rank determined successfully, saturating...
87 (0, 1, [])
```

$D = 64$:

```
1 sage: load("example.sage")
2 --------------------
3 Defining the elliptic curve:
4 E = EllipticCurve([0,0,0,0,64])
5 E = Elliptic Curve defined by y^2 = x^3 + 64 over Rational Field
6 discriminant = -1769472 = -1 * 2^16 * 3^3
7 E minimal? False
8 --------------------
9 Cremona Data:
10 36a1
```

```
11 conductor = 36 = 2^2 * 3^2
12 --------------------
13 Local Data:
14 [Local data at Principal ideal (2) of Integer Ring:
15 Reduction type: bad additive
16 Local minimal model: Elliptic Curve defined by y^2 = x^3 + 1 over Rational
       Field
17 Minimal discriminant valuation: 4
18 Conductor exponent: 2
19 Kodaira Symbol: IV
20 Tamagawa Number: 3, Local data at Principal ideal (3) of Integer Ring:
21 Reduction type: bad additive
22 Local minimal model: Elliptic Curve defined by y^2 = x^3 + 1 over Rational
       Field
23 Minimal discriminant valuation: 3
24 Conductor exponent: 2
25 Kodaira Symbol: III
26 Tamagawa Number: 2]
27 --------------------
28 Torsion:
29 Torsion Subgroup isomorphic to Z/6 associated to the Elliptic Curve defined
      by y^2 = x^3 + 64 over Rational Field
30 Torsion Points: [(0 : 1 : 0), (8 : 24 : 1), (0 : 8 : 1), (-4 : 0 : 1), (0 :
      -8 : 1), (8 : -24 : 1)]
31 --------------------
32 Rank:
33 rank = 0
34 --------------------
35 two_descent_by_two_isogeny:
36
37 2-isogeny
38
39 changing coordinates
40
41 new curve is y^2 == x( x^2 + (-12)x + (48) )
42 new isogenous curve is y^2 == x( x^2 + (24)x + (-48) )
43 Found small global point, quartic (1,0,-12,0,48)
44 Found small global point, quartic (3,0,-12,0,16)
45 Found small global point, quartic (1,0,24,0,-48)
46 Found small global point, quartic (-3,0,24,0,16)
47
48 Results:
49 2 <= #E(Q)/phi'(E'(Q)) <= 2
50 2 <= #E'(Q)/phi(E(Q)) <= 2
51 #Sel^(phi')(E'/Q) = 2
52 #Sel^(phi)(E/Q) = 2
53 1 <= #Sha(E'/Q)[phi'] <= 1
54 1 <= #Sha(E/Q)[phi] <= 1
55 1 <= #Sha(E/Q)[2], #Sha(E'/Q)[2] <= 1
56 0 <= rank of E(Q) = rank of E'(Q) <= 0
57 (2, 2, 2, 2)
58 --------------------
59 simon_two_descent:
60  Elliptic curve: Y^2 = x^3 + 64
61  E[2] = [[0], [-4, 0]]
62  Elliptic curve: Y^2 = x^3 - 12*x^2 + 48*x
63  trivial points on E(Q) = [[0, 0], [1, 1, 0], [0, 0], [0, 0]]
64
65  points on E(Q) = [[0, 0]]
66
67 [E(Q):phi'(E'(Q))]   = 2
68 #S^(phi')(E'/Q)      = 2
69 #III(E'/Q)[phi']     = 1
70
71  trivial points on E'(Q) = [[0, 0], [1, 1, 0], [0, 0], [0, 0]]
72
73  points on E'(Q) = [[0, 0]]
74  points on E(Q) = [[0, 0]]
75
76 [E'(Q):phi(E(Q))]    = 2
77 #S^(phi)(E/Q)        = 2
78 #III(E/Q)[phi]       = 1
79
80 #III(E/Q)[2]         = 1
```

```
81 #E(Q)[2]            = 2
82 #E(Q)/2E(Q)         = 2
83 rank                = 0
84
85 points = [[0, 0]]
86 Rank determined successfully, saturating...
87 (0, 1, [])
```

$D = 125$:

```
 1 sage: load("example.sage")
 2 --------------------
 3 Defining the elliptic curve:
 4 E = EllipticCurve([0,0,0,0,125])
 5 E = Elliptic Curve defined by y^2 = x^3 + 125 over Rational Field
 6 discriminant = -6750000 = -1 * 2^4 * 3^3 * 5^6
 7 E minimal? True
 8 --------------------
 9 Cremona Data:
10 900b1
11 conductor = 900 = 2^2 * 3^2 * 5^2
12 --------------------
13 Local Data:
14 [Local data at Principal ideal (2) of Integer Ring:
15 Reduction type: bad additive
16 Local minimal model: Elliptic Curve defined by y^2 = x^3 + 125 over Rational
    Field
17 Minimal discriminant valuation: 4
18 Conductor exponent: 2
19 Kodaira Symbol: IV
20 Tamagawa Number: 1, Local data at Principal ideal (3) of Integer Ring:
21 Reduction type: bad additive
22 Local minimal model: Elliptic Curve defined by y^2 = x^3 + 125 over Rational
    Field
23 Minimal discriminant valuation: 3
24 Conductor exponent: 2
25 Kodaira Symbol: III
26 Tamagawa Number: 2, Local data at Principal ideal (5) of Integer Ring:
27 Reduction type: bad additive
28 Local minimal model: Elliptic Curve defined by y^2 = x^3 + 125 over Rational
    Field
29 Minimal discriminant valuation: 6
30 Conductor exponent: 2
31 Kodaira Symbol: I0*
32 Tamagawa Number: 2]
33 --------------------
34 Torsion:
35 Torsion Subgroup isomorphic to Z/2 associated to the Elliptic Curve defined
    by y^2 = x^3 + 125 over Rational Field
36 Torsion Points: [(0 : 1 : 0), (-5 : 0 : 1)]
37 --------------------
38 Rank:
39 rank = 0
40 --------------------
41 two_descent_by_two_isogeny:
42
43 2-isogeny
44
45 changing coordinates
46
47 new curve is y^2 == x( x^2 + (-15)x + (75) )
48 new isogenous curve is y^2 == x( x^2 + (30)x + (-75) )
49 Found small global point, quartic (1,0,-15,0,75)
50 Found small global point, quartic (3,0,-15,0,25)
51 Found small global point, quartic (1,0,30,0,-75)
52 Found small global point, quartic (-3,0,30,0,25)
53
54 Results:
55 2 <= #E(Q)/phi'(E'(Q)) <= 2
56 2 <= #E'(Q)/phi(E(Q)) <= 2
57 #Sel^(phi')(E'/Q) = 2
58 #Sel^(phi)(E/Q) = 2
59 1 <= #Sha(E'/Q)[phi'] <= 1
60 1 <= #Sha(E/Q)[phi] <= 1
```

```
1 <= #Sha(E/Q)[2], #Sha(E'/Q)[2] <= 1
0 <= rank of E(Q) = rank of E'(Q) <= 0
(2, 2, 2, 2)
--------------------
simon_two_descent:
 Elliptic curve: Y^2 = x^3 + 125
 E[2] = [[0], [-5, 0]]
 Elliptic curve: Y^2 = x^3 - 15*x^2 + 75*x
 trivial points on E(Q) = [[0, 0], [1, 1, 0], [0, 0], [0, 0]]

 points on E(Q) = [[0, 0]]

[E(Q):phi'(E'(Q))]   = 2
#S^(phi')(E'/Q)      = 2
#III(E'/Q)[phi']     = 1

 trivial points on E'(Q) = [[0, 0], [1, 1, 0], [0, 0], [0, 0]]

 points on E'(Q) = [[0, 0]]
 points on E(Q) = [[0, 0]]

[E'(Q):phi(E(Q))]    = 2
#S^(phi)(E/Q)        = 2
#III(E/Q)[phi]       = 1

#III(E/Q)[2]         = 1
#E(Q)[2]             = 2
#E(Q)/2E(Q)          = 2
rank                 = 0

points = [[0, 0]]
Rank determined successfully, saturating...
(0, 1, [])
```

# Corrections to the Master Thesis
## *Using Selmer Groups*
## *to compute Mordell-Weil Groups*
## *of Elliptic Curves*

Anika Behrens

October 23, 2016

## Section 10.1

My `SageMath` script "`quad-res-2.sage`" given in Section 10.1 computes quadratic residues, but we need quartic residues. The following script "`quartic-res-2.sage`" uses the theorem of Gauss (Quartic Reciprocity, Proposition 8.3):

```
pmax=input('pmax = ');
for p in [1..pmax]:
  Q = quadratic_residues(p);
  if ((is_prime(p)) & (2 in Q) & (mod(p,8)==1) ):
    squares = two_squares(p)
    if(mod(squares[0]*squares[1],8)==0):
      E = EllipticCurve([0,0,0,p,0])
      print E
      print 'p=', p
      print 'rank=', E.rank(only_use_mwrank=False)
      print '___________________'
```

We get the following result, showing that rank 0 and 2 occur.

```
sage: load("quartic-res-2.sage")
pmax = 600
Elliptic Curve defined by y^2 = x^3 + 73*x over Rational Field
p= 73
rank= 2
___________________
Elliptic Curve defined by y^2 = x^3 + 89*x over Rational Field
p= 89
rank= 2
___________________
Elliptic Curve defined by y^2 = x^3 + 113*x over Rational Field
p= 113
rank= 2
___________________
Elliptic Curve defined by y^2 = x^3 + 233*x over Rational Field
p= 233
rank= 2
___________________
Elliptic Curve defined by y^2 = x^3 + 257*x over Rational Field
p= 257
rank= 0
___________________
Elliptic Curve defined by y^2 = x^3 + 281*x over Rational Field
p= 281
rank= 2
___________________
Elliptic Curve defined by y^2 = x^3 + 337*x over Rational Field
```

```
28  p= 337
29  rank= 2
30  --------------------
31  Elliptic Curve defined by y^2 = x^3 + 353*x over Rational Field
32  p= 353
33  rank= 2
34  --------------------
35  Elliptic Curve defined by y^2 = x^3 + 577*x over Rational Field
36  p= 577
37  rank= 0
38  --------------------
39  Elliptic Curve defined by y^2 = x^3 + 593*x over Rational Field
40  p= 593
41  rank= 2
42  --------------------
```

So, for $p \leq 600$ the rank is 2 if $p = 73, 89, 113, 233, 281, 337, 353$ or $593$ and the rank is 0 if $p = 257$ or $577$.

For $p \leq 10000$ we get the following results:

**Rank 0:**   257 577 1097 1201 1217 1481 1721 2441 2657 2833 2857 3121 3449 3761 4001 4057 4177 4217 4297 4409 4481 4657 4721 4817 4937 5297 5569 5737 6121 6481 6521 6793 6841 6857 7121 7129 7793 7817 7841 8081 8161 8761 9001 9137 9209 9241 9281 9697 9769

**Rank 2:**   73 89 113 233 281 337 353 593 601 617 881 937 1033 1049 1153 1193 1249 1289 1433 1553 1601 1609 1753 1777 1801 1889 1913 2089 2113 2129 2273 2281 2393 2473 2593 2689 2969 3049 3089 3137 3217 3257 3313 3361 3529 3673 3833 4049 4153 4201 4273 4289 4457 4513 4801 4993 5081 5113 5209 5233 5393 5689 5881 6089 6353 6361 6449 6529 6553 6569 6689 6761 7393 7481 7489 7529 7577 7753 7993 8209 8233 8273 8369 8537 8609 8713 8969 9337 9377 9473 9521 9601 9649 9721